\documentstyle[11pt,psfig,amssymb,amsmath]{article}
\newcommand{\be}{\begin{equation}} 
\newcommand{\ee}{\end{equation}}   
\newcommand{\bea}{\begin{eqnarray}}
\newcommand{\eea}{\end{eqnarray}}

\newtheorem{theorem}{Theorem}

\newtheorem{corollary}{Corollary}
\newtheorem{lemma}{Lemma}
\newtheorem{example}{Example}

\newtheorem{definition}{Definition}   

\def\1#1{^{(#1)}}

\begin{document}
\title{Frobenius Problem for Semigroups ${\sl S}\left(d_1,d_2,d_3\right)$}

\author{Leonid G. Fel\\
\\Department of Civil and Environmental Engineering,\\
Technion, Haifa 3200, Israel\\
\vspace{-.3cm}
\\{\sl e-mail: lfel@techunix.technion.ac.il}}
\date{\today}
\maketitle
\def\be{\begin{equation}}
\def\ee{\end{equation}}
\def\p{\prime}  
\vspace{-1cm}
\begin{abstract}
The matrix representation of the set $\Delta({\bf d}^3)$, ${\bf d}^3=(d_1,d_2,
d_3)$, of the integers which are unrepresentable by $d_1,d_2,d_3$ is found. The 
diagrammatic procedure of calculation of the generating function $\Phi\left(
{\bf d}^3;z\right)$ for the set $\Delta({\bf d}^3)$ is developed. The Frobenius 
number $F\left({\bf d}^3\right)$, genus $G\left({\bf d}^3\right)$ and Hilbert 
series $H({\bf d}^3;z)$ of a graded subring for non--symmetric and symmetric 
semigroups ${\sf S}\left({\bf d}^3\right)$ are found. The upper bound for the 
number of non--zero coefficients in the polynomial numerators of Hilbert series 
$H({\bf d}^m;z)$ of graded subrings for non--symmetric semigroups ${\sf S}
\left({\bf d}^m\right)$ of dimension, $m\geq 4$, is established.
\vskip .1 cm
\noindent
\begin{tabbing}
{\sf Key words}:\hspace{.1in}  \=  Restricted partitions, Frobenius problem, 
Non--symmetric and symmetric semigroups, \\
               \> Hilbert series of a graded subring.
\end{tabbing}
\end{abstract}

\vspace{.6cm}
\newpage
\tableofcontents

\newpage   
\section{Introduction}
\label{intro}

Let ${\sf S}\left(d_1,\ldots,d_m\right)\subset {\mathbb N}$ be 
the subsemigroup 
generated by a set of integers $\{d_1,\ldots,d_m\}$ such that 
\begin{equation}
1<d_1<\ldots <d_m\;,\;\;\;\gcd(d_1,\ldots,d_m)=1\;.
\label{definn0}
\end{equation}
The set $\{d_1,\ldots,d_m\}$ is called {\em minimal} if there are no  
nonnegative integers $b_{i,j}$ for which the following linear dependence 
holds
\begin{equation}
d_i=\sum_{j\neq i}^mb_{i,j}d_j\;,\;\;\;
b_{i,j}\in\{0,1,\ldots\}\;\;\;\mbox{for any}\;\;i\leq m\;.
\label{defin2}
\end{equation}
For short we denote the tuple $(d_1,\ldots,d_m)$ by ${\bf d}^m$ where $m$ is 
{\em the dimension} of ${\bf d}^m$. Henceforth ${\bf d}^m$ will be a minimal 
generating set of ${\sf S}\left({\bf d}^m\right)$. {\em The conductor} $c\left
({\bf d}^m\right)$ of ${\sf S}\left({\bf d}^m\right)$ is defined by $c\left(
{\bf d}^m\right):=\min\left\{s\in {\sf S}\left({\bf d}^m\right)\;|\;s+{\mathbb 
N}\cup \{0\}\subset {\sf S}\left({\bf d}^m\right)\right\}$. {\em The genus} 
$G\left({\bf d}^m\right)$ of ${\sf S}\left({\bf d}^m\right)$ is defined as the 
cardinality ($\#$) of its complement $\Delta$ in ${\mathbb N}$, i.e. 
$\Delta\left({\bf d}^m\right)={\mathbb N}\setminus 
{\sf S}\left({\bf d}^m\right)$ and 
\begin{eqnarray}
G\left({\bf d}^m\right):=\#\Delta\left({\bf d}^m\right)\;. 
\label{defgen}  
\end{eqnarray}
For the reason explained in Section \ref{nomina0}, it is worth to introduce 
the generating function $\Phi\left({\bf d}^m;z\right)$ for the set $\Delta(
{\bf d}^m)$ of unrepresentable integers in accordance with \cite{stan96} 
\begin{eqnarray}
\Phi\left({\bf d}^m;z\right)=\sum_{s\;\in\;\Delta\left({\bf d}^m\right)} z^s\;,
\label{def1}
\end{eqnarray}
Let ${\bf d}^m$ and ${\bf d}^n$ be two tuples of different dimensions, $m>n$. 
Define a relation ${\bf d}^n\prec {\bf d}^m$ if ${\bf d}^n$ is an initial 
segment of ${\bf d}^m$ as a {\em word} in $d_1,\ldots,d_m$.  If for any $2<n<m$ 
the tuple ${\bf d}^n$ is $(d_1,\ldots,d_n)$ then ${\bf d}^2\prec \ldots\prec 
{\bf d}^n\prec \ldots\prec {\bf d}^m$. This implies an ordering in the 
following three sequences
\begin{eqnarray}
\Delta({\bf d}^{m}) \subset \ldots\subset \Delta({\bf d}^2)\;,\;\;\;\;
G({\bf d}^{m})<\ldots <G({\bf d}^2)\;,\;\;\;\;
c({\bf d}^{m})\leq \ldots \leq c({\bf d}^2)\;.
\label{def1b}
\end{eqnarray}
The semigroup ring ${\sf k}\left[X_1,\ldots,X_m\right]$ over a field ${\sf k}$ 
of characteristic 0 associated with ${\sf S}\left({\bf d}^m\right)$ is a 
polynomial subring graded by $\deg X_i=d_i,1\leq i\leq m$ and generated by all 
monomials $z^{d_i}$. The Hilbert series $H({\bf d}^m;z)$ of a graded subring 
${\sf k}\left[z^{d_1},\ldots,z^{d_m}\right]$ is defined by \cite{stan96}
\begin{equation}
H({\bf d}^m;z)=\sum_{s\;\in\;{\sf S}\left({\bf d}^m\right)} z^s=
\frac{Q({\bf d}^m;z)}{\prod_{j=1}^m \left(1-z^{d_j}\right)}\;,
\label{hilb0}
\end{equation}
where $Q({\bf d}^m;z)$ is a polynomial in $z$. 
The number 
\begin{equation}
F\left({\bf d}^m\right):=-1+c\left({\bf d}^m\right)
\label{condfrob}
\end{equation} 
is referred to as {\em Frobenius number} in honor of G. Frobenius who, 
according to \cite{brau42}, repeatedly raised the following question in 
his lectures: 
determine (or bound) $F\left({\bf d}^m\right)$. Actually, all three entities, 
$F\left({\bf d}^m\right)$, $G\left({\bf d}^m\right)$ and $H({\bf d}^m;z)$, are 
originated by the same semigroup ${\sf S}\left({\bf d}^m\right)$ and have a 
strong algebraic relationship (see Section \ref{nomina0}). Due to this reason 
the determination of $F\left({\bf d}^m\right)$, $G\left({\bf d}^m\right)$ and 
$H({\bf d}^m;z)$ will be called the $m$--dimensional (mD) {\em Frobenius 
problem}.

Let $R={\sf k}\left[X_1,\ldots,X_m\right]$ be the ring of polynomials over a 
field ${\sf k}$ and $\pi:R\longmapsto {\sf k}\left[{\sf S}\left({\bf d}^m\right)
\right]$ be the projection induced by $\pi\left(X_i\right)=z^{d_i}$. Then 
${\sf k}\left[{\sf S}\left({\bf d}^m\right)\right]$ has a presentation ${\sf k}
\left[X_1,\ldots,X_m\right]/{\cal I}_m$ where ${\cal I}_m$ is the kernel of 
the map $\pi$. The semigroup ${\sf S}\left({\bf d}^m\right)$ is called {\em 
symmetric} iff for all $s\in {\sf S}\left({\bf d}^m\right)$ the following holds 
$F({\bf 
d}^m)-s\not\in {\sf S}\left({\bf d}^m\right)$. This kind of semigroups is of 
high importance due to Kunz's theorem \cite{kunz70} which asserts that ${\sf k}
\left[{\sf S}\left({\bf d}^m\right)\right]$ is a Gorenstein ring iff ${\sf S}
\left({\bf d}^m\right)$ is symmetric. It is classically known that in small 
dimensions $m=2,3$ the situation is even simpler. For every ${\bf d}^2$, 
${\sf S}\left({\bf d}^2\right)$ is a symmetric semigroup \cite{brau42} and 
${\sf k}\left[{\sf S}\left({\bf d}^2\right)\right]$ is a complete intersection 
\cite{mora87}. The kernel  ${\cal I}_2$ is principal and has the generator 
$p=X_1^{c_1}-X_2^{c_2}$ where $c_i={\sf lcm}(d_1,d_2)/d_i$. For $m=3$, Herzog 
\cite{herz70} has proved that ${\sf k}\left[{\sf S}\left({\bf d}^3\right)\right
]$ is a complete intersection iff ${\sf S}\left({\bf d}^3\right)$ is symmetric. 

For larger $m$ the generic semigroup ${\sf S}\left({\bf d}^m\right)$ is mostly 
non--symmetric, e.g. ${\sf S}\left({\bf d}^3\right)$ minimally generated by 
three pairwise relatively prime elements $d_i$ is such a semigroup \cite
{frob87}. Concerning the Frobenius numbers, a theorem of Curtis \cite{curt90} 
asserts that, for $m=3$, there is no non--zero {\em polynomial} $P\in {\mathbb 
C}\left(Y_1,Y_2,Y_3,Z\right)$ such that $P\left(d_1,d_2,d_3,F\left({\bf d}^3
\right)\right)=0$ for {\em all minimal sets} $\left(d_1,d_2,d_3\right)$ where 
$d_1,d_2$ are primes not dividing $d_3$. In other words, $F\left({\bf d}^3
\right)$ cannot be determined for {\em all minimal sets} $\left(d_1,d_2,d_3
\right)$ by any set of closed formulas which could be reduced to a finite set 
of {\em polynomials} 
\footnote{The words "{\em all}" and "{\em polynomial}" are essential here, 
since there exist {\em infinitely} many triples $d_1,d_2,d_3$ of primes in 
arithmetic progression \cite{corp39} constituting a minimal set and, according 
to Roberts \cite{rober56}, the Frobenius number associated to them can be 
presented in closed, but {\em not polynomial} formula $F(d,d+p,d+2p)=d\left
\lfloor\frac{d-2}{2}\right\rfloor+(d-1)p,\;\gcd(d,p)=1$. A standard notation 
$\left\lfloor a\right\rfloor$ is used for the integer part of a real number $a$.
}. As for Hilbert series, for any $m\geq 4$, there is no way to write the 
rational function $H({\bf d}^m;z)$ so that its polynomial numerator $Q({\bf 
d}^m;z)$ has a bounded number of non--zero terms for all choices of $d_1,
\ldots,d_m$ \cite{worm86}. The semigroup ${\sf S}\left({\bf d}^3\right)$ 
presents the first nontrivial and most elaborated case.

Our main results are the expressions for the Frobenius number $F\left({\bf d}^3
\right)$, genus $G\left({\bf d}^3\right)$ and the numerator $Q({\bf d}^3;z)$ of 
Hilbert series for both symmetric and non--symmetric semigroups ${\sf S}\left(
{\bf d}^3\right)$. In order to present them introduce auxiliary notions. 
Following Johnson \cite{john60} define {\em the 1st minimal relation} 
${\cal R}_1\left({\bf d}^3\right)$ for given ${\bf d}^3=(d_1,d_2,d_3)$ as 
follows
\begin{eqnarray}
{\cal R}_1\left({\bf d}^3\right):\;\;
a_{11}d_1=a_{12}d_2+a_{13}d_3\;,\;\;\;a_{22}d_2=a_{21}d_1+a_{23}d_3\;,\;\;\;
a_{33}d_3=a_{31}d_1+a_{32}d_2\;,
\label{herznon1a}
\end{eqnarray}
where
\begin{eqnarray}
a_{11}&=&\min\left\{v_{11}\;\bracevert\;v_{11}\geq 2,\;v_{11}d_1=v_{12}d_2+
v_{13}d_3,\;v_{12},v_{13}\in {\mathbb N}\cup\{0\}\right\}\;,\nonumber\\
a_{22}&=&\min\left\{v_{22}\;\bracevert\;v_{22}\geq 2,\;v_{22}d_2=v_{21}d_1+
v_{23}d_3,\;v_{21},v_{23}\in {\mathbb N}\cup\{0\}\right\}\;,\label{herznon3}\\
a_{33}&=&\min\left\{v_{33}\;\bracevert\;v_{33}\geq 2,\;v_{33}d_3=v_{31}d_1+
v_{32}d_2,\;v_{31},v_{32}\in {\mathbb N}\cup\{0\}\right\}\;.
\nonumber
\end{eqnarray}
The uniquely defined values of $v_{ij},i\neq j$ which give $a_{ii}$ will be 
denoted by $a_{ij},i\neq j$. Note that due to minimality of the set $(d_1,
d_2,d_3)$ the elements $a_{ij},i,j\geq 3$ satisfy 
\begin{eqnarray}
\gcd(a_{11},a_{12},a_{13})=1\;,\;\;\;
\gcd(a_{21},a_{22},a_{23})=1\;,\;\;\;
\gcd(a_{31},a_{32},a_{33})=1\;.\label{herznon3a}
\end{eqnarray}
The procedure defined in (\ref{herznon3}) completely determines the elements 
$a_{ij},i,j\geq 3$ of ${\cal R}_1\left({\bf d}^3\right)$ as the functions 
$a_{ij}=a_{ij}(d_1,d_2,d_3)$. 

For $F\left({\bf d}^3\right)$ , $G\left({\bf d}^3\right)$ and $Q({\bf d}^3;z)$ 
we get the following formulas
\begin{eqnarray}
F\left({\bf d}^3\right)&=&\frac{1}{2}\left[\langle{\bf a},{\bf d}\rangle+
J\left({\bf d}^3\right)\right]-\sum_{i=1}^3d_i\;,\;\;\;\;
\langle{\bf a},{\bf d}\rangle=\sum_{i=1}^3a_{ii}d_i\;,\label{kra00}\\
G\left({\bf d}^3\right)&=&\frac{1}{2}\left(1+\langle{\bf a},{\bf d}\rangle
-\prod_{i=1}^3a_{ii}-\sum_{i=1}^3d_i\right)\;,\label{herznon5}\\
Q({\bf d}^3;z)&=&1-\sum_{i=1}^3 z^{a_{ii}d_i}+
z^{1/2\left[\langle{\bf a},{\bf d}\rangle-J\left({\bf d}^3\right)\right]}+  
z^{1/2\left[\langle{\bf a},{\bf d}\rangle+J\left({\bf d}^3\right)\right]}\;,
\label{nomi4}
\end{eqnarray}
where $J\left({\bf d}^3\right)$ is the positive integer
\begin{eqnarray}
J\left({\bf d}^3\right)=\sqrt{\langle{\bf a},{\bf d}\rangle^2-
4\sum_{i>j}^3a_{ii}a_{jj}d_id_j+4d_1d_2d_3}\;.
\label{kra02}
\end{eqnarray}
Formula (\ref{nomi4}) for the numerator $Q({\bf d}^3;z)$ of Hilbert series of
non--symmetric semigroup ${\sf S}\left({\bf d}^3\right)$ is new and was not 
obtained earlier. Formula (\ref{herznon5}) is in full agreement with the genus 
of generic monomial space curves found by Kraft \cite{kraf85} while formula 
(\ref{kra00}) can be reduced to the known expressions for symmetric and 
non--symmetric semigroups ${\sf S}\left({\bf d}^3\right)$ obtained by Herzog 
\cite{herz70} and Fr\"oberg \cite{frob94}. 

We also prove a theorem (Theorem \ref{theo10}) on the upper bound of the number 
of non--zero coefficients in numerator $Q({\bf d}^m;z)$ of Hilbert series for 
non--symmetric semigroup ${\sf S}\left({\bf d}^m\right),m\geq 4$, that 
essentially enhances the result obtained in \cite{worm86}.
\section{3D Frobenius problem: brief review}
\label{revi}
Start with the 2D Frobenius problem for which $F({\bf d}^2)$, $G({\bf d}^2)$ 
and $Q({\bf d}^2;z)$ were known already to J. Sylvester \cite{sylv84}
\begin{equation}
F({\bf d}^2)=d_1 d_2-d_1-d_2\;,\;\;\;\;G({\bf d}^2)=
\frac{1}{2}(d_1-1)(d_2-1)\;,\;\;\;\;Q({\bf d}^2;z)=1-z^{d_1d_2}\;.
\label{sylves1}
\end{equation}
Recall some basic results on the Frobenius problem for semigroup ${\sf S}\left(
{\bf d}^3\right)$ following \cite{herz70}, \cite{kraf85}, \cite{frob94} and 
\cite{denh03}. Let ${\sf S}\left({\bf d}^3\right)$ be a non--symmetric 
semigroup with the 1st minimal relation ${\cal R}_1\left({\bf d}^3\right)$ 
defined by (\ref{herznon1a}), (\ref{herznon3}). Such relations always exist due 
to the finiteness of the Frobenius number $F\left({\bf d}^2\right)$. Then by 
\cite{herz70} the kernel ${\cal I}_3$ is generated by $p_1,p_2,p_3$, where
\begin{eqnarray}
p_1=X_1^{a_{11}}-X_2^{a_{12}}X_3^{a_{13}}\;,\;\;\;
p_2=X_2^{a_{22}}-X_1^{a_{21}}X_3^{a_{23}}\;,\;\;\;
p_3=X_3^{a_{33}}-X_1^{a_{31}}X_2^{a_{32}}\;,\;\;\;\pi(p_i)=0\;.
\label{herznon01}
\end{eqnarray}
Represent (\ref{herznon1a}) as a matrix equation
\begin{eqnarray} 
\widehat {\cal A}_3\left(\begin{array}{r}
d_1\\d_2\\d_3 \end{array}\right)=\left(\begin{array}{r}
0\\0\\0\end{array}\right)\;,\;\;\;
\widehat {\cal A}_{3}=\left(\begin{array}{rrr}
a_{11} & -a_{12} & -a_{13} \\
-a_{21} & a_{22} & -a_{23} \\
-a_{31} & -a_{32} & a_{33} \end{array}\right)\;,\;\;\;
\left\{\begin{array}{r}
\gcd(a_{11},a_{12},a_{13})=1\\
\gcd(a_{21},a_{22},a_{23})=1\\
\gcd(a_{31},a_{32},a_{33})=1\end{array}\right.\;,
\label{herznon2}
\end{eqnarray}  
and establish {\em the standard forms} of the matrix $\widehat {\cal A}_{3}$ 
satisfying (\ref{herznon3}), (\ref{herznon2}). 
\subsection{3D non--symmetric semigroups}
\label{nonsymsub}
First, let all off--diagonal entries of $\widehat {\cal A}_{3}$ be negative 
integers, $a_{ij}\in \{1,2,\ldots\},i\neq j$ i.e. omitting 0. Then, as was 
shown by Johnson \cite{john60}, it leads necessarily to the following
\begin{eqnarray}
&&a_{11}=a_{21}+a_{31}\;,\;\;\;a_{22}=a_{12}+a_{32}\;,\;\;\;
a_{33}=a_{13}+a_{23}\;,\label{herznon3aaa}\\
&&\det \left(\begin{array}{rr}a_{22} & -a_{23} \\
-a_{32} & a_{33}\end{array}\right)=d_1\;,\;\;
\det \left(\begin{array}{rr}a_{11} & -a_{13} \\
-a_{31} & a_{33}\end{array}\right)=d_2\;,\;\;
\det \left(\begin{array}{rr}a_{11} & -a_{12} \\
-a_{21} & a_{22}\end{array}\right)=d_3\;.\label{herznon3aa}
\end{eqnarray}
The ordering (\ref{definn0}) of integers, $d_1<d_2<d_3$, imposes additional 
constraints on the elements $a_{ij}$
\begin{eqnarray}
a_{11}>a_{12}+a_{13}\;,\;\;\;a_{22}>a_{23}\;,\;\;\;a_{33}<a_{31}+a_{32}\;.
\label{addineq1}
\end{eqnarray}
Denote $\widehat {\cal A}_3$ satisfying (\ref{herznon3aaa}) and 
(\ref{herznon3aa}) by $\widehat {\cal A}_3^{(n)}$ and call it {\em the 
standard form} for non--symmetric semigroup ${\sf S}\left({\bf d}^3\right)$. 
Formula (\ref{herznon1a}) together with (\ref{herznon3aaa}) and (\ref{herznon3aa}) 
make it possible to show that at least one of the $a_{ii}$ exceeds 2. The proof 
is obtained by way of contradiction. Let all $a_{ii}=2$. Then due to 
(\ref{herznon3aaa}) we have $a_{ij}=1,i\neq j$, or in accordance with 
(\ref{herznon1a})
$$ 
2d_1=d_2+d_3\;,\;\;2d_2=d_3+d_1\;,\;\;2d_3=d_1+d_2\;\;\;\;\longrightarrow\;\;\;
\;d_1=d_2=d_3\;,
$$
that violates the minimality of $(d_1,d_2,d_3)$. This implies an inequality 
$a_{11}a_{22}a_{33}\geq 12$. Note that the above consideration does not exclude 
the possibility that the diagonal elements $a_{ii}$ coincide in pairs, e.g. 
$a_{11}=a_{22}$, and, moreover, to completely coincide, $a_{ii}=a$. The latter 
kind of degeneration reduces significantly the number of different admissible 
triples $d_1,d_2,d_3$ being a minimal set and satisfying (\ref{herznon3aaa}), 
(\ref{herznon3aa}) (see Appendix \ref{appendix1}). 

The Frobenius number $F\left({\bf d}^3\right)$ for non--symmetric semigroup was 
found for the first time in \cite{herz70} (see also \cite{frob94}) calculating 
only the largest degree of $H({\bf d}^m;z)$ (without calculating Hilbert 
series itself)
\begin{eqnarray}
F\left({\bf d}^3\right)+\sum_{i=1}^3d_i&=&\max\left\{a_{11}d_1+a_{32}d_2;
a_{22}d_2+a_{31}d_1\right\}=\max\left\{a_{22}d_2+a_{13}d_3;a_{33}d_3+a_{12}d_2
\right\}\nonumber\\&=&\max\left\{a_{33}d_3+a_{21}d_1;a_{11}d_1+a_{23}d_3
\right\}\;.\label{herznon1}
\end{eqnarray}
The genus $G\left({\bf d}^3\right)$ of non--symmetric semigroup was calculated 
in algebraic geometry \cite{kraf85}. Dealing with the singularity degrees of 
the monomial space curve whose corresponding semigroup is
${\sf S}\left({\bf d}^3\right)$, Kraft \cite{kraf85} was able to calculate its 
Milnor number $\mu\left({\bf d}^3\right)$ which in the unibranch case is twice 
larger than $G\left({\bf d}^3\right)$ and given by (\ref{herznon5}). Thus, 
(\ref{herznon5}) gives a generalization of the Milnor number for the monomial 
plane curves presented in \cite{miln67}
\begin{eqnarray}
\mu\left({\bf d}^2\right)=1+\sum_{i=1}^2a_{ii}d_i-\prod_{i=1}^2a_{ii}-
\sum_{i=1}^2d_i=(d_1-1)(d_2-1)\;,\;\;\;\gcd(d_1,d_2)=1\;.\label{herznon6}
\end{eqnarray}
As for Hilbert series, partial progress was achieved by
Sz\'ekely and Wormald \cite{worm86} who proved that $Q\left({\bf d}^3;z\right)$
consists of only a limited number of terms (at most twelve) independent of the
values of $d_1,d_2$ and $d_3$. Recently this result was essentially refined by 
Denham \cite{denh03} who gave an algorithm to compute the Hilbert series $H(
{\bf d}^3;z)$ of a graded subring ${\sf k}\left[{\sf S}\left({\bf d}^3\right)
\right]$ for non--symmetic semigroups and established a universal property of 
these series: $Q\left({\bf d}^3;z\right)$ has exactly six terms where the first 
four of them read $1-\sum_{i=1}^3z^{a_{ii}d_i}$. The attempts \cite{barv02} to 
extend further the algorithmic procedure to higher $m$ results only in the 
estimation of the polynomial time of computation of $H\left({\bf d}^m;z
\right)$. 
\subsection{3D symmetric semigroups}
\label{symsub}
The number of independent entries $a_{ij}$ in (\ref{herznon2}) can be reduced 
if at least one off--diagonal element of $\widehat {\cal A}_{3}$ vanishes, e.g. 
$a_{13}=0$ and therefore $a_{11}d_1=a_{12}d_2$. Due to {\em minimality} of the 
last relation we have from (\ref{herznon1a}) the following equalities and 
consequently the matrix representation \cite{herz70}
\begin{eqnarray}
a_{11}=a_{21}=\frac{{\sf lcm}(d_1,d_2)}{d_1},\;\;
a_{12}=a_{22}=\frac{{\sf lcm}(d_1,d_2)}{d_2},\;\;a_{23}=0,\;\;
\widehat {\cal A}_3^{(s)}=\left(\begin{array}{rrr}
a_{11} & -a_{22} & 0 \\
-a_{11} & a_{22} & 0 \\
-a_{31} & -a_{32} & a_{33} \end{array}\right).
\label{herznon3d}
\end{eqnarray}
Call $\widehat {\cal A}_3^{(s)}$ the standard form for the symmetric semigroup 
${\sf S}\left({\bf d}^3\right)$. The kernel ${\cal I}_3$ has 2 generators 
\cite{herz70}
\begin{eqnarray}
p_1=-p_2=X_1^{a_{11}}-X_2^{a_{22}}\;,\;\;\;
p_3=X_3^{a_{33}}-X_1^{a_{31}}X_2^{a_{32}}\;,
\label{herznon3e}
\end{eqnarray}
and the Frobenius number $F\left({\bf d}^3\right)$ looks like \cite{herz70}
\begin{eqnarray}
F\left({\bf d}^3\right)=
a_{11}d_1+a_{33}d_3-\sum_{i=1}^3d_i=a_{22}d_2+a_{33}d_3-\sum_{i=1}^3d_i=
{\sf lcm}(d_1,d_2)+a_{33}d_3-\sum_{i=1}^3d_i\;.
\label{herznon4}
\end{eqnarray}
The corresponding genus $G\left({\bf d}^3\right)$ can be also simplified (see 
formula (\ref{ari2a}) in Section \ref{nomina0}). If ${\sf S}\left({\bf d}^3
\right)$ is symmetic semigroup then 
${\sf k}\left[{\sf S}\left({\bf d}^3\right)\right]$ is a complete intersection 
\cite{herz70} and the Hilbert series $H_s({\bf d}^3;z)$ reads \cite{stanl79}
\begin{eqnarray}
H_s({\bf d}^3;z)=\frac{(1-z^{a_{22}d_2})(1-z^{a_{33}d_3})}
{(1-z^{d_1})(1-z^{d_2})(1-z^{d_3})}\;.
\label{sylves10}
\end{eqnarray}
It is interesting to interpret (\ref{herznon4}) in the sense of Johnson's 
formula \cite{john60} when $\gcd(d_1,d_2)=k\geq 1$
\begin{eqnarray}
F\left(d_1,d_2,d_3\right)=kF\left(\frac{d_1}{k},\frac{d_2}{k},d_3\right)+
(k-1)d_3\;.
\label{herznon4a}
\end{eqnarray}
Comparison of (\ref{herznon4}) and (\ref{herznon4a}) gives
\begin{eqnarray}
kF\left(\frac{d_1}{k},\frac{d_2}{k},d_3\right)={\sf lcm}(d_1,d_2)-d_1-d_2+
(a_{33}-k)d_3\;.
\label{herznon4c}
\end{eqnarray}
In case $\gcd(d_1,d_2)=1$, this leads to $F\left(d_1,d_2,d_3\right)=F\left(
d_1,d_2\right)+(a_{33}-1)d_3$. Recalling the inequality (\ref{def1b}) for 
conductors $c({\bf d}^3)\leq c({\bf d}^2)$ and their connection with the 
Frobenius numbers we have $F\left(d_1,d_2,d_3\right)\leq F\left(d_1,d_2\right)$ 
that results together with (\ref{sylves1}) and (\ref{herznon4c}), in $a_{33}=1$,
i.e., $d_3$ is representable by $d_1$ and $d_2$. Thus, every semigroup 
generated by three pairwise relatively prime elements cannot be symmetric 
\cite{frob87}.

We finish this Section noting that the minimal set $\{d_1,d_2,d_3\}$, which 
generates the semigroup ${\sf S}\left({\bf d}^3\right)$, cannot 
include 2 as an element. Indeed, assume the opposite, that $d_1=2$ and the 
other two $d_2<d_3$ are both odd integers. Then $d_3-d_2$ is divisible by 2, 
and therefore such set is not minimal in accordance with (\ref{defin2}). In the 
case, when one of $d_2,d_3$ does represent an even integer, the claim is 
clear. Henceforth, we assume that the elements of the minimal set ${\bf d}^3$ 
satisfy 
\begin{eqnarray}
3\leq d_1<d_2<d_3\;. 
\label{assump3}
\end{eqnarray}
Further generalization of (\ref{assump3}) to non--symmetric semigroups 
${\sf S}\left({\bf d}^m\right)$ of higher dimension, $m\geq 4$, will be 
given in Section \ref{appl0}.
\section{Matrix representation of the set $\Delta\left({\bf d}^2\right)$}
\label{repr2}
In this Section we construct the matrix representation of the set $\Delta\left(
{\bf d}^2\right)$ of integers $t$ which are unrepresentable by $d_1,d_2$. We 
start with the important statement about matrix representation which dates back 
to A. Brauer \cite{brau42} and results partly from his discussion with I. Schur
\footnote{We quote from \cite{brau42} : {\em The Theorems in \S  3--5 result 
partly from discussions of Schur and the author. It was formerly intended to 
publish these results in a joint paper. I conform with Schur's wishes that the 
publishing be not longer postponed and that I publish the paper alone}. The 
paper \cite{brau42} was submitted for publication in November 25, 1940, less 
than two months before Schur's death, and was published two years after.}. 
\begin{lemma}{\rm (\cite{brau42})}
Let $d_1$ and $d_2$ be relatively prime positive integers. Then every positive 
integer $s$ not divisible by $d_1$ or by $d_2$ is representable either in the 
form $s=xd_1+yd_2,\;x>0,y>0$ or in the form $s=d_1d_2-pd_1-qd_2,\;p>0,q>0;p,q
\in {\mathbb N}$.
\label{lem1}
\end{lemma}
\begin{definition}
\label{definit1}
Let integers $2<d_1<d_2$ be given. Define function $\sigma(p,q)$ as follows
\begin{eqnarray}
\sigma(p,q):=d_1 d_2-pd_1-qd_2\;.\label{sylves0}
\end{eqnarray}
\end{definition}
The next Lemma specifies the bounds on the values of $p$ and $q$ introduced in 
Lemma \ref{lem1} above. 
\begin{lemma}
Let $t$ be an integer and $d_2>d_1,\;\gcd(d_1,d_2)=1$. Then $t\in\Delta\left(
{\bf d}^2\right)$ iff $t$ is uniquely representable as
\begin{eqnarray}
t=\sigma(p,q)\;,\;\;\;\;\mbox{where}
\label{sylves2}
\end{eqnarray}
\begin{eqnarray}
1\leq p\leq \left\lfloor d_2-\frac{d_2}{d_1}\right\rfloor\;\;\;\;
1\leq q\leq d_1-1\;,\;\;\;\;\mbox{and}\;\;\;\;
d_1-1\leq \left\lfloor d_2-\frac{d_2}{d_1}\right\rfloor\;.\label{sylves3}
\end{eqnarray}
\label{lem2}
\end{lemma}
{\sf Proof} $\;\;\;$
In accordance with Lemma \ref{lem1} every integer $t$ which is unrepresentable 
by $d_1,d_2,\;\gcd(d_1,d_2)=1$ is representable by (\ref{sylves0}) and 
(\ref{sylves2}). Thus, the Frobenius number is $F\left({\bf d}^2\right)=\sigma(
1,1)$. The restrictions (\ref{sylves3}) come from simple considerations
\begin{eqnarray}
\sigma(p,1)> 0\;\;\rightarrow\;\;p\leq \left\lfloor 
d_2-\frac{d_2}{d_1}\right\rfloor\;\;,\;\;\;\;
\sigma(1,q)> 0\;\;\rightarrow\;\;q\leq
\left\lfloor d_1-\frac{d_1}{d_2}\right\rfloor=d_1-1\;.\nonumber
\end{eqnarray}
The presentation of $\sigma(p,q)$ by (\ref{sylves0}) is unique. A standard
proof of uniqueness of (\ref{sylves0}) is to assume, by way of contradiction,
that there are two such representations $\sigma(p_1,q_1)$, $\sigma(p_2,q_2)$
and consequently, 
$$
(p_1-p_2)d_1=(q_2-q_1)d_2,\;\;1\leq p_1,p_2\leq \left\lfloor
d_2-\frac{d_2}{d_1}\right\rfloor,\;\;
1\leq q_1,q_2\leq d_1-1\;\rightarrow\;|p_1-p_2|<d_2,\;\;
|q_2-q_1|<d_1.
$$
But this is impossible since $d_1$ and $d_2$ have no common factors. 

Finally prove the last inequality in (\ref{sylves3}). Assuming $d_2\geq d_1+2$ 
we get
$$
\left\lfloor d_2-\frac{d_2}{d_1}\right\rfloor -(d_1-1)\geq
d_2-\frac{d_2}{d_1}-1-(d_1-1)=\frac{(d_2-d_1)(d_1-1)}{d_1}-1\geq
2\frac{d_1-1}{d_1}-1=1-\frac{2}{d_1}> 0\;.
$$
In the case $d_2=d_1+1$ we obtain
$$
\left\lfloor d_2-\frac{d_2}{d_1}\right\rfloor-(d_1-1)=
\left\lfloor d_1+1-\frac{d_1+1}{d_1}\right\rfloor-(d_1-1)=
\left\lfloor d_1-\frac{1}{d_1}\right\rfloor-(d_1-1)=0\;.
$$
Thus, combining both cases we arrive at (\ref{sylves3}). This completes the
proof of the Lemma.$\;\;\;\;\;\;\Box$

Show that the integers $\sigma(p,q)$ given by (\ref{sylves2}), (\ref{sylves3})
exhaust all integers unrepresentable by $d_1$ and $d_2$, or, in other words,  
they give the genus $G({\bf d}^2)$ obtained by Sylvester \cite{sylv84} and 
given in (\ref{sylves1}). Indeed, counting the number of integers $\sigma(
p,q)$ with the above properties (\ref{sylves2}) and (\ref{sylves3}) 
successively over the index set $q=1,\ldots, d_1-1$ one gets
\begin{eqnarray}
G({\bf d}^2)&=&\sum_{q=1}^{d_1-1}\left\lfloor d_2-q\frac{d_2}{d_1}\right\rfloor=
\sum_{q=d_1-1}^{1}\left\lfloor q\frac{d_2}{d_1}\right\rfloor=\frac{1}{2}\left(
\sum_{q=1}^{d_1-1}\left\lfloor q\frac{d_2}{d_1}\right\rfloor+
\sum_{q=1}^{d_1-1}\left\lfloor (d_1-q)\frac{d_2}{d_1}\right\rfloor\right)
\nonumber\\
&=&\frac{1}{2}\sum_{q=1}^{d_1-1}\left(\left\lfloor q\frac{d_2}{d_1}\right
\rfloor+\left\lfloor d_2-q\frac{d_2}{d_1}\right\rfloor\right)=\frac{(d_2-1)
(d_1-1)}{2}\;,\nonumber
\end{eqnarray}
that follows from the equalities
\begin{eqnarray}
q\frac{d_2}{d_1}=\left\lfloor q\frac{d_2}{d_1}\right\rfloor+
\left\{ q\frac{d_2}{d_1}\right\},\;\;
\left\lfloor d_2-q\frac{d_2}{d_1}\right\rfloor=
\left\lfloor d_2-\left\{q\frac{d_2}{d_1}\right\}\right\rfloor
-\left\lfloor q\frac{d_2}{d_1}\right\rfloor=d_2-1-
\left\lfloor q\frac{d_2}{d_1}\right\rfloor,\;\;q<d_1\;.\label{schur44}
\end{eqnarray}
In (\ref{schur44}) we denote by $\{b\}$ the fractional part of a real number $b$.

\begin{figure}[h]
\centerline{\psfig{figure=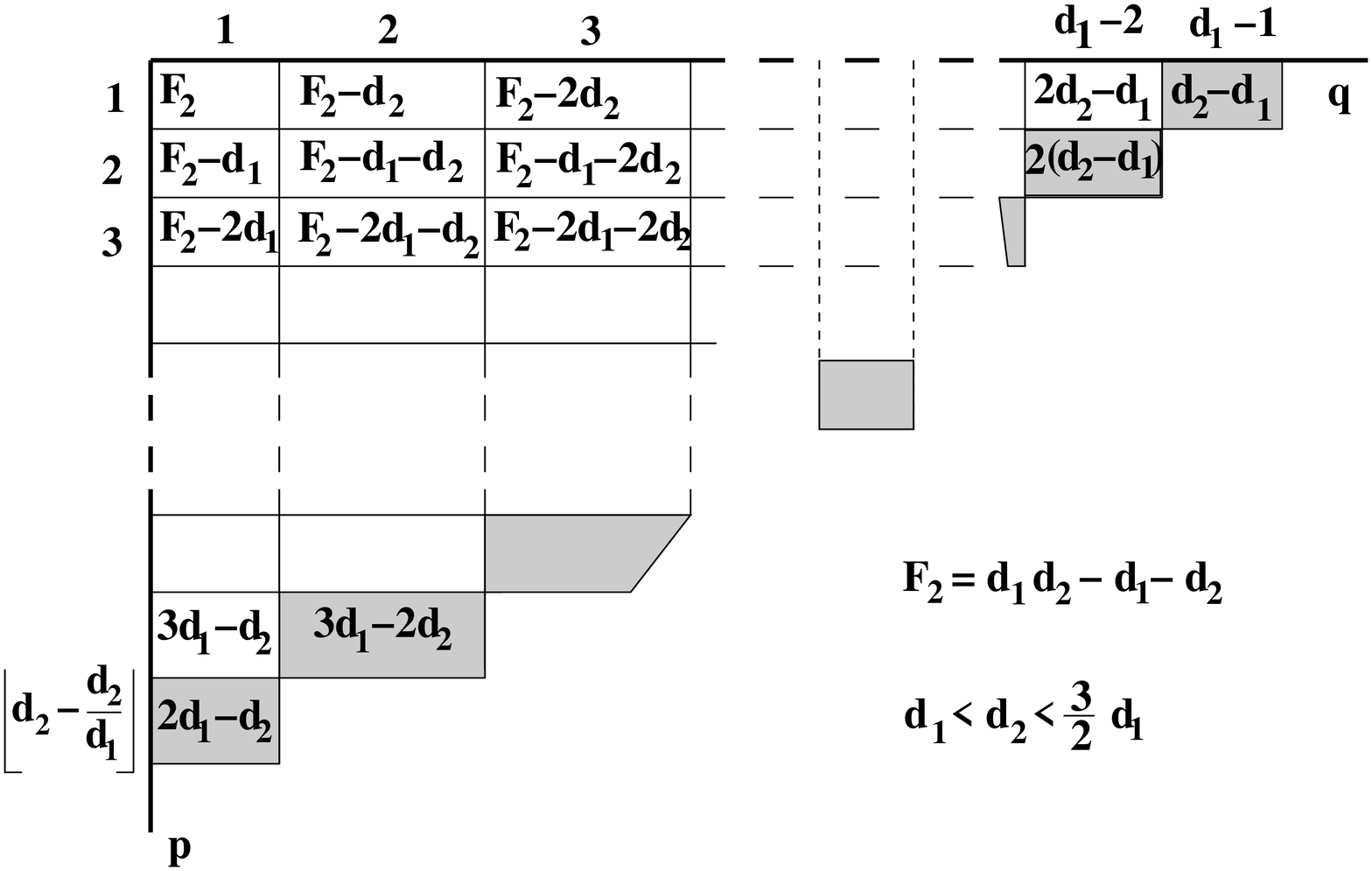,height=5.5cm,width=14cm}}
\caption{Typical matrix representation $M\left\{\Delta({\bf d}^2)\right\}$ of 
the set $\Delta(d_1,d_2)$ for the case $d_1<d_2<3/2 d_1$. The bottom layer 
${\sf BL}_M\left\{\Delta({\bf d}^2)\right\}$ is marked in {\em gray color}. 
The top layer ${\sf TL}_M\left\{\Delta({\bf d}^2)\right\}$ coincides with the 
highest row of diagram.} 
\label{repr1}
\end{figure}

The representation (\ref{sylves2}) of all integers $\sigma(p,q)\in \Delta(d_1,
d_2)$ is called {\em the matrix representation} of the set $\Delta(d_1,d_2)$ 
and is denoted by $M\left\{\Delta({\bf d}^2)\right\}$ (see Figure \ref{repr1})
\begin{eqnarray}
\sigma\left\{M\left\{\Delta({\bf d}^2)\right\}\right\}=\Delta({\bf d}^2)\;.
\label{schur55}
\end{eqnarray}
$\sigma(p,q)$ is the integer which occurs in the row $p$ and the column 
$q$ of $M\left\{\Delta({\bf d}^2)\right\}$, e.g. $\sigma(1,1)=d_1d_2-d_1-d_2$. 

Based on 
$M\left\{\Delta({\bf d}^2)\right\}$ introduce two sets which will be important 
in the coming Sections. We call the totality of the lowest cells in every 
column of $M\left\{\Delta({\bf d}^2)\right\}$ {\em the bottom layer of $M\left\{
\Delta({\bf d}^2)\right\}$} and denote it by ${\sf BL}_M\left\{\Delta({\bf d}^2
)\right\}$. We also call the totality of the highest cells in every column of 
$M\left\{\Delta({\bf d}^2)\right\}$ {\em the top layer of $M\left\{\Delta({\bf 
d}^2)\right\}$} and denote it by ${\sf TL}_M\left\{\Delta({\bf d}^2)\right\}$.

We establish the structure of both sets ${\sf BL}_M\left\{\Delta({\bf d}^2)
\right\}$ and ${\sf TL}_M\left\{\Delta({\bf d}^2)\right\}$. First, prove 
the following Lemma.
\begin{lemma}
\label{lem3}
For every number $k,\;1\leq k\leq d_1-1$ there exists $(p,q)\in {\sf BL}_M
\left\{\Delta({\bf d}^2)\right\}$ such that $\sigma(p,q)=k$.
\end{lemma}
{\sf Proof} $\;\;\;$Let $1\leq q\leq d_1-1$ and define 
\begin{eqnarray}
p_b(q)=\max\left\{1\leq p\;\bracevert\;\sigma(p,q)>0\right\}\;,
\label{sylves7}
\end{eqnarray}
where subscript $"b"$ stands for "bottom". We may consider $p_b$ as a function 
of $q$. Then $(p_b(q),q)\in {\sf BL}_M\left\{\Delta({\bf d}^2)\right\}$. Let us 
derive the function $p_b(q)$. It follows from (\ref{sylves0}) that
$$
\sigma(p,q)>0\;\;\rightarrow\;\;p<d_2-\frac{d_2}{d_1}q\;\;\rightarrow\;\;
p_b(q)=\left\lfloor d_2-\frac{d_2}{d_1}q\right\rfloor\;.
$$
Hence, according to (\ref{schur44}) we get $p_b(q)=d_2-1-\left\lfloor q
\frac{d_2}{d_1}\right\rfloor$, and further
\begin{eqnarray}
\sigma(p_b(q),q)=d_1-qd_2+d_1\left\lfloor q\frac{d_2}{d_1}
\right\rfloor=d_1-d_1\left\{q\frac{d_2}{d_1}\right\}\;.
\label{sylves7a}
\end{eqnarray}
For $1\leq q\leq d_1-1$ we have $1/d_1\leq \left\{q\frac{d_2}{d_1}\right\}\leq
(d_1-1)/d_1$. Combining this with (\ref{sylves7a}) gives the bounds 
for $\sigma(p_b(q),q)$
\begin{eqnarray}
1\leq \sigma(p_b,q)\leq d_1-1\;.\label{sylves8}
\end{eqnarray}
Due to the uniqueness of the presentation of $\sigma(p,q)$ by (\ref{sylves0}), 
the bounds (\ref{sylves8}) lead to the conclusion that ${\sf BL}_M\left\{\Delta
({\bf d}^2)\right\}$ is occupied exclusively by the integers $1,\ldots,d_1-1$ 
not in a necessarily consecutive order. This proves the Lemma. 
$\;\;\;\;\;\;\Box$

As for the top layer, ${\sf TL}_M\left\{\Delta({\bf d}^2)\right\}$ coincides 
with the highest row in $M\left\{\Delta({\bf d}^2)\right\}$. Thus, finally we 
can write, in accordance with (\ref{schur55}),
\begin{eqnarray}
\sigma\left\{{\sf BL}_M\left\{\Delta({\bf d}^2)\right\}\right\}=\{1,\ldots,
d_1-1\}\;,\;\;\;\sigma\left\{{\sf TL}_M\left\{\Delta({\bf d}^2)\right\}\right\}=
\{\sigma(1,1),\ldots,\sigma(1,d_1-1)\}\;.
\label{sylves9}   
\end{eqnarray}
where $\sigma(1,1)=F({\bf d}^2)$ and $\sigma(1,d_1-1)=d_2-d_1$.
\section{Matrix representation of the set $\Delta\left({\bf d}^3\right)$}
\label{repr3}
In this Section we construct the set $\Delta\left({\bf d}^3\right)$ out of the 
set $\Delta\left({\bf d}^2\right)$. Introduce new objects -- {\it associated 
sets} $\Omega_{d_3}^k({\bf d}^2)$ and, based on it, define the matrix 
representation of the set $\Delta\left({\bf d}^3\right)$. Here $k$ is an 
integer variable (see discussion later). This construction paves the way 
to solve the 3D Frobenius problem for non--symmetric semigroup ${\sf S}\left(
{\bf d}^3\right)$.

Following (\ref{def1b}) and the definition (\ref{condfrob}) of the Frobenius 
number we have $F({\bf d}^3)\leq F({\bf d}^2)$. In the coming Lemma we show 
that equality does not occur.
\begin{lemma}
\label{lem4}
Let ${\bf d}^2$ be given, ${\bf d}^2=(d_1,d_2)$ and $d_1,d_2,\;\gcd(d_1,d_2)=1$.
Every integer $d_3\in\Delta({\bf d}^2)$ gives rise to the minimal generating set
$\{d_1,d_2,d_3\}$, which generates the semigroup ${\sf S}\left({\bf d}^3\right)$
such that 
\begin{eqnarray}
F({\bf d}^3)<F({\bf d}^2)\;.
\label{essen2}
\end{eqnarray}
\end{lemma}
{\sf Proof}$\;\;\;$Let $d_3\in \Delta({\bf d}^2)$. Then due to Lemma \ref{lem2}
with $t=d_3$ there exist $p_{d_3}$ and $q_{d_3}$ satisfying (\ref{sylves3}) 
such that
\begin{eqnarray}
d_3=d_1d_2-p_{d_3}d_1-q_{d_3}d_2\;,
\label{d3ext1}
\end{eqnarray}
and $d_3$ is unrepresentable by $d_1,d_2$ (see (\ref{sylves0})). Thus, the 
triple $\{d_1,d_2,d_3\}$ represents the minimal set generating ${\sf S}\left(
{\bf d}^3\right)$ in accordance with (\ref{defin2}). Define the set $\Omega_{
d_3}^1({\bf d}^2)$ of integers $A_1$ in $\Delta({\bf d}^2)$ representable by 
$d_1,d_2,d_3$ as follows
\begin{eqnarray}
\Omega_{d_3}^1({\bf 
d}^2)=\left\{A_1\;\bracevert\;A_1=u_1d_1+u_2d_2+d_3,\;0\leq u_1\leq p_{d_3}-1,
0\leq u_2\leq q_{d_3}-1\right\}\;.
\label{d3ext2}
\end{eqnarray}
$A_1$ depends on $u_1$ and $u_2$, hence we shall write $A_1=A_1(u_1,u_2)$. 
Since $\Delta({\bf d}^3)$ consists of the integers unrepresentable by 
$d_1,d_2,d_3$ it is clear that 
\begin{eqnarray}  
\Omega_{d_3}^1({\bf d}^2)\bigcap \Delta({\bf d}^3)=\emptyset\;.
\label{d3ext3}
\end{eqnarray}
It follows from (\ref{sylves0}) that $A_1(u_1,u_2)=\sigma(p_{d_3}-u_1,q_{d_3}-
u_2)$. By expressions (\ref{sylves1}) and (\ref{d3ext1}) we have $A_1(p_{d_3}-
1,q_{d_3}-1)=F({\bf d}^2)$. In particular, 
\begin{eqnarray}
F({\bf d}^2)\in \Omega_{d_3}^1({\bf d}^2)\;. 
\label{d3ext4}  
\end{eqnarray}
Since $F\left({\bf d}^2\right)\stackrel{def}=\max\left\{t\in\Delta({\bf d}^2)
\right\}$ and $F({\bf d}^2)\in \Omega_{d_3}^1({\bf d}^2)$ by (\ref{d3ext4}), 
hence due to $\Omega_{d_3}^1({\bf d}^2)\subset \Delta({\bf d}^2)$ we get 
$$
\max\left\{t\in \Omega_{d_3}^1({\bf d}^2)\right\}=\max\left\{t\in\Delta({\bf 
d}^2)\right\}\;.
$$
By (\ref{def1b}) we have $\Delta({\bf d}^3)\subset\Delta({\bf d}^2)$, hence 
$\max\left\{t\in \Delta({\bf d}^3)\right\}\leq \max\left\{t\in\Delta({\bf 
d}^2)\right\}$. However, $F\left({\bf d}^2\right)\stackrel{def}=\max\left\{t\in
\Delta({\bf d}^2)\right\}$ and $F({\bf d}^2)\in \Omega_{d_3}^1({\bf d}^2)$ by 
(\ref{d3ext4}), so it follows from (\ref{d3ext3}) that 
$$
\max\left\{t\in \Delta({\bf d}^3)\right\}<\max\left\{t\in\Delta({\bf d}^2)
\right\}
$$ 
that proves Lemma.$\;\;\;\;\;\;\Box$

It may happen that the set $\Omega_{d_3}^1({\bf d}^2)$ described in Lemma 
\ref{lem4} does not exhaust all elements of $\Delta({\bf d}^2)$ which are 
representable by $d_1,d_2,d_3$. 
\begin{lemma} 
\label{lem5}
If the integers $d_3$ and $2d_3$ are unrepresentable by $d_1,d_2$ then $2d_3
\not\in\Omega_{d_3}^1({\bf d}^2)$.
\end{lemma}
{\sf Proof} $\;\;\;$The proof follows by way of contradiction. Let $2d_3\in
\Omega_{d_3}^1({\bf d}^2)$, then due to (\ref{d3ext2}) there exist nonnegative 
integers $u_1$ and $u_2$ such that
$$
2d_3=u_1d_1+u_2d_2+d_3\;,\;\;\;\mbox{or}\;\;\;d_3=u_1d_1+u_2d_2\;,
$$
that violates the minimality of $\{d_1,d_2,d_3\}$.$\;\;\;\;\;\;\Box$
\subsection{Associated sets $\Omega_{d_3}^k({\bf d}^2)$}
\label{asssso1}
In order to account for all integers which contribute to the construction of 
$\Delta({\bf d}^3)$ we have to extend the set $\Omega_{d_3}^1({\bf d}^2)$. 
First, recall from (\ref{herznon1a}) and (\ref{herznon3}) one of the 1st minimal
relations ${\cal R}_1\left({\bf d}^3\right)$ for a given ${\bf d}^3$: $a_{33}
d_3=a_{31}d_1+a_{32}d_2$, where $a_{33}\geq 2$. 
\begin{definition}
\label{definit2}
Let $d_3\in \Delta({\bf d}^2)$ with representation $d_3=d_1d_2-p_{d_3}d_1-
q_{d_3}d_2$ where $p_{d_3}$ and $q_{d_3}$ satisfy (\ref{sylves3}). Let $k$ be 
a positive integer, $1\leq k< a_{33}$. Define the set $\Omega_{d_3}^k({\bf 
d}^2)$ of integers $A_k$ in $\Delta({\bf d}^2)$
\begin{eqnarray}
\Omega_{d_3}^k({\bf d}^2)=\left\{A_k\;\bracevert\;A_k=u_1d_1+u_2d_2+kd_3\;,
0\leq u_1\leq p_{kd_3}-1,0\leq u_2\leq q_{kd_3}-1\right\}\;.
\label{assoc1}
\end{eqnarray}
Call $\Omega_{d_3}^k({\bf d}^2)$ a $kd_3$--associated set.
\end{definition}
$A_k$ depends on $u_1$ and $u_2$, hence we shall write $A_k=A_k(u_1,u_2)$. 
Taking $u_1=p_{kd_3}-1,u_2=q_{kd_3}-1$ gives 
\begin{eqnarray}
F\left({\bf d}^2\right)=A_k\left(p_{kd_3}-1,q_{kd_3}-1\right)\in \Omega_{d_3}
^k({\bf d}^2)\;,\;\;\;1\leq k<a_{33}\;.
\label{assoc1a}
\end{eqnarray}

\begin{figure}[h]
\centerline{\psfig{figure=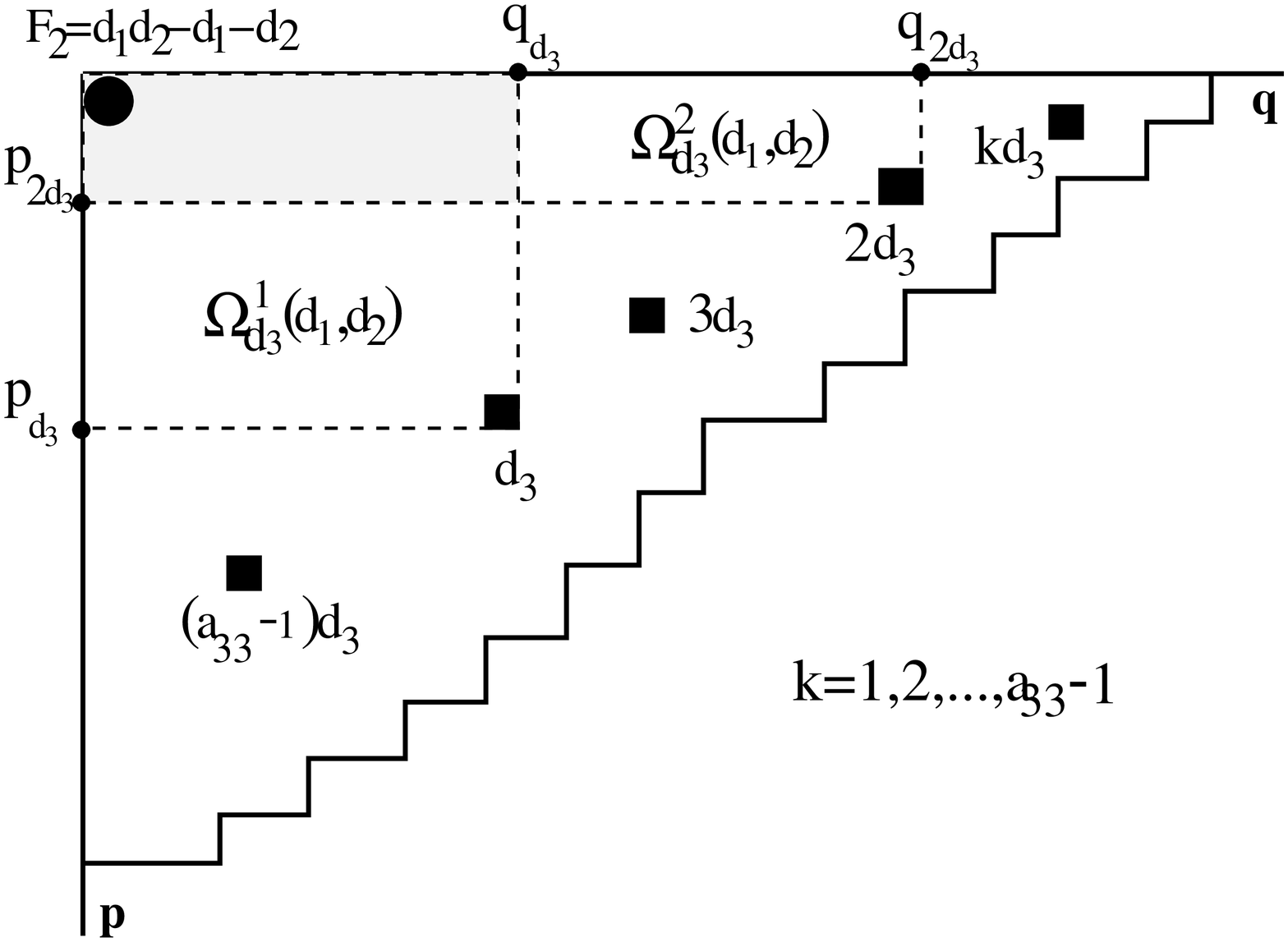,height=5.5cm,width=14cm}}
\caption{Typical matrix representation $M\left\{\Delta({\bf d}^2)\right\}$ of 
the set $\Delta({\bf d}^2)$. Matrix representations $M\left\{\Omega_{d_3}
^k({\bf d}^2)\right\}$ of the $kd_3$--associated sets $\Omega_{d_3}^1({\bf d}
^2)$ and $\Omega_{d_3}^2({\bf d}^2)$ are drawn by {\em dashed lines}. Their 
intersection ({\em gray color}) contains the Frobenius number $F\left({\bf d}^
2\right)$ which is marked by {\em black oval}. The integers $kd_3\not\in 
\Delta({\bf d}^3),1\leq k<a_{33}$ are shown by {\em black boxes}.}
\label{repr2}
\end{figure}

It follows from (\ref{assoc1a}) that the intersection of any two associated 
sets $\Omega_{d_3}^j({\bf d}^2)$ and $\Omega_{d_3}^k({\bf d}^2),1\leq j,k<
a_{33}$ is non--empty set. 
As follows from (\ref{assoc1}) the matrix representation $M\left\{\Omega_{d_3}
^k({\bf d}^2)\right\}$ of $\Omega_{d_3}^k({\bf d}^2)$ is assigned by the 
rectangle $[1,p_{kd_3}]\times [1,q_{kd_3}]$ inside $M\left\{\Delta({\bf d}^2)
\right\}$ with cardinality $\#\Omega_{d_3}^k({\bf d}^2)=p_{kd_3}q_{kd_3}$ 
(see Figure \ref{repr2}). 
\begin{theorem}
\label{theo1}
Let $\{d_1,d_2,d_3\}$ be a minimal generating set of the semigroup ${\sf S}
\left({\bf d}^3\right)$ and let $\overline{\sigma}$ be an integer $\overline{
\sigma}\in \Delta({\bf d}^2)$ representable by $d_1,d_2,d_3$. Then there exists 
at least one $k,\;1\leq k<a_{33}$ such that $\overline{\sigma}\in 
\Omega_{d_3}^k({\bf d}^2)$.
\end{theorem}
Before giving the proof of Theorem \ref{theo1} let us present an auxiliary 
Lemma based on the theory of restricted partition function $W(\overline{
\sigma},{\bf d}^m)$. Recall the main recursion relation \cite{fel02} for $W(
\overline{\sigma},{\bf d}^m)$ which gives the number of partitions of 
$\overline{\sigma}$ into positive integers $d_1,\ldots,d_m$ each not greater
than $\overline{\sigma}$. Then
\begin{eqnarray}
W(\overline{\sigma},{\bf d}^m) - W(\overline{\sigma}-d_m,{\bf d}^m) =
W(\overline{\sigma},{\bf d}^{m-1})\;,\;\;\;\;
{\bf d}^{m-1}=\{d_1,\ldots,d_{m-1}\}\;.\label{SW_recursion}
\end{eqnarray}
\begin{lemma}
\label{lem6}
Let $\{d_1,\ldots,d_m\}$ be a minimal generating set of ${\sf S}\left({\bf d}
^m\right)$ and let $\Delta({\bf d}^m)$ be the corresponding set of 
unrepresentable integers. If $\overline{\sigma}\in \Delta({\bf d}^m)$ and 
$\overline{\sigma}-d_k>0,\;1\leq k\leq m$, then necessarily $\overline{\sigma}-
d_k\in \Delta({\bf d}^m)$. 
\end{lemma}
{\sf Proof}$\;\;\;$Assume first that $k=m$. If $\overline{\sigma}\in\Delta(
{\bf d}^m)$ then $W(\overline{\sigma},{\bf d}^m)=W(\overline{\sigma},{\bf d}
^{m-1})=0$ and consequently $W(\overline{\sigma}-d_m,{\bf d}^m)=0$ due to 
(\ref{SW_recursion}). The latter implies $\overline{\sigma}-d_m\in 
\Delta({\bf d}^m)$. 

Now let $k$ be arbitrary, $1\leq k<m$. 
The validity of the relation (\ref{SW_recursion}) does not depend on the 
position of $d_k$ in the tuple ${\bf d}^m$. Thus, resorting the tuple ${\bf d}
^m$ in such a way that $d_k$ becomes the last in the list $d_1,
\ldots,d_m$ and repeating the above consideration, we come to the proof of the 
Lemma. $\;\;\;\;\;\;\Box$

Note that Lemma \ref{lem6} states the necessary but not sufficient 
requirement 
for $\overline{\sigma}$, i.e. an opposite implication $\overline{\sigma}-d_k\in 
\Delta({\bf d}^m)\;\rightarrow\;\overline{\sigma}\in \Delta({\bf d}^m)$ is 
not true. 

Now we return to the proof of Theorem \ref{theo1}.

\vspace{.2cm}
\noindent
{\sf Proof of Theorem \ref{theo1}}$\;\;\;$Let ${\cal R}_1\left({\bf d}^3\right)$ 
be the 1st minimal relation defined in (\ref{herznon3}). Then $kd_3\in \Delta(
{\bf d}^2),\;1\leq k<a_{33}$. Consider an integer $\overline{\sigma}\in \Delta(
{\bf d}^2)$ representable by $d_1,d_2,d_3$
\begin{eqnarray}
\overline{\sigma}=\alpha_1 d_1+\alpha_2 d_2+\alpha_3 d_3\;,\;\;\;
\alpha_1,\alpha_2,\alpha_3\in {\mathbb N}\cup \{0\}\;.\label{frac1}
\end{eqnarray}
It follows from (\ref{herznon3}) that $\alpha_3$ is not divisible by $a_{33}$, 
otherwise
$$
\overline{\sigma}=\alpha_1 d_1+\alpha_2 d_2+\left\lfloor\frac{\alpha_3}{a_{33}}
\right\rfloor a_{33}d_3=\left(\alpha_1+\left\lfloor\frac{\alpha_3}{a_{33}}
\right\rfloor a_{31}\right)d_1+\left(\alpha_2+\left\lfloor
\frac{\alpha_3}{a_{33}}\right\rfloor a_{32}\right)d_2\not\in\Delta({\bf 
d}^2)\;, 
$$
that contradicts our assumption $\overline{\sigma}\in \Delta({\bf d}^2)$. 

We are going to show that $\overline{\sigma}\in\Omega_{d_3}^{\alpha_3}({\bf 
d}^2)$. To this end we have to show 
\begin{eqnarray}
0\leq \alpha_1\leq p_{\alpha_3 d_3}-1\;,\;\;\;0\leq\alpha_2\leq
q_{\alpha_3d_3}-1\;.\label{frac1a}
\end{eqnarray}
Applying Lemma \ref{lem6} with $m=3$, $\alpha_1$ times with $k=1$ and $\alpha_2$ 
times with $k=2$ we get $\alpha_3 d_3\in \Delta({\bf d}^2)$. Consider 2 cases. 
First, let $1\leq\alpha_3<a_{33}$, then substituting (\ref{d3ext1}) into 
(\ref{frac1}) we obtain
\begin{eqnarray}
\overline{\sigma}=\alpha_1 d_1+\alpha_2 d_2+d_1d_2-p_{\alpha_3 d_3}d_1-
q_{\alpha_3 d_3}d_2=d_1d_2-(p_{\alpha_3 d_3}-\alpha_1)d_1-
(q_{\alpha_3 d_3}-\alpha_2)d_2\;.\nonumber
\end{eqnarray}
Applying Lemma \ref{lem2} to the last representation of $\overline{\sigma}$ 
we get
\begin{eqnarray}
p_{\alpha_3 d_3}-\alpha_1\geq 1\;,\;\;\;q_{\alpha_3d_3}-\alpha_2\geq 1\;,
\nonumber
\end{eqnarray}
and combining this with $\alpha_1,\alpha_2\in {\mathbb N}\cup\{0\}$ in 
(\ref{frac1}) one concludes that (\ref{frac1a}) does hold. This 
leads to $\overline{\sigma}\in \Omega_{d_3}^{\alpha_3}({\bf d}^2)$ in accordance
with Definition \ref{definit2}. 

In the second case, consider $\alpha_3>a_{33}$ and represent $\alpha_3 d_3$ as 
follows
\begin{eqnarray}
\alpha_3d_3=a_{33}d_3\left\lfloor\frac{\alpha_3}{a_{33}}\right\rfloor+
a_{33}d_3\left\{\frac{\alpha_3}{a_{33}}\right\}\;.\label{frac2}
\end{eqnarray}
Substituting $a_{33}d_3$ from (\ref{herznon1a}) into (\ref{frac2}) we get
\begin{eqnarray}
\alpha_3 d_3=\left\lfloor
\frac{\alpha_3}{a_{33}}\right\rfloor \left(a_{31}d_1+a_{32}d_2\right)+
a_{33}\left\{\frac{\alpha_3}{a_{33}}\right\}d_3\;.\nonumber
\end{eqnarray}
Further, substituting the above result into (\ref{frac1}), we obtain
\begin{eqnarray}
\overline{\sigma}=\xi_1d_1+\xi_2d_2+\xi_3d_3\;,\label{frac4} 
\end{eqnarray}
where 
\begin{eqnarray}
\xi_1=\alpha_1+a_{31}\left\lfloor\frac{\alpha_3}{a_{33}}\right\rfloor,\;\;
\xi_2=\alpha_2+a_{32}\left\lfloor\frac{\alpha_3}{a_{33}}\right\rfloor,\;\;
\xi_3=a_{33}\left\{\frac{\alpha_3}{a_{33}}\right\}<a_{33},\;\;\;
\xi_1,\xi_2,\xi_3 \in {\mathbb N}\cup \{0\}\;.\label{frac5}
\end{eqnarray}
Comparing (\ref{frac4}), (\ref{frac5}) with (\ref{frac1}) one concludes that 
the second case ($\alpha_3>a_{33},\xi_3<a_{33}$) is reduced to the first one 
($\alpha_3<a_{33}$) and therefore $\overline{\sigma}\in \Omega_{d_3}^{\xi_3}(
{\bf d}^2)$ with $\xi_3$ instead of $\alpha_3$. This completes the proof 
of the Theorem.$\;\;\;\;\;\;\Box$

Finally we are ready to prove the main theorem of this Section.
\begin{theorem}
\label{theo2}
Let ${\bf d}^3$ be given, ${\bf d}^3=(d_1,d_2,d_3)$, and the 1st minimal 
relation ${\cal R}_1\left({\bf d}^3\right)$ is defined by (\ref{herznon1a}). 
The set $\Delta({\bf d}^3)$ coincides with the complement of the union of all 
associated sets $\Omega_{d_3}^k({\bf d}^2)$ in the set of unrepresentable 
integers $\Delta({\bf d}^2)$ where $k=1,\ldots,a_{33}-1$.
\begin{eqnarray}
\Delta({\bf d}^3)=\Delta({\bf d}^2)\setminus \left\{
\bigcup_{k=1}^{a_{33}-1} \Omega_{d_3}^k({\bf d}^2)\right\}\;.
\label{emp2}
\end{eqnarray}
\end{theorem}
{\sf Proof}$\;\;\;$First, we show that 
\begin{eqnarray}
\Delta({\bf d}^2)\setminus \left\{\bigcup_{k=1}^{a_{33}-1} 
\Omega_{d_3}^k({\bf d}^2)\right\}\subseteq \Delta({\bf d}^3)\;.
\label{sub1}
\end{eqnarray}
Let $\overline{\sigma}\in \Delta({\bf d}^2)\setminus \left\{\bigcup_{k=1}^
{a_{33}-1}\Omega_{d_3}^k({\bf d}^2)\right\}$ and suppose $\overline{\sigma}
\not\in \Delta({\bf d}^3)$. Then by definition of $\Delta({\bf d}^3)$ 
$\overline{\sigma}$ is representable by $d_1,d_2,d_3$. Hence, by Theorem 
\ref{theo1}, we have $\overline{\sigma}\in \bigcup_{k=1}^{a_{33}-1}\Omega_{d_3}
^k({\bf d}^2)$ that contradicts our assumption on $\overline{\sigma}$. 
Consequently, (\ref{sub1}) holds true.

Finally we show that
\begin{eqnarray}
\Delta({\bf d}^2)\setminus \left\{\bigcup_{k=1}^{a_{33}-1}
\Omega_{d_3}^k({\bf d}^2)\right\}\supseteq \Delta({\bf d}^3)\;.
\label{sub2}
\end{eqnarray}
Let $\overline{\sigma}\in \Delta({\bf d}^3)$. Then $\overline{\sigma}\in 
\Delta({\bf d}^2)$ by (\ref{def1b}). Suppose $\overline{\sigma}\not\in
\Delta({\bf d}^2)\setminus \left\{\bigcup_{k=1}^{a_{33}-1}\Omega_{d_3}^k(
{\bf d}^2)\right\}$. Then $\overline{\sigma}\in \bigcup_{k=1}^{a_{33}-1}
\Omega_{d_3}^k({\bf d}^2)$. But then $\overline{\sigma}$ is representable 
by $d_1,d_2,d_3$ by Definition \ref{definit2} that again contradicts our 
assumption on $\overline{\sigma}$. Hence (\ref{sub2}) holds true and the 
Theorem is proved.$\;\;\;\;\;\;\Box$

\begin{figure}[h]
\centerline{\psfig{figure=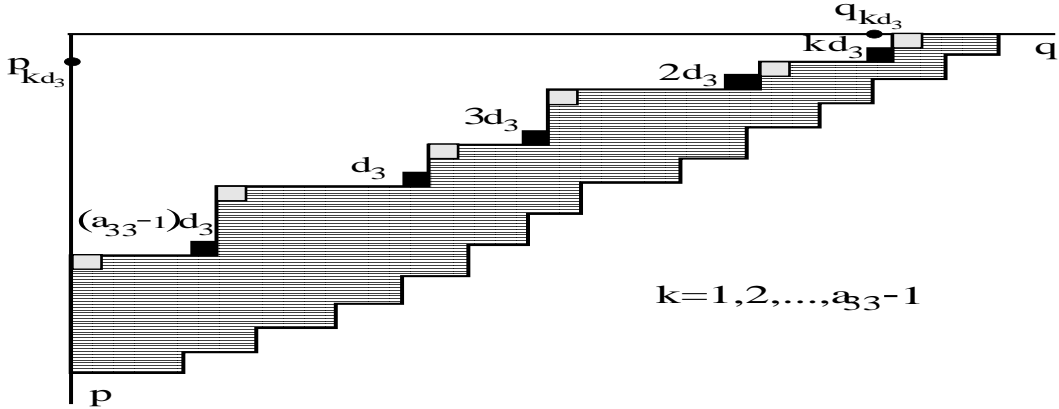,height=5.5cm,width=14cm}}
\caption{Typical matrix representation $M\left\{\Delta({\bf d}^3)\right\}$ of
the set $\Delta({\bf d}^3)$ ({\em striped area}) inside $M\left\{\Delta({\bf
d}^2)\right\}$. ${\sf TL}_M\left\{\Delta({\bf d}^3)\right\}$ has $a_{33}+1$ 
convex corners ({\em gray boxes}). The integers $kd_3\not\in\Delta({\bf d}^3),1
\leq k<a_{33}$ occupy the concave corners ({\em black boxes}) of the union 
$\bigcup_{k=1}^{a_{33}-1}\Omega_{d_3}^k({\bf d}^2)$, which are adjacent
to the concave corners of ${\sf TL}_M\left\{\Delta({\bf d}^3)\right\}$.}
\label{repr3}
\end{figure}

Generalizing ${\sf BL}_M\left\{\Delta({\bf d}^2)\right\}$ and ${\sf TL}_M\left\{
\Delta({\bf d}^2)\right\}$ on $m=3$ call the totalities of the lowest and top 
cells in every column of $M\left\{\Delta({\bf d}^3)\right\}$ {\em the bottom} 
and {\em top layers} of $M\left\{\Delta({\bf d}^3)\right\}$, respectively, with 
corresponding notations, ${\sf BL}_M\left\{\Delta({\bf d}^3)\right\}$ and 
${\sf TL}_M\left\{\Delta({\bf d}^3)\right\}$
\begin{eqnarray}
{\sf BL}_M\left\{\Delta({\bf d}^3)\right\}\subset M\left\{\Delta({\bf d}^3)
\right\}\;,\;\;\;{\sf TL}_M\left\{\Delta({\bf d}^3)\right\}\subset 
M\left\{\Delta({\bf d}^3)\right\}\;.\label{map2aa}
\end{eqnarray}
In Figure \ref{repr3} we present the typical matrix representation $M\left\{
\Delta({\bf d}^3)\right\}$ of the set $\Delta({\bf d}^3)$ inside $M\left\{
\Delta({\bf d}^2)\right\}$. The bottom layer ${\sf BL}_M\left\{\Delta({\bf d}
^3)\right\}$ of this diagram coincides with ${\sf BL}_M\left\{\Delta({\bf d}^2)
\right\}$ presented in Figure \ref{repr1} (see (\ref{map2b}) in Section 
\ref{nomina0}). The top layer ${\sf TL}_M\left\{\Delta({\bf d}^3)\right\}$ of 
this diagram is much more intricate than ${\sf TL}_M\left\{\Delta({\bf d}^2)
\right\}$ given in (\ref{sylves9}), e.g. ${\sf TL}_M\left\{\Delta({\bf 
d}^3)\right\}$ has $a_{33}+1$ convex corners. 
\section{Diagrammatic calculation on the set $\Delta\left({\bf d}^3\right)$}
\label{nomina0}
A straightforward reconstruction of the Hilbert series $H({\bf d}^3;z)$ of a 
graded ring ${\sf k}\left[z^{d_1},z^{d_2},z^{d_3}\right]$ out of the set $\Delta
\left({\bf d}^3\right)$ is a difficult problem. In order to overcome this 
difficulty we develop the procedure of {\em diagrammatic calculation} in 
$\Delta\left({\bf d}^3\right)$ in the present Section. This procedure will be 
applied in Section \ref{qqq1} to calculate $Q({\bf d}^3;z)$ and to give a 
complete solution of the 3D Frobenius problem. The diagrammatic calculation is 
also useful in higher dimensions $m\geq 4$ and enables us to estimate the upper 
bound for the number of non--zero coefficients in the polynomial $Q({\bf d}^m;
z)$ (see Section \ref{appl0}).

The algebraic approach to the Frobenius problem is based on a strong 
relationship between Hilbert series $H({\bf d}^m;z)$ of a graded ring ${\sf k}
\left[z^{d_1},\ldots,z^{d_m}\right]$ over a field of characteristic 0, and the 
generating function $\Phi\left({\bf d}^m;z\right)$ for the set $\Delta({\bf 
d}^m)$ \cite{stan96}
\begin{eqnarray}
H({\bf d}^m;z)+\Phi({\bf d}^m;z)=\frac{1}{1-z}\;,
\label{nonrep1}
\end{eqnarray}
where $\Phi({\bf d}^m;z)$ and $H({\bf d}^m;z)$ are defined in (\ref{def1}) and 
(\ref{hilb0}), respectively. Being evaluated at a special value of $z$ the 
function $\Phi\left({\bf d}^m;z\right)$ gives the Frobenius number and the 
genus of semigroups in any dimension $m$. Indeed, according to 
the definitions (\ref{defgen}) and (\ref{condfrob}) we have
\begin{eqnarray}
F\left({\bf d}^m\right)=\deg \Phi\left({\bf d}^m;z\right)\;,\;\;\;\;
G\left({\bf d}^m\right)=\sum_{s\;\in\;\Delta\left({\bf d}^m\right)} 1^s=
\Phi\left({\bf d}^m;1\right)\;.
\label{integr1ab}
\end{eqnarray}
Making use of (\ref{hilb0}) and (\ref{nonrep1}) formulas (\ref{integr1ab}) can 
be represented in more analytical way 
\begin{eqnarray}
F\left({\bf d}^m\right)&=&\deg Q({\bf d}^m;z)-\sum_{j=1}^md_j
\;,\label{integr1a}\\
G\left({\bf d}^m\right)&=&\lim_{z\to 1}\frac{\prod_{j=1}^m\left(1-z^{d_j}
\right)-(1-z)Q({\bf d}^m;z)}{(1-z)\prod_{j=1}^m\left(1-z^{d_j}\right)}=
\frac{\partial_z^{m+1}\left[\prod_{j=1}^m \left(1-z^{d_j}\right)-(1-z)
Q({\bf d}^m;z)\right]_{|\;z=1}}
{\partial_z^{m+1}\left[(1-z)\prod_{j=1}^m\left(1-z^{d_j}\right)\right]_{|\;z=1}}
\nonumber\\
&=&\frac{(-1)^{m+1}}{(m+1)!\prod_{j=1}^md_j}\partial_z^{m+1}
\left[\prod_{j=1}^m\left(1-z^{d_j}\right)-(1-z)Q({\bf d}^m;z)\right]_{|\;z=1},
\label{integr1b}
\end{eqnarray}
where $\partial_z^n={\sl d^n/dz^n}$ stands for the usual derivative of $n$th 
order. As one can see from (\ref{integr1a}) and (\ref{integr1b}), the Frobenius
problem is reduced to finding the numerator $Q({\bf d}^m;z)$ of Hilbert series
which follows if one substitutes (\ref{hilb0}) into (\ref{nonrep1})
\begin{eqnarray}
Q({\bf d}^m;z)=\prod_{j=2}^{m-1}\left(1-z^{d_j}\right)\left[\sum_{k=0}^
{d_1-1}z^k-\left(1-z^{d_1}\right)\Phi\left({\bf d}^m;z\right)\right]\;.
\label{integr1c}
\end{eqnarray}
A straightforward reconstruction of the numerator $Q\left({\bf d}^m;z\right)$ 
out of the set $\Delta\left({\bf d}^m\right)$ is a very difficult problem. The 
main difficulty arises when we are going to handle the term 
\begin{eqnarray}
\sum_{k=0}^{d_1-1}z^k-\left(1-z^{d_1}\right)\Phi\left({\bf d}^m;z\right)\;.
\label{integr1d}
\end{eqnarray}
However, it appears that in dimension $m=3$ one can elaborate an effective
procedure to calculate $Q({\bf d}^3;z)$ via geometrical transformations ({\em 
shifts}) of a diagram of the matrix representation $M\left\{\Delta\left({\bf d}^
3\right)\right\}$. We call such procedure {\em diagrammatic calculation}. It 
turns out, that diagrammatic calculation in dimensions $m=3$ reduces the 
determination of $Q({\bf d}^3;z)$ to the calculation of ${\sf TL}_M\left\{
\Delta({\bf d}^3)\right\}$ and ${\sf BL}_M\left\{\Delta({\bf d}^3)\right\}$ but 
not of the entire matrix $M\left\{\Delta\left({\bf d}^3\right)\right\}$.

First, introduce two functions, $\tau$ and its inverse $\tau^{-1}$, where $\tau$ 
maps each polynomial $\sum c_kz^k\in{\mathbb N}[z]$ with $c_k\in\{0,1\}$ onto 
the set of degrees $\{k\in {\mathbb N}\;\bracevert\;c_k\neq 0\}$. In particular, 
it follows from (\ref{def1})
\begin{eqnarray}
\tau\left[\Phi\left({\bf d}^3;z\right)\right]=\Delta\left({\bf d}^3\right)\;
\;\;\mbox{and}\;\;\;\tau^{-1}\left[\Delta\left({\bf d}^3\right)\right]=
\Phi\left({\bf d}^3;z\right)\;.
\label{map1a}
\end{eqnarray}
Observe that since all coefficients of the polynomial $\Phi\left({\bf d}^3;z
\right)$ are 1 or 0, we can uniquely reconstruct $\Phi\left({\bf d}^3;z\right)$ 
from $\Delta\left({\bf d}^3\right)$ and vice versa. In this sense $\tau$ is an 
isomorphic map. The map $\tau$ is also linear in the following sense: 

Let ${\bf d}^3$ be given and a set $\Delta\left({\bf d}^3\right)$ of all 
unrepresentable integers be related to its generating function $\Phi\left(
{\bf d}^3;z\right)$  by the isomorphic map $\tau$ defined in (\ref{map1a}). 
Let two sets $\Delta_1\left({\bf d}^3\right)$ and $\Delta_2\left({\bf d}^3
\right)$ be given such that
\begin{eqnarray}
\Delta_1\left({\bf d}^3\right),\Delta_2\left({\bf d}^3\right)\subset 
\Delta\left({\bf d}^3\right)\;,\;\;\;\Delta_1\left({\bf d}^3\right)\bigcap 
\Delta_2\left({\bf d}^3\right)=\emptyset\;.\label{map10a}
\end{eqnarray}
Then the following holds
\begin{eqnarray}
\tau^{-1}\left[\Delta_1\left({\bf d}^3\right)\bigcup \Delta_2\left({\bf 
d}^3\right)\right]=\tau^{-1}\left[\Delta_1\left({\bf d}^3\right)\right]+
\tau^{-1}\left[\Delta_2\left({\bf d}^3\right)\right]\;.
\label{map10b}
\end{eqnarray}
Recalling (\ref{schur55}) and (\ref{map1a}) we present below the relations 
between three main entities $M\left\{\Delta\left({\bf d}^3\right)\right\}$, 
$\Delta({\bf d}^3)$ and $\Phi\left({\bf d}^3;z\right)$ which are concerned with 
the integers that are unrepresentable by $d_1,d_2,d_3$. These relations are 
carried out by two maps, $\sigma$ and $\tau$,
\begin{eqnarray}
M\left\{\Delta\left({\bf d}^3\right)\right\}\stackrel{\sigma}\longrightarrow 
\Delta({\bf d}^3)\stackrel{\tau}\longleftarrow \Phi\left({\bf d}^3;z\right)\;.
\label{map1a1}
\end{eqnarray}
\subsection{Construction of the set $\tau\left[\sum_{k=0}^{d_1-1}z^k-
\left(1-z^{d_1}\right)\Phi\left({\bf d}^3;z\right)\right]$}
\label{1ststep}
Introduce {\em an upward shift operator} $\widehat {\sf U}_1$ which shifts the 
diagram of the matrix representation $M\left\{\Delta({\bf d}^3)\right\}$ one 
step upwards. We define 
\begin{eqnarray}
\widehat {\sf U}_1\;\sigma(p,q):=\sigma(p-1,q)\;.\label{map1b}  
\end{eqnarray}
Thus, by (\ref{sylves0}) $\sigma(p-1,q)=\sigma(p,q)+d_1$ and if we denote by 
$\widehat {\sf U}_1\;\Delta({\bf d}^3)$ the set of all integers $\widehat 
{\sf U}_1\;\sigma(p,q)$ such that $\sigma(p,q)\in\Delta({\bf d}^3)$ and define 
$\Delta^{\prime}\left({\bf d}^3\right)=\widehat {\sf U}_1\Delta\left({\bf d}^3
\right)$ then $\Delta^{\prime}\left(p,q\right)=\Delta\left(p-1,q\right)$ and 
\begin{eqnarray}
\Delta^{\prime}\left({\bf d}^3\right)=\widehat {\sf U}_1\;\Delta({\bf d}^3)=
\bigcup_{(p,q)\in M\left\{\Delta({\bf d}^3)\right\}}\widehat {\sf U}_1\;
\sigma(p,q)\;.\label{map1c}
\end{eqnarray}
For the determination of the term (\ref{integr1d}) via diagrammatic calculation 
we need the following results.

Let $\{d_1,d_2,d_3\}$ be the minimal generating set of ${\sf S}\left({\bf d}^3
\right)$ and let $\Phi\left({\bf d}^3;z\right)$ be the generating function for 
the set $\tau\left[\Phi\left({\bf d}^3;z\right)\right]$ of unrepresentable 
integers.
For our purpose here it is important that the construction of $M\left\{\Delta(
{\bf d}^3)\right\}$ via Theorem \ref{theo2} does not affect ${\sf BL}_M\left\{
\Delta({\bf d}^2)\right\}$ (the gray cells in Figure \ref{repr1}). Indeed, it is
clear that the $d_1-1$ integers $1,\ldots,d_1-1$ are unrepresentable by $d_1,d_
2,d_3$. Due to the first equality in (\ref{sylves9}) this leads to the important
result about the bottom layer ${\sf BL}_M\left\{\Delta({\bf d}^3)\right\}$
\begin{eqnarray}
\sigma\left\{{\sf BL}_M\left\{\tau\left[\Phi\left({\bf d}^3\right)\right]
\right\}\right\}=\sigma\left\{{\sf BL}_M\left\{\tau\left[\Phi\left({\bf d}^2
\right)\right]\right\}\right\}=\left\{1,2,\ldots,d_1-1\right\}\;.\label{map2b}
\end{eqnarray}
Consider the top layer ${\sf TL}_M\left\{\tau\left[\Phi\left({\bf d}^3\right)
\right]\right\}$.  It is given by
\begin{eqnarray}
\sigma\left\{{\sf TL}_M\left\{\tau\left[\Phi\left({\bf d}^3\right)\right]
\right\}\right\}=\left\{\sigma\left(p_{t{\sf 3}}(q),q\right)\right\}\;,\;\;\;
q=1,\ldots,d_1-1\;,\label{map2a}
\end{eqnarray}
where subscript $"t{\sf 3}"$ stands for top of $M\left\{\Delta({\bf 
d}^3)\right\}$ and $p_{t{\sf 3}}(q)$ is defined as
$$
p_{t{\sf 3}}(q)=\min\left\{1\leq p\;\bracevert\;\sigma(p,q)\in\Delta({\bf 
d}^3)\right\}\;.
$$
\begin{lemma}   
\label{lem7}
\begin{eqnarray}
\sigma\left\{{\sf TL}_M\left\{\tau\left[z^{d_1}\Phi\left({\bf d}^3\right)
\right]\right\}\right\}=\left\{\sigma\left(p_{t{\sf 3}}(q)-1,q\right)\right\}
\;,\;\;\;q=1,\ldots,d_1-1\;.\label{map3a}
\end{eqnarray}
\end{lemma}
{\sf Proof} $\;\;\;$
Consider the polynomial $z^{d_1}\Phi\left({\bf d}^3;z\right)$. By (\ref{def1}), 
(\ref{map1b}) and (\ref{map1c}) we obtain
\begin{eqnarray}
z^{d_1}\Phi\left({\bf d}^3;z\right)=\sum_{s\;\in\;\Delta\left({\bf d}^3\right)}
z^{s+d_1}=\sum_{s\;\in\;\widehat {\sf U}_1\Delta\left({\bf d}^3\right)}z^s\;.
\label{lemmmx1}
\end{eqnarray} 
Acting on it by the map $\tau$ we get
\begin{eqnarray}
\tau\left[z^{d_1}\Phi\left({\bf d}^3\right)\right]=\widehat {\sf U}_1\Delta
\left({\bf d}^3\right)=\widehat {\sf U}_1\;\tau\left[\Phi\left({\bf d}^3\right)
\right] \;.\label{lemmmx2}
\end{eqnarray}
By (\ref{map1b}), (\ref{map1c}) and (\ref{map2a}) this leads to (\ref{map3a}).
$\;\;\;\;\;\;\Box$
\begin{corollary} 
\label{corol1}
Let ${\bf d}^3$ be given, ${\bf d}^3=(d_1,d_2,d_3)$, and the 1st minimal 
relation ${\cal R}_1\left({\bf d}^3\right)$ is defined by (\ref{herznon1a}). 
Then
\begin{eqnarray}
kd_3\in \sigma\left\{{\sf TL}_M\left\{\tau\left[z^{d_1}\Phi\left({\bf 
d}^3\right)\right]\right\}\right\}\;,\;\;\;k=1,\ldots,a_{33}-1\;.
\label{ccc1}
\end{eqnarray}
\end{corollary}
{\sf Proof} $\;\;\;$Let ${\cal R}_1\left({\bf d}^3\right)$ be the 1st minimal 
relation for the given ${\bf d}^3$. Then $kd_3\not\in\Delta\left({\bf d}^3\right),
1\leq k< a_{33}$. First, consider one of such integers, $kd_3$, and show that
$kd_3-d_1\in \Delta\left({\bf d}^3\right)$. Let, by way of contradiction, 
$kd_3-d_1\not\in\Delta\left({\bf d}^3\right)$, then there exist $\rho_1,\rho_2,
\rho_3\in  {\mathbb N}\cup\{0\}$ such that 
$$
kd_3-d_1=\rho_1d_1+\rho_2d_2+\rho_3d_3\;\;\rightarrow\;\;(k-\rho_3)d_3=
(\rho_1+1)d_1+\rho_2d_2\;,\;\;1\leq k\leq a_{33}-1\;,
$$
violating the minimality of the relation ${\cal R}_1\left({\bf d}^3\right)$ 
given by (\ref{herznon1a}). Now for every $k=1,\ldots,a_{33}-1$ we have
\begin{eqnarray}
kd_3\not\in\Delta\left({\bf d}^3\right)\;\;\;\;\mbox{and}\;\;\;\;
kd_3-d_1\in \Delta\left({\bf d}^3\right)\;.
\label{ccc10}
\end{eqnarray}
Comparing (\ref{ccc10}) with (\ref{map2a}) and (\ref{map3a}) we conclude that 
the integers $kd_3-d_1,1\leq k< a_{33}$ occupy ${\sf TL}_M\left\{\tau\left[\Phi
\left({\bf d}^3\right)\right]\right\}$ while the integers $kd_3,1\leq k<a_{33}$ 
occupy ${\sf TL}_M\left\{\tau\left[z^{d_1}\Phi\left({\bf d}^3\right)\right]
\right\}$. This proves the Corollary.$\;\;\;\;\;\;\Box$

\begin{figure}[h]
\centerline{\psfig{figure=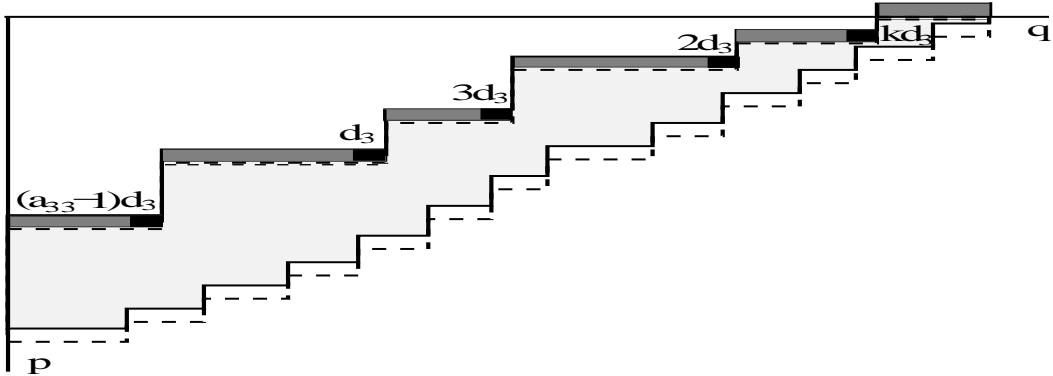,height=5cm,width=14cm}}
\caption{Matrix representation of two sets: $\tau\left[z^{d_1}\Phi\left({\bf
d}^3;z\right)\right]$ ({\em inside plain frame}) and $\tau\left[\Phi\left(
{\bf d}^3;z\right)\right]$ ({\em inside dashed frame}). Their intersection
$\Pi_1\left({\bf d}^3\right)$ is marked by {\em bright gray color}. The top
${\sf TL}_M\left\{\tau\left[z^{d_1}\Phi\left({\bf d}^3\right)\right]\right\}$
and bottom layer ${\sf BL}_M\left\{\tau\left[\Phi\left({\bf d}^3\right)\right]
\right\}$ are marked by {\em dark gray} and {\em white colors}, respectively.  
The integers $kd_3\in\sigma\left\{{\sf TL}_M\left\{\tau\left[z^{d_1}\Phi\left(
{\bf d}^3\right)\right]\right\}\right\},1\leq k<a_{33}$ are shown by {\em black
boxes}.}
\label{repr4}
\end{figure}

In Figure \ref{repr4} we show the matrix representations of two sets $\tau\left[
z^{d_1}\Phi\left({\bf d}^3;z\right)\right]$ and $\tau\left[\Phi\left({\bf d}^3;
z\right)\right]$ with their intersection 
\begin{eqnarray}
\Pi_1\left({\bf d}^3\right):=\tau\left[z^{d_1}\Phi\left({\bf d}^3;z\right)
\right]\bigcap\tau\left[\Phi\left({\bf d}^3;z\right)\right]\;.
\label{lemmmx2}
\end{eqnarray}
$M\left\{\tau\left[z^{d_1}\Phi\left({\bf d}^3;z\right)\right]\right\}$ is 
shifted one step upwards with respect to $M\left\{\tau\left[\Phi\left({\bf d}^3
;z\right)\right]\right\}$. 

From this presentation follows  
\begin{eqnarray}
\tau\left[\Phi\left({\bf d}^3;z\right)\right]&=&\Pi_1\left({\bf d}^3\right)
\bigcup\sigma\left\{{\sf BL}_M\left\{\tau\left[\Phi\left({\bf d}^3\right)\right]
\right\}\right\}\;,\label{lemmmx3a}\\
\tau\left[z^{d_1}\Phi\left({\bf d}^3;z\right)\right]&=&\Pi_1\left({\bf d}^3
\right)\bigcup\sigma\left\{{\sf TL}_M\left\{\tau\left[z^{d_1}\Phi\left(
{\bf d}^3\right)\right]\right\}\right\}\;,\;\;\;\mbox{and}\label{lemmmx3b}
\end{eqnarray}
\begin{eqnarray}
\Pi_1\left({\bf d}^3\right)\bigcap \sigma\left\{{\sf BL}_M\left\{\tau\left[\Phi
\left({\bf d}^3\right)\right]\right\}\right\}=\emptyset\;,\;\;\;\Pi_1\left(
{\bf d}^3\right)\bigcap \sigma\left\{{\sf TL}_M\left\{\tau\left[z^{d_1}\Phi
\left({\bf d}^3\right)\right]\right\}\right\}=\emptyset\;.\label{lemmmx3c}
\end{eqnarray}
Denote the integers occupying the top layer ${\sf TL}_M\left\{\tau\left[z^
{d_1}\Phi\left({\bf d}^3\right)\right]\right\}$ by $\lambda_q$. Thus, we have
\begin{eqnarray}
\sigma\left\{{\sf TL}_M\left\{\tau\left[z^{d_1}\Phi\left({\bf 
d}^3\right)\right]\right\}\right\}=\left\{\lambda_q\;\bracevert\;\lambda_q=
\sigma\left(p_{t{\sf 3}}(q)-1,q\right),\;1\leq q\leq d_1-1\right\}\;.
\label{lemmmx4b}
\end{eqnarray}
Now we are ready to prove the main Theorem of this Section
\begin{theorem}
\label{theo3}
\begin{eqnarray}
\left(1-z^{d_1}\right)\Phi\left({\bf d}^3;z\right)=\sum_{q=1}^{d_1-1}z^q-
\sum_{q=1}^{d_1-1}z^{\lambda_q}\;.\label{map4a}
\end{eqnarray}
\end{theorem}
{\sf Proof} $\;\;\;$Consider the two polynomials
\begin{eqnarray}
\Phi\left({\bf d}^3;z\right)=\tau^{-1}\left[\Delta\left({\bf d}^3\right)\right]
\;,\;\;\;z^{d_1}\Phi\left({\bf d}^3\right)=\tau^{-1}\left[\widehat{\sf U}_1
\Delta\left({\bf d}^3\right)\right]\;,\label{map4b} 
\end{eqnarray}
and construct their difference $K_1=\left(1-z^{d_1}\right)\Phi\left({\bf d}^3;z
\right)$ acting on (\ref{lemmmx3a}) and (\ref{lemmmx3b}) by $\tau^{-1}$ 
\begin{eqnarray}
K_1=\tau^{-1}\left[\Pi_1\left({\bf d}^3\right)\bigcup \sigma\left\{
{\sf BL}_M\left\{\tau\left[\Phi\left({\bf d}^3\right)\right]\right\}\right\}
\right]-\tau^{-1}\left[\Pi_1\left({\bf d}^3\right)\bigcup \sigma\left\{
{\sf TL}_M\left\{\tau\left[z^{d_1}\Phi\left({\bf d}^3\right)\right]\right\}
\right\}\right]\;.\nonumber
\end{eqnarray}
Making use of (\ref{map10a}), (\ref{map10b}) and (\ref{lemmmx3c}) we obtain
\begin{eqnarray}
K_1=\tau^{-1}\left[\sigma\left\{{\sf BL}_M\left\{\tau\left[\Phi\left(
{\bf d}^3\right)\right]\right\}\right\}\right]-
\tau^{-1}\left[\sigma\left\{{\sf TL}_M\left\{\tau\left[z^{d_1}\Phi\left(
{\bf d}^3\right)\right]\right\}\right\}\right]\;.
\label{map4c}
\end{eqnarray}
Substituting (\ref{map2b}) and (\ref{lemmmx4b}) into (\ref{map4c}) we come 
to (\ref{map4a}) that finishes the proof of the Theorem.$\;\;\;\;\;\;\Box$

\begin{figure}[h]
\centerline{\psfig{figure=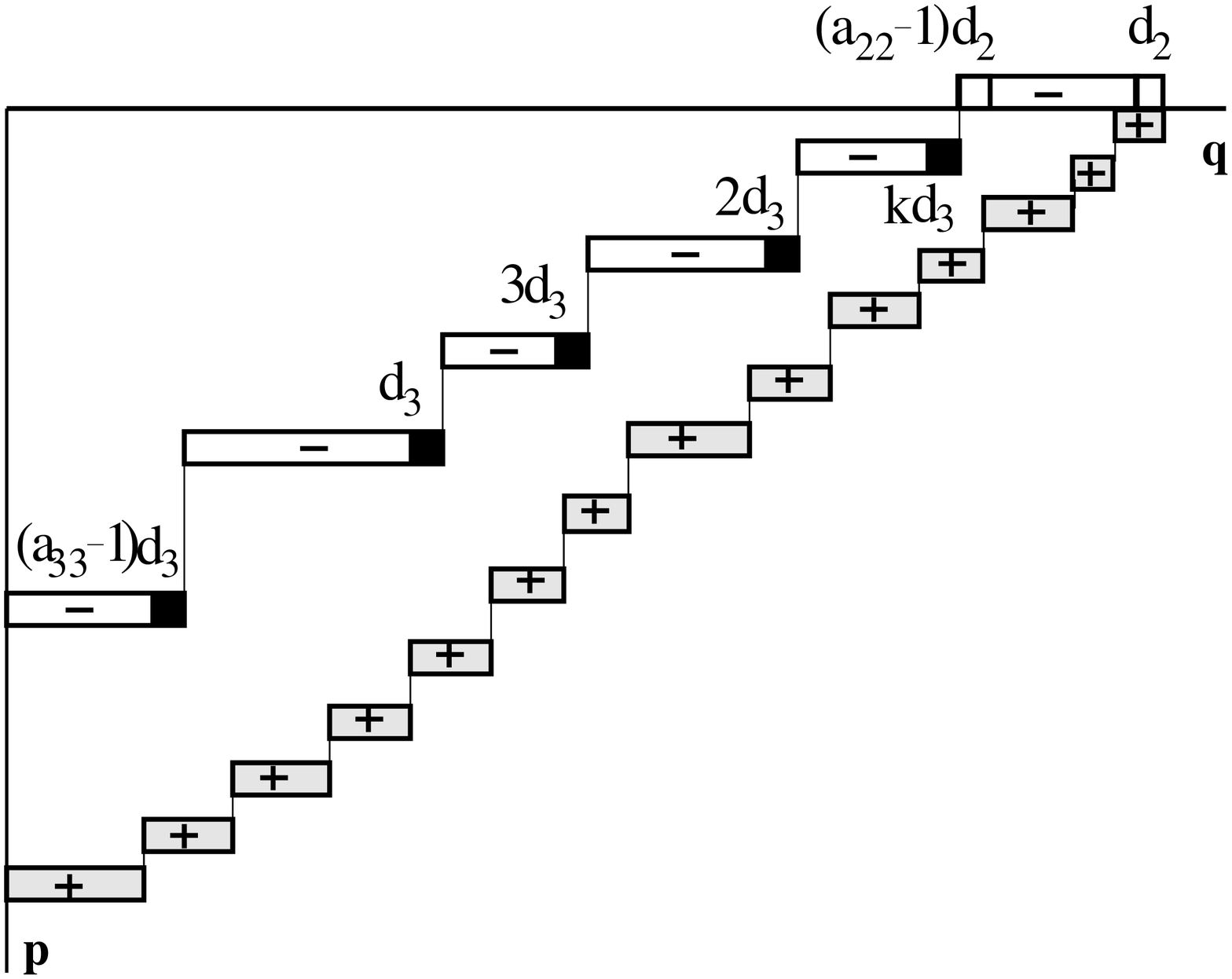,height=5cm,width=14cm}}
\caption{Matrix representation of the set $\tau\left[\left(1-z^{d_1}\right)
\Phi\left({\bf d}^3;z\right)\right]$. The signs "-" and "+" in the cells mark
the corresponding terms $-z^{\lambda_q},\lambda_q\in \sigma\left\{{\sf TL}_M
\left\{\tau\left[z^{d_1}\Phi\left({\bf d}^3\right)\right]\right\}\right\}$ and 
$z^q,q=1,\ldots,d_1-1$, respectively, which enter into polynomial (\ref{map4a}).
The integers $kd_3\not\in\Delta({\bf d}^3),1\leq k<a_{33}$ are shown by 
{\em black boxes}.}
\label{repr5} 
\end{figure}

In Figure \ref{repr5} we show the matrix representation of the set 
$\tau\left[\left(1-z^{d_1}\right)\Phi\left({\bf d}^3;z\right)\right]$.

\noindent
Finally we arrive at the term (\ref{integr1d}) which will be calculated in 
the next Theorem.
\begin{theorem}
\label{theo4}
\begin{eqnarray}
\sum_{q=0}^{d_1-1}z^q-\left(1-z^{d_1}\right)\Phi\left({\bf d}^3;z\right)=
\sum_{q=0}^{d_1-1}z^{\lambda_q}\;,\;\;\;\lambda_0=0\;.
\label{viadia1}
\end{eqnarray}
\end{theorem}
{\sf Proof} $\;\;\;$The proof follows immediately from Theorem \ref{theo3}.
$\;\;\;\;\;\;\Box$

For application in the next Section we  introduce the following notation
\begin{eqnarray}
\Lambda\left({\bf d}^3\right):=\sigma\left\{{\sf TL}_M\left\{\tau\left[z^{d_1}
\Phi\left({\bf d}^3\right)\right]\right\}\right\}\;\cup\;\left\{0\right\}\;.
\label{lemmmx4a}
\end{eqnarray}
The basic properties of the set $\Lambda\left({\bf d}^3\right)$ 
follow form (\ref{lemmmx4b}) and (\ref{lemmmx4a}) 
\begin{eqnarray}
\Lambda\left({\bf d}^3\right)\not\subset \Delta\left({\bf d}^3\right)\;,\;\;\;
\#\Lambda\left({\bf d}^3\right)=d_1\;,\;\;\;\Lambda\left({\bf d}^3\right)=
\tau\left[\sum_{q=0}^{d_1-1}z^{\lambda_q}\right]\;.\label{lemmmx4c}
\end{eqnarray}  
The structure of $\Lambda\left({\bf d}^3\right)$ is very intricate. The set 
$\Lambda\left({\bf d}^3\right)$ includes the integers $\lambda_q$  which do 
not even belong to the set $\Delta({\bf d}^2)$. Those are
\begin{eqnarray}
\overline{\lambda}_q=qd_2\;,\;\;\;1\leq q\leq a_{22}-1\;.\label{viadia2}
\end{eqnarray}
Indeed, from (\ref{sylves9}) and Lemma \ref{lem7} follows $\overline{\lambda}_q
\in\Lambda\left({\bf d}^3\right)$. On the other hand, (\ref{viadia2}) means
that $\overline{\lambda}_q$ are representable by $d_2$ and therefore 
$\overline{\lambda}_q\not\in\Delta({\bf d}^2)$. 
\subsection{The polynomial $\left(1-z^{d_2}\right)\left\{\sum_{k=0}^{d_1-1}z^k-
\left(1-z^{d_1}\right)\Phi\left({\bf d}^3;z\right)\right\}$}
\label{2ststep}
In the previous Section we have found the set $\Lambda\left({\bf d}^3\right)$ 
of integers $\lambda_q$ which contribute to the polynomial $\sum_{k=0}^{d_1-1}
z^k-\left(1-z^{d_1}\right)\Phi\left({\bf d}^3;z\right)$. Here we continue to 
construct the numerator $Q({\bf d}^3;z)$ according to (\ref{integr1c}). This 
will be done by further successive application of diagrammatic calculation on 
$\Lambda\left({\bf d}^3\right)$. 

The diagrammatic representation of the set $\Lambda\left({\bf d}^3\right)$ (see 
Figure \ref{repr5}) is not convenient to deal with. We shall start with a matrix 
representation of the integers $\lambda_q\in\Lambda\left({\bf d}^3 \right)$ 
which essentially simplifies the procedure of calculation. 
\begin{definition}
\label{definit3}
Let integers $2<d_1<d_2<d_3$ be given. Define the function $\lambda(v_2,v_3)$ as 
follows
\begin{eqnarray}
\lambda(v_2,v_3):=v_2d_2+v_3d_3\;,\;\;\;v_2,v_3\in {\mathbb N}\cup\{0\}\;.
\label{juh1}
\end{eqnarray}
\end{definition}
The next Lemma specifies the restrictions on the domain of $(v_2,v_3)$ 
introduced in Definition \ref{definit3}. This is the hardest part of the paper.
\begin{lemma}
Let ${\bf d}^3$ be given, ${\bf d}^3=(d_1,d_2,d_3)$, with the 1st minimal 
relation ${\cal R}_1\left({\bf d}^3\right)$ defined by (\ref{herznon1a}). Let 
$r$ be an integer. Then $r\in\Lambda\left({\bf d}^3\right)$ iff $r$ is uniquely 
representable as
\begin{eqnarray}
r=\lambda(v_2,v_3)\;,\;\;\;\;\mbox{where}\label{nomi44a}
\end{eqnarray}
\begin{eqnarray}
(v_2,v_3)\in \left(\left[0,a_{22}-1\right]\times \left[0,a_{13}-1\right]\right)
\cup \left(\left[0,a_{12}-1\right]\times \left[a_{13},a_{33}-1\right]\right)
\;.\label{nomi44} 
\end{eqnarray}
\label{lem8}
\end{lemma}
\vspace{-.5cm}
{\sf Proof} $\;\;\;$Observe that according to (\ref{viadia2}) and Corollary 
\ref{corol1} in Section \ref{1ststep} the following holds, respectively
\begin{eqnarray}
v_2d_2\in \Lambda\left({\bf d}^3\right)\;,\;\;0\leq v_2<a_{22}\;\;\;\;\;
\mbox{and}\;\;\;\;\;
v_3d_3\in \Lambda\left({\bf d}^3\right)\;,\;\;0\leq v_3<a_{33}\;.
\label{horiz0}
\end{eqnarray}
Fix $v_3$ such that $0\leq v_3<a_{33}$ and consider the sequence of integers 
\begin{eqnarray}
v_3d_3,\;d_2+v_3d_3,\;2d_2+v_3d_3,\ldots,\;\overline{v_2}d_2+v_3d_3\in
\Omega_{d_3}^{v_3}({\bf d}^2)\;,
\label{horiz1}
\end{eqnarray}
where the maximal element $\overline{v_2}d_2+v_3d_3$ of the sequence 
(\ref{horiz1}) is defined by 
\begin{eqnarray}
\overline{v_2}d_2+v_3d_3=\max\{v_2d_2+v_3d_3\bracevert\;v_2d_2+v_3d_3\in
\Lambda\left({\bf d}^3\right),0\leq v_2<a_{22}\}\;.
\label{horiz1a}
\end{eqnarray}
Note that the integers of the sequence (\ref{horiz1}) occupy continuously all 
the cells of the corresponding $v_3$-th horizontal row in the diagram in Figure 
\ref{repr5} (from the right to the left, without jumps).

In order to calculate $\overline{v_2}d_2+v_3d_3$ we should formulate the 
requirements it has to satisfy. They are based on two facts which follow from 
(\ref{horiz1a}).

First, according to definition (\ref{horiz1a}) the element 
$\overline{v_2}d_2+v_3d_3$ is contained in $\Lambda\left({\bf d}^3\right)$.

Second, $\overline{v_2}d_2+v_3d_3+d_2$ belongs neither to $\Lambda\left(
{\bf d}^3\right)$ (since the element $\overline{v_2}d_2+v_3d_3$ is the maximal 
in the sequence (\ref{horiz1})) nor to $\Delta\left({\bf d}^3\right)$ (since 
$\overline{v_2}d_2+v_3d_3+d_2$ is representable by $d_2,d_3$).

Recalling definition (\ref{lemmmx4a}) of the set $\Lambda\left({\bf d}^3\right)$
and Lemma \ref{lem7} we summarize the requirements as follow
\begin{eqnarray}
\overline{v_2}d_2+v_3d_3\in\Lambda\left({\bf d}^3\right)\;\;&\rightarrow&\;\;
\overline{v_2}d_2+v_3d_3-d_1\in \Delta\left({\bf d}^3\right)\;,
\label{horiz2}\\
\left\{\begin{array}{l}
(\overline{v_2}+1)d_2+v_3d_3\not\in\Lambda\left({\bf d}^3\right)\\
(\overline{v_2}+1)d_2+v_3d_3\not\in\Delta\left({\bf d}^3\right)\end{array}
\right.\;&\rightarrow&\;(\overline{v_2}+1)d_2+v_3d_3-d_1\not\in\Delta\left({\bf 
d}^3\right)\;.
\label{horiz2b}
\end{eqnarray}
The requirement (\ref{horiz2b}) provides the following representation
\begin{eqnarray}
(\overline{v_2}+1)d_2+v_3d_3-d_1=\gamma_1d_1+\gamma_2d_2+\gamma_3d_3\;\;\;
\mbox{for some}\;\;\;\gamma_1,\gamma_2,\gamma_3\in {\mathbb N}\cup \{0\}\;.
\label{horiz2c}
\end{eqnarray}
On the other hand, the integer $\overline{v_2}d_2+v_3d_3-d_1$ is not 
representable by $d_1,d_2,d_3$ due to (\ref{horiz2}). This can happen only if 
$\gamma_2=0$, otherwise $\overline{v_2}d_2+v_3d_3-d_1$ is always representable 
due to (\ref{horiz2c})
$$
\overline{v_2}d_2+v_3d_3-d_1=\gamma_1d_1+(\gamma_2-1)d_2+\gamma_3d_3\;,
\;\;\;\gamma_2\geq 1\;.
$$
Return to (\ref{horiz2c}) and consider its solution, inserting $\gamma_2=0$. 

First, consider $v_3$ in the interval $0\leq v_3<a_{13}$ and rewrite 
(\ref{horiz2c}) in the form
\begin{eqnarray}
(\overline{v_2}+1)d_2=(\gamma_1+1)d_1+(\gamma_3-v_3)d_3\;.\label{horiz4}
\end{eqnarray}
Comparing it with the 1st minimal relation $a_{22}d_2=a_{21}d_1+a_{23}d_3$ we 
get, by uniqueness (see (\ref{herznon3})), the maximal value $\overline{v_2}$
\begin{eqnarray}
\overline{v_2}=a_{22}-1\;,\;\;\;\gamma_1=a_{21}-1\;,\;\;\;\gamma_3=v_3+a_{23}\;.
\label{horiz5}
\end{eqnarray}
Thus, the first kind of the integers $\lambda(v_2,v_3)\in\Lambda\left({\bf 
d}^3\right)$ has a representation
\begin{eqnarray}
\lambda(v_2,v_3)=v_2d_2+v_3d_3\;,\;\;\mbox{where}\;\;\;\;\left(\left[0,a_{22}-
1\right]\times \left[0,a_{13}-1\right]\right)\;.\label{horiz5a}
\end{eqnarray}

Next, consider the solution of (\ref{horiz2c}) for $v_3$ in the interval $a_{13}
\leq v_3<a_{33}$.  Summation of (\ref{horiz2c}) and the 1st minimal relation 
$a_{11}d_1-a_{12}d_2-a_{13}d_3=0$ leads to the identity
$$
a_{11}d_1+(\overline{v_2}+1-a_{12})d_2+(v_3-a_{13})d_3=
(\gamma_1+1)d_1+\gamma_3d_3\;,
$$
which, by uniqueness, gives the maximal value $\overline{v_2}$
\begin{eqnarray}
\overline{v_2}=a_{12}-1\;,\;\;\;\gamma_1=a_{11}-1\;,\;\;\;
\gamma_3=v_3-a_{13}\;.\label{horiz6}
\end{eqnarray}
Thus, the second kind of the integers $\lambda(v_2,v_3)\in\Lambda\left({\bf
d}^3\right)$ has a representation
\begin{eqnarray}
\lambda(v_2,v_3)=v_2d_2+v_3d_3\;,\;\;\mbox{where}\;\;\;\;\left(\left[0,a_{12}-
1\right]\times \left[a_{13},a_{33}-1\right]\right)\;.\label{horiz6a}
\end{eqnarray}  
Combining (\ref{horiz5a}) and (\ref{horiz6a}) we come to (\ref{nomi44}).

Finally, it remains to prove the uniqueness of (\ref{nomi44a}). Indeed, 
the presentation of $\lambda(v_2,v_3)$ by (\ref{nomi44a}) is unique. A standard 
proof of uniqueness is to assume, by way of contradiction, that there are two 
such representations $\lambda(v_2^a,v_3^a)$, and $\lambda(v_2^b,v_3^b)$, 
$v_2^a\neq v_2^b,v_3^a\neq v_3^b$, and consequently,
\begin{eqnarray}
(v_2^b-v_2^a)d_2=(v_3^a-v_3^b)d_3\;.\label{horiz7j}
\end{eqnarray}
Making use of (\ref{herznon3aa}) $a_{22}a_{11}=d_3-a_{12}a_{21}\;\rightarrow
\;a_{22}<d_3$ and $a_{33}a_{11}=d_2-a_{13}a_{31}\;\rightarrow\;a_{33}<d_2$, 
and recalling the necessary constraints (\ref{horiz0}) imposed on $v_2,v_3$ we 
get
$$
|v_2^b-v_2^a|<d_3\;,\;\;\;|v_3^a-v_3^b|<d_2\;.
$$
Thus, we come to the conclusion that (\ref{horiz7j}) has no nontrivial 
solutions, since $d_2$ and $d_3$ have no common factors. This completes 
the proof of the Lemma. $\;\;\;\;\;\;\Box$

Note that due to uniqueness of the matrix representation (\ref{juh1}) of all 
integers $\lambda(v_2,v_3)\in\Lambda\left({\bf d}^3\right)$ the following 
inequalities hold for non--symmetric semigroups
\begin{eqnarray}
a_{ii}d_i\neq a_{jj}d_j\;\;\;\mbox{for}\;\;\;i\neq j\;,\;\;1\leq i,j\leq 3\;.
\label{horiz7i}
\end{eqnarray}
The case of symmetric semigroups admits only one equality in (\ref{horiz7i}) 
(see Section \ref{case1} for details). 

\begin{figure}[h]
\psfig{figure=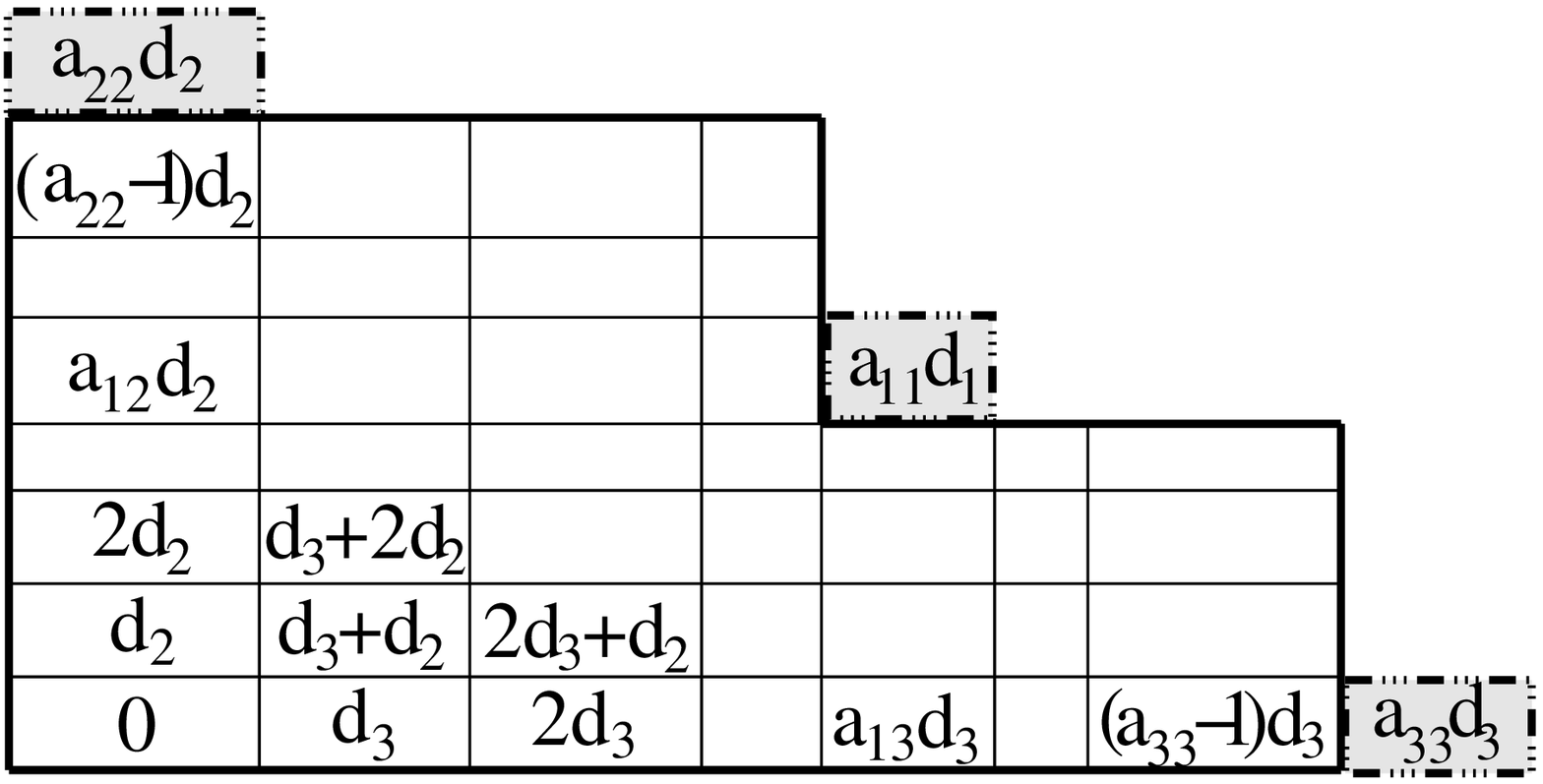,height=5cm,width=14cm}
\caption{Typical matrix representation $M\left\{\Lambda\left({\bf d}^3\right)
\right\}$ of the set $\Lambda\left({\bf d}^3\right)$. All cells in the diagram, 
which are marked in {\em white color}, are occupied by integers $\lambda(v_2,
v_3)\in\Lambda\left({\bf d}^3\right)$. The integers $a_{ii}d_i\not\in\Lambda
\left({\bf d}^3\right),i=1,2,3$, occupy three cells marked in {\em gray color}.}
\label{repr6}
\end{figure}

The representation (\ref{nomi44a}) of all integers $\lambda(v_2,v_3)\in 
\Lambda\left({\bf d}^3\right)$ is called {\em the matrix representation} 
of the set $\Lambda\left({\bf d}^3\right)$ and is denoted by $M\left\{
\Lambda\left({\bf d}^3\right)\right\}$ (see Figure \ref{repr6})
\begin{eqnarray}
\lambda\left\{M\left\{\Lambda\left({\bf d}^3\right)\right\}\right\}=
\Lambda\left({\bf d}^3\right)\;.\label{horiz7}
\end{eqnarray}
$\lambda(v_2,v_3)$ is the integer which occurs in row $v_2$ and column $v_3$ of 
$M\left\{\Lambda\left({\bf d}^3\right)\right\}$.

In  Figure \ref{repr6} we show the matrix representation 
$M\left\{\Lambda\left(
{\bf d}^3\right)\right\}$ of the set $\Lambda\left({\bf d}^3\right)$ for the 
non--symmetric semigroup ${\sf S}\left({\bf d}^3\right)$. This diagram appeared 
for the first time in \cite{schok62} for algorithmic calculation of $F\left(
{\bf d}^3\right)$. Later it was also used in \cite{frob87} and \cite{selm77} for
the same purpose. 

Lemma \ref{lem8} has also two interesting corollaries related to the diagram 
in Figures \ref{repr5} and \ref{repr6}.
\begin{corollary}
\label{corol2}
The length of horizontal rows in the matrix representation ${\sf TL}_M\left\{
\tau\left[z^{d_1}\Phi\left({\bf d}^3\right)\right]\right\}$ in Figure 
\ref{repr5}, which are covered continuously by integers $\lambda(v_2,v_3)$, 
is either $a_{22}-1$, $a_{22}$ or $a_{12}$ (defined in unit cells).
\end{corollary}
{\sf Proof} $\;\;\;$Consider the upper horizontal row of the matrix 
representation ${\sf TL}_M\left\{\tau\left[z^{d_1}\Phi\left({\bf d}^3\right)
\right]\right\}$ in Figure \ref{repr5}. According to (\ref{viadia2}) this is 
the unique row of the length $a_{22}-1$ (in unit cells) covered continuously 
by integers $\overline{\lambda}_q=qd_2,1\leq q< a_{22}$ which does not belong 
to $\Delta\left({\bf d}^2\right)$. All the other rows are contained in 
$\Delta\left({\bf d}^2\right)$ and in accordance with Corollary \ref{corol1} 
their furthest right cells are occupied by one of the integers $kd_3,1\leq k<
a_{33}$ (see Figure \ref{repr5}). Observe that all these horizontal rows are 
mapped in one to one manner into vertical columns in the diagram of the matrix 
representation $M\left\{\Lambda\left({\bf d}^3\right)\right\}$ (see Figure 
\ref{repr6}). Thus, we conclude in accordance with Lemma \ref{lem8} that their 
length is either $a_{22}$ or $a_{12}$ (defined in unit cells).$\;\;\;\;\;\;
\Box$
\begin{corollary}
\label{corol3}
Let $\{d_1,d_2,d_3\}$ be a minimal generating set of a non--symmetric 
semigroup ${\sf S}\left({\bf d}^3\right)$ and let the 1st minimal relation 
be defined by (\ref{herznon1a}). Then
\begin{eqnarray}
a_{22}+a_{33}\leq d_1+1\leq a_{22}a_{33}\;,\;\;\;
a_{33}+a_{11}\leq d_2+1\leq a_{33}a_{11}\;,\;\;\;
a_{11}+a_{22}\leq d_3+1\leq a_{11}a_{22}.\label{ari14}
\end{eqnarray}
\end{corollary}
{\sf Proof} $\;\;\;$The right hand sides of (\ref{ari14}) follow from 
(\ref{herznon3aa}): 
$$
a_{22}a_{33}=d_1+a_{23}a_{32}\geq d_1+1\;,\;
a_{33}a_{11}=d_2+a_{31}a_{13}\geq d_2+1\;,\;
a_{11}a_{22}=d_3+a_{12}a_{21}\geq d_3+1\;.
$$
The proof of the left hand sides of (\ref{ari14}) follows from 
(\ref{herznon3aaa}) and (\ref{herznon3aa}), e.g.
\begin{eqnarray}
d_1+1-(a_{22}+a_{33})=a_{22}a_{33}-a_{23}a_{32}+1-(a_{22}+a_{33})=\nonumber\\
a_{22}a_{33}-(a_{22}-a_{12})(a_{33}-a_{13})+1-(a_{22}+a_{33})=
1+a_{22}(a_{13}-1)+a_{33}(a_{12}-1)-a_{12}a_{13}\geq\nonumber\\
1+(a_{12}+1)(a_{13}-1)+(a_{13}+1)(a_{12}-1)-a_{12}a_{13}=
a_{12}a_{13}-1\geq 0\;.\nonumber
\end{eqnarray}
Thus, the Corollary is proved.$\;\;\;\;\;\;\Box$

Notice that the relations (\ref{ari14}) are survived as invariants under 
permutations of the elements $d_i$ in the generating set $\{d_1,d_2,d_3\}$, 
This is not completely obvious from the first glance since the ordering, 
$d_1<d_2<d_3$, should break such invariance.

We move on to the calculation of the polynomial $\left(1-z^{d_2}\right)\left\{
\sum_{k=0}^{d_1-1}z^k-\left(1-z^{d_1}\right)\Phi\left({\bf d}^3;z\right)
\right\}$ and apply the technique of diagrammatic calculation in the same way as
it was done in Section \ref{1ststep}. For this purpose call the totality of the 
lowest 
and top cells in every column of $M\left\{\Lambda({\bf d}^3)\right\}$ {\em the 
bottom} and {\em top layers} of $M\left\{\Lambda({\bf d}^3)\right\}$, 
respectively, and denote them ${\sf BL}_M\left\{\Lambda({\bf d}^3)\right\}$ and 
${\sf TL}_M\left\{\Lambda({\bf d}^3)\right\}$, correspondingly. 
As one can see from Figure \ref{repr6}
\begin{eqnarray}
\lambda\left\{{\sf BL}_M\left\{\Lambda({\bf d}^3)\right\}\right\}&=&
\left\{0,d_3,\ldots,(a_{33}-1)d_3\right\}\;,\;\;\label{horiz8}\\
\lambda\left\{{\sf TL}_M\left\{\Lambda({\bf d}^3)\right\}\right\}&=&
\left\{\lambda(a_{22}-1,v_3);\;0\leq v_3<a_{13}\right\}\;\cup\;
\left\{\lambda(a_{12}-1,v_3);\;a_{13}\leq v_3<a_{33}\right\}.\label{horiz8a}
\end{eqnarray}
Introduce {\em an upward shift operator} $\widehat {\sf U}_2$ which shifts the 
diagram of the matrix representation $M\left\{\Lambda({\bf d}^3)\right\}$ one 
step upwards. We define
\begin{eqnarray}
\widehat {\sf U}_2\;\lambda(v_2,v_3)=\lambda(v_2+1,v_3)\;.
\label{horiz9}
\end{eqnarray}
Thus, by (\ref{juh1}) $\lambda(v_2+1,v_3)=\lambda(v_2,v_3)+d_2$ and if we 
denote by $\widehat {\sf U}_2\;\Lambda({\bf d}^3)$ the set of all integers
$\widehat {\sf U}_2\;\lambda(v_2,v_3)$ such that $\lambda(v_2,v_3)\in\Lambda(
{\bf d}^3)$ and define $\Lambda^{\prime}\left({\bf d}^3\right)=\widehat{\sf U}_
2\Lambda\left({\bf d}^3\right)$ then $\Lambda^{\prime}\left(v_2,v_3\right)=
\Lambda\left(v_2+1,v_3\right)$ and 
\begin{eqnarray}
\Lambda^{\prime}\left({\bf d}^3\right)=\widehat {\sf U}_2\;\Lambda({\bf d}^3)=
\bigcup_{(v_2,v_3)\in M\left\{\Lambda({\bf d}^3)\right\}}\widehat {\sf U}_2\;
\lambda(v_2,v_3)\;.\label{horiz10}
\end{eqnarray}
Let $\{d_1,d_2,d_3\}$ be a minimal generating set of ${\sf S}\left({\bf d}^3
\right)$ and let $\Phi\left({\bf d}^3;z\right)$ be a generating function for the
set $\tau\left[\Phi\left({\bf d}^3;z\right)\right]$ of unrepresentable integers.
This implies, by Theorem \ref{theo4}, that $\sum_{q=0}^{d_1-1}z^q-\left(1-
z^{d_1}\right)\Phi\left({\bf d}^3;z\right)$ is a generating function for 
the set $\Lambda({\bf d}^3)$. 
\begin{lemma}
\label{lem9}
\begin{eqnarray}
\lambda\left\{{\sf TL}_M\left\{\tau\left[z^{d_2}\left(\sum_{k=0}^{d_1-1}z^k-
\left(1-z^{d_1}\right)\Phi\left({\bf d}^3;z\right)\right)\right]\right\}\right\}
=\Gamma_1({\bf d}^3)\;\cup\;\Gamma_2({\bf d}^3)\;,\label{horiz11}
\end{eqnarray}
where
\begin{eqnarray}
\Gamma_1({\bf d}^3)=\left\{\lambda(a_{22},v_3);\;0\leq v_3<a_{13}\right\},
\;\;\Gamma_2({\bf d}^3)&=&\left\{\lambda(a_{12},v_3);\;a_{13}\leq v_3<a_{33}
\right\}\;,\label{horiz12}
\end{eqnarray}
\end{lemma}
{\sf Proof} $\;\;\;$
Consider the polynomial $z^{d_2}\left(\sum_{k=0}^{d_1-1}z^k-\left(1-z^{d_1}
\right)\Phi\left({\bf d}^3;z\right)\right)$. By (\ref{viadia1}), (\ref{horiz9}) 
and (\ref{horiz10}) we obtain 
\begin{eqnarray}
z^{d_2}\left(\sum_{k=0}^{d_1-1}z^k-\left(1-z^{d_1}\right)\Phi\left({\bf d}^3;
z\right)\right)=\sum_{q=0}^{d_1-1}z^{\lambda_q+d_2}=\sum_{\lambda\in\Lambda(
{\bf d}^3)}z^{\lambda+d_2}=\sum_{\lambda\;\in\;\widehat{\sf U}_2\Lambda\left(
{\bf d}^3\right)}z^{\lambda}\;.\label{horiz13}
\end{eqnarray}
Acting on it by the map $\tau$ we get
\begin{eqnarray}
\tau\left[z^{d_2}\left(\sum_{k=0}^{d_1-1}z^k-\left(1-z^{d_1}\right)\Phi\left(
{\bf d}^3;z\right)\right)\right]=\widehat{\sf U}_2\Lambda\left({\bf d}^3\right)
\;.\label{horiz14}
\end{eqnarray}
Thus, the proof of (\ref{horiz11}) is reduced to finding a set of the integers 
$\lambda\left\{{\sf TL}_M\left\{\widehat{\sf U}_2\Lambda\left({\bf d}^3\right)
\right\}\right\}$ occupying the top layer of the matrix representation 
${\sf TL}_M\left\{\widehat{\sf U}_2\Lambda\left({\bf d}^3\right)\right\}$. 
Making successive use of (\ref{horiz9}), (\ref{horiz10}) and (\ref{horiz8a})
we obtain
\begin{eqnarray}
\lambda\left\{{\sf TL}_M\left\{\widehat{\sf U}_2\Lambda({\bf d}^3)\right\}
\right\}=\left\{\lambda(a_{22},v_3);\;0\leq v_3<a_{13}\right\}\;\cup\;\left\{
\lambda(a_{12},v_3);\;a_{13}\leq v_3<a_{33}\right\}.\label{horiz14a}
\end{eqnarray}
Introducing in accordance with (\ref{horiz12}) the notations of two 
non--intersecting sets $\Gamma_1({\bf d}^3)$ and $\Gamma_2({\bf d}^3)$, 
$\Gamma_1({\bf d}^3)\cap \Gamma_2({\bf d}^3)=\emptyset$, we arrive at the 
proof of the Lemma .$\;\;\;\;\;\;\Box$

Note that according to (\ref{horiz8}) and (\ref{horiz14a})
\begin{eqnarray}
\lambda\left\{{\sf BL}_M\left\{\Lambda({\bf d}^3)\right\}\right\}\cap 
\lambda\left\{{\sf TL}_M\left\{\widehat{\sf U}_2\Lambda({\bf d}^3)\right\}
\right\}=\emptyset\;.\label{horiz14c}
\end{eqnarray}
Denote the intersection of the sets $\Lambda\left({\bf d}^3\right)$ and 
$\widehat{\sf U}_2\Lambda\left({\bf d}^3\right)$ by $\Pi_2\left({\bf d}^3
\right)$
\begin{eqnarray}
\Pi_2\left({\bf d}^3\right):=\Lambda\left({\bf d}^3\right)\cap\widehat{\sf U}_
2\Lambda\left({\bf d}^3\right)\;.\label{horiz14b}
\end{eqnarray}
Observe that the following presentation holds :
\begin{eqnarray}
\Lambda\left({\bf d}^3\right)&=&\Pi_2\left({\bf d}^3\right)\cup \lambda\left\{
{\sf BL}_M\left\{\Lambda({\bf d}^3)\right\}\right\}\;,\label{horiz15u}\\
\widehat{\sf U}_2\Lambda\left({\bf d}^3\right)&=&\Pi_2\left({\bf d}^3\right)
\cup\lambda\left\{{\sf TL}_M\left\{\widehat{\sf U}_2\Lambda\left({\bf d}^3
\right)\right\}\right\}\;,\label{horiz15a}
\end{eqnarray} 
where
\begin{eqnarray}
\Pi_2\left({\bf d}^3\right)\cap \lambda\left\{{\sf BL}_M\left\{\Lambda({\bf d}
^3)\right\}\right\}=\emptyset\;,\;\;\;
\Pi_2\left({\bf d}^3\right)\cap \lambda\left\{{\sf TL}_M\left\{\widehat{\sf U}_
2\Lambda\left({\bf d}^3\right)\right\}\right\}=\emptyset\;.\label{horiz15b}
\end{eqnarray}
Show that (\ref{horiz14b}) necessarily follows from (\ref{horiz15u}), 
(\ref{horiz15a}) and (\ref{horiz14c}). Indeed, a straightforward 
calculation gives
\begin{eqnarray}
\Lambda\left({\bf d}^3\right)\cap \widehat{\sf U}_2\Lambda\left({\bf d}^3\right)
&=&\left(\Pi_2\left({\bf d}^3\right)\cup \lambda\left\{{\sf BL}_M\left\{
\Lambda({\bf d}^3)\right\}\right\}\right)\cap \left(\Pi_2\left({\bf d}^3
\right)\cup\lambda\left\{{\sf TL}_M\left\{\widehat{\sf U}_2\Lambda\left(
{\bf d}^3\right)\right\}\right\}\right)\nonumber\\
&=&\Pi_2\left({\bf d}^3\right)\cup \left(\lambda\left\{{\sf BL}_M\left\{
\Lambda({\bf d}^3)\right\}\right\}\cap \lambda\left\{{\sf TL}_M\left\{
\widehat{\sf U}_2\Lambda\left({\bf d}^3\right)\right\}\right\}\right)=
\Pi_2\left({\bf d}^3\right)\;.\nonumber
\end{eqnarray}
Prove the important theorem.
\begin{theorem}
\label{theo5}
\begin{eqnarray}
\left(1-z^{d_2}\right)\left(\sum_{k=0}^{d_1-1}z^k-\left(1-z^{d_1}\right)\Phi
\left({\bf d}^3;z\right)\right)=\sum_{\lambda\in\lambda\left\{{\sf BL}_M\left\{
\Lambda\left({\bf d}^3\right)\right\}\right\}}z^{\lambda}-
\sum_{\lambda\in\lambda\left\{{\sf TL}_M\left\{\widehat{\sf U}_2\Lambda
\left({\bf d}^3\right)\right\}\right\}}z^{\lambda}\;.
\label{horiz16}
\end{eqnarray}
\end{theorem}
{\sf Proof} $\;\;\;$The proof is similar to that given in Theorem \ref{theo3}. 
Consider the polynomials
\begin{eqnarray}
\sum_{k=0}^{d_1-1}z^k-\left(1-z^{d_1}\right)\Phi\left({\bf d}^3;z\right)=
\tau^{-1}\left[\Lambda\left({\bf d}^3\right)\right]\;,\;\;
z^{d_2}\left(\sum_{k=0}^{d_1-1}z^k-\left(1-z^{d_1}\right)\Phi\left({\bf d}^3;z
\right)\right)=\tau^{-1}\left[\widehat {\sf U}_2\Lambda\left({\bf d}^3\right)
\right]\nonumber
\end{eqnarray}
and construct their difference $K_2=\left(1-z^{d_2}\right)\left(\sum_{k=0}^
{d_1-1}z^k-\left(1-z^{d_1}\right)\Phi\left({\bf d}^3;z\right)\right)$ acting 
on (\ref{horiz15u}) and (\ref{horiz15a}) by $\tau^{-1}$ 
\begin{eqnarray}
K_2=\tau^{-1}\left[\Pi_2\left({\bf d}^3\right)\cup \lambda\left\{{\sf BL}_M
\left\{\Lambda\left({\bf d}^3\right)\right\}\right\}\right]-\tau^{-1}\left[\Pi_2
\left({\bf d}^3\right)\cup \lambda\left\{{\sf TL}_M\left\{\widehat{\sf U}_2\Lambda
\left({\bf d}^3\right)\right\}\right\}\right]\;.\nonumber
\end{eqnarray}
Making use of (\ref{horiz15b}) and (\ref{map10a}), (\ref{map10b}) we obtain
\begin{eqnarray}
K_2=\tau^{-1}\left[\lambda\left\{{\sf BL}_M\left\{\Lambda\left({\bf d}^3\right)
\right\}\right\}\right]-\tau^{-1}\left[\lambda\left\{{\sf TL}_M\left\{\widehat{
\sf U}_2\Lambda\left({\bf d}^3\right)\right\}\right\}\right]\;,\label{horiz17}
\end{eqnarray}
that leads to (\ref{horiz16}) in accordance with definition (\ref{map1a}) of 
the inverse map $\tau^{-1}$. $\;\;\;\;\;\;\Box$

\begin{figure}[h]
\psfig{figure=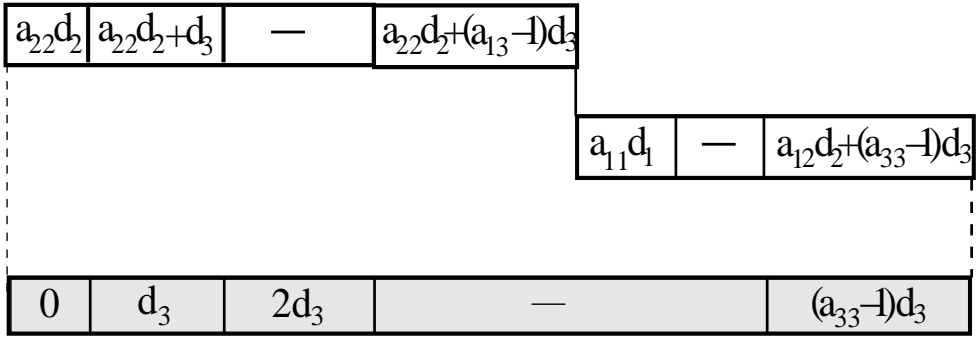,height=4.5cm,width=13cm}
\caption{Matrix representation of the set $\tau\left[\left(1-z^{d_2}\right)
\left\{\sum_{k=0}^{d_1-1}z^k-\left(1-z^{d_1}\right)\Phi\left({\bf d}^3;z\right)
\right\}\right]$ for semigroup ${\sf S}\left({\bf d}^3\right)$. Positive and 
negative contributions to polynomial (\ref{horiz16}) of the terms $z^{\lambda_
q}$ with $\lambda_q$ occupying the cells are marked in {\em gray} and {\em 
white} colors, respectively.}
\label{repr7}
\end{figure}

The result of diagrammatic calculation is shown in Figure \ref{repr7}.

\subsection{The polynomial $\left(1-z^{d_3}\right)\left(1-z^{d_2}
\right)\left\{\sum_{k=0}^{d_1-1}z^k-\left(1-z^{d_1}\right)\Phi\left({\bf d}^3;z
\right)\right\}$}
\label{3ststep}
In this Section we finish to calculate $\left(1-z^{d_3}\right)\left(1-z^{d_2}\right)
\left\{\sum_{k=0}^{d_1-1}z^k-\left(1-z^{d_1}\right)\Phi\left({\bf d}^3;z\right)
\right\}$ relying on the results obtained in Section \ref{2ststep}.

Let $\{d_1,d_2,d_3\}$ be a minimal generating set of ${\sf S}\left({\bf d}^3
\right)$ and let $\Phi\left({\bf d}^3;z\right)$ be the generating function for 
the set $\tau\left[\Phi\left({\bf d}^3;z\right)\right]$ of unrepresentable 
integers.
\begin{theorem}
\label{theo6}
\begin{eqnarray}
Q({\bf d}^3;z)=
\left(1-z^{d_3}\right)\left(1-z^{d_2}\right)\left(\sum_{k=0}^{d_1-1}z^k-
\left(1-z^{d_1}\right)\Phi\left({\bf d}^3;z\right)\right)=1-\sum_{i=1}^3z^{a_{ii}d_i}
+z^{L_1}+z^{L_2}\label{horiz18}
\end{eqnarray}
where
\begin{eqnarray}
L_1=a_{12}d_2+a_{33}d_3\;,\;\;\;\;L_2=a_{22}d_2+a_{13}d_3\;.\label{horiz19}
\end{eqnarray}
\end{theorem}
\noindent
{\sf Proof} $\;\;\;$The first equality is due to (\ref{integr1c}). By Theorem 
\ref{theo5} and Lemma \ref{lem9} we obtain
\begin{eqnarray}
\left(1-z^{d_3}\right)\left(1-z^{d_2}\right)\left(\sum_{k=0}^{d_1-1}z^k-\left(1-
z^{d_1}\right)\Phi\left({\bf d}^3;z\right)\right)=T_1-T_2-T_3\;,\label{horiz20}
\end{eqnarray}
where
\begin{eqnarray}
T_1=\left(1-z^{d_3}\right)\sum_{\lambda\;\in\;\lambda\left\{{\sf BL}_M
\left\{\Lambda\left({\bf d}^3\right)\right\}\right\}}z^{\lambda},\;\;\;
T_2=\left(1-z^{d_3}\right)\sum_{\lambda\;\in\;\Gamma_1\left({\bf d}^3\right)}
z^{\lambda},\;\;\;T_3=\left(1-z^{d_3}\right)\sum_{\lambda\;\in\;
\Gamma_2\left({\bf d}^3\right)}z^{\lambda}\;,\nonumber
\end{eqnarray}
and the sets $\lambda\left\{{\sf BL}_M\left\{\Lambda\left({\bf d}^3\right)
\right\}\right\}$ and $\Gamma_1\left({\bf d}^3\right)$, $\Gamma_2\left({\bf d}^3
\right)$ are given in (\ref{horiz8}) and (\ref{horiz12}), respectively. 
Calculating the terms $T_1,T_2$ and $T_3$ separately we get
\begin{eqnarray}
T_1&=&\left(1-z^{d_3}\right)\sum_{k=0}^{a_{33}-1}z^{kd_3}=1-z^{a_{33}d_3}\;,\;\;
T_2=\left(1-z^{d_3}\right)\sum_{k=0}^{a_{13}-1}z^{a_{22}d_2+kd_3}=z^{a_{22}d_2}
-z^{a_{22}d_2+a_{13}d_3},\nonumber\\
T_3&=&\left(1-z^{d_3}\right)\sum_{k=a_{13}}^{a_{33}-1}z^{a_{12}d_2+kd_3}=
z^{a_{11}d_1}-z^{a_{12}d_2+a_{33}d_3}\;.\label{horiz21}
\end{eqnarray}
Substituting (\ref{horiz21}) into (\ref{horiz20}) we arrive at (\ref{horiz18}).
$\;\;\;\;\;\;\Box$

\begin{figure}[h]
\centerline{\psfig{figure=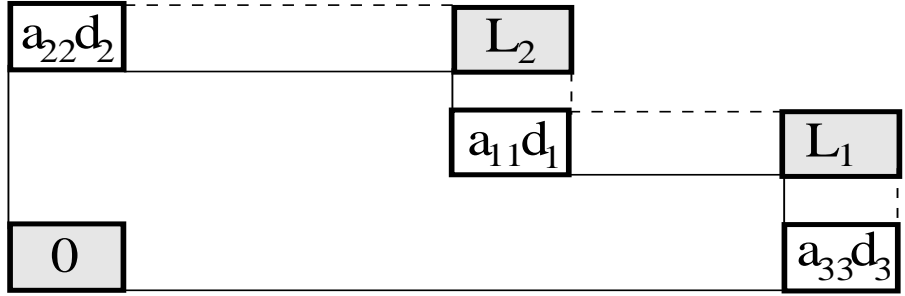,height=4cm,width=12cm}}
\caption{Matrix representation of the set $\tau\left[\left(1-z^{d_3}\right)
\left(1-z^{d_2}\right)\left\{\sum_{k=0}^{d_1-1}z^k-\left(1-z^{d_1}\right)\Phi
\left({\bf d}^3;z\right)\right\}\right]$ for non--symmetric semigroup ${\sf S}
\left({\bf d}^3\right)$. Positive and negative contributions of six terms $z^{
\lambda},\;\lambda=0,a_{11}d_1,a_{22}d_2,a_{33}d_3,L_1,L_2$ to the polynomial 
(\ref{horiz18}) are marked in {\em gray} and {\em white} colors, respectively.}
\label{repr8}
\end{figure}

The matrix representation of the set $\tau\left[\left(1-z^{d_3}\right)\left(1-z^{d_2}
\right)\left\{\sum_{k=0}^{d_1-1}z^k-\left(1-z^{d_1}\right)\Phi\left({\bf d}^3;z
\right)\right\}\right]$ for non--symmetric semigroup ${\sf S}\left({\bf d}^3\right)$
is shown at Figure \ref{repr8}.

Observe that the number of the terms contributing to (\ref{horiz18}) coincides 
with the number of corners of the polygon, which assigned the 
matrix 
representation of the set $\Lambda\left({\bf d}^3\right)$ (see Figure 
\ref{repr6}).

Below we consider two important results on the integers $L_1,L_2$ defined in 
(\ref{horiz19}). Let ${\bf d}^3$ be given, ${\bf d}^3=(d_1,d_2,d_3)$, and the 
1st minimal relation ${\cal R}_1\left({\bf d}^3\right)$ for semigroup ${\sf S}
\left({\bf d}^3\right)$ is defined by (\ref{herznon1a}). We show that 
\begin{eqnarray}
L_1\neq L_2\;.\label{lneql}
\end{eqnarray}
Assume, by way of contradiction, that the opposite is true, $L_1=L_2$. Then by 
(\ref{horiz19}) $a_{12}d_2+a_{33}d_3=a_{22}d_2+a_{13}d_3$, hence $(a_{22}-
a_{12})d_2=(a_{33}-a_{13})d_3$. By (\ref{herznon3aaa}) $a_{22}-a_{12}=a_{32}$
and $a_{33}-a_{13}=a_{23}$, hence $a_{32}d_2=a_{23}d_3$. Also by 
(\ref{herznon3aaa}) $a_{23}<a_{33}$ and $a_{32}<a_{22}$. But then $a_{32}d_2$ 
and $a_{23}d_3$ be in $\Lambda\left({\bf d}^3\right)$. Hence by the uniqueness 
of representation of elements of $\Lambda\left({\bf d}^3\right)$ (see Lemma 
\ref{lem8}) we have $a_{32}=a_{23}=0$. 
Consider the non--symmetric semigroup ${\sf S}\left({\bf d}^3\right)$ with 
$a_{23}=a_{32}=0$. The matrix $\widehat {\cal A}_{3}$ of the 1st minimal 
relation ${\cal R}_1\left({\bf d}^3\right)$ has necessarily $a_{13}=a_{12}=0$ 
that leads to $a_{11}=0$ and contradicts (\ref{herznon3}). Note that $L_1\neq 
L_2$ holds for non--symmetric and symmetric semigroups as well.

The next Lemma is related to non--symmetric semigroups ${\sf S}\left({\bf d}^3
\right)$ only.
\begin{lemma} 
\label{lem10}
\begin{eqnarray}
&&L_1\geq a_{11}d_1+d_3\;,\;\;L_1\geq a_{33}d_3+d_2\;\;\;L_1\geq 
a_{22}d_2+d_1\;,\label{lneql0}\\
&&L_2\geq a_{11}d_1+d_2\;,\;\;L_2\geq a_{33}d_3+d_1\;\;\;L_2\geq
a_{22}d_2+d_3\;.\nonumber
\end{eqnarray}
\end{lemma}
{\sf Proof} $\;\;\;$It clearly follows from diagram in Figure \ref{repr8} 
$$
L_1\geq a_{11}d_1+d_3\;,\;\;L_1\geq a_{33}d_3+d_2\;\;\;\mbox{and}\;\;\;
L_2\geq a_{11}d_1+d_2\;,\;\;L_2\geq a_{22}d_2+d_3\;.
$$
One can show that the rest two inequalities, $L_1\geq a_{22}d_2+d_1$ and $L_2
\geq a_{33}d_3+d_1$, are also true. Indeed, in accordance with (\ref{horiz19}) 
we have
$$
L_1-a_{22}d_2=a_{12}d_2+a_{33}d_3-a_{22}d_2=a_{33}d_3-a_{32}d_2=a_{31}d_1\geq 
d_1\;.
$$
In the last two equalities we used (\ref{herznon3aaa}). The last inequality, 
$L_2\geq a_{33}d_3+d_1$, can be proved in the similar manner.$\;\;\;\;\;\;\Box$

Both results, (\ref{lneql}) and Lemma \ref{lem10}, will be used later, in 
Section \ref{qqq1}.
\section{Hilbert series, Frobenius number and genus of monomial curve}
\label{qqq1}
In this Section we give a complete solution of the 3D Frobenius problem, i.e. 
calculate the numerator $Q({\bf d}^3;z)$ of Hilbert series for both 
non--symmetric and symmetric semigroups ${\sf S}\left({\bf d}^3\right)$ and on 
its basis determine the Frobenius number and genus. 
\subsection{Frobenius problem for non--symmetric semigroup ${\sf S}\left({\bf 
d}^3\right)$}
\label{case2}
Hilbert series $H\left({\bf d}^3;z\right)$, the Frobenius number $F\left({\bf d}
^3\right)$ and genus $G\left({\bf d}^3\right)$ of semigroup ${\sf S}\left(
{\bf d}^3\right)$ are invariants under permutations of elements $d_i$ in the 
generating set $\{d_1,d_2,d_3\}$. Therefore we have to find a more 
symmetrical representation for the integers $L_1$ and $L_2$ which were defined 
in (\ref{horiz19}). This leads to the main Theorem of the Section.
\begin{theorem}
\label{theo7}
Let ${\bf d}^3$ be given, ${\bf d}^3=(d_1,d_2,d_3)$, and the 1st minimal 
relation ${\cal R}_1\left({\bf d}^3\right)$ for non--symmetric semigroup 
${\sf S}\left({\bf d}^3\right)$ be defined by (\ref{herznon1a}). Then the 
numerator $Q({\bf d}^3;z)$ of Hilbert series reads
\begin{eqnarray}
Q({\bf d}^3;z)&=&1-\sum_{i=1}^3 z^{a_{ii}d_i}+
z^{1/2\left[\langle{\bf a},{\bf d}\rangle-J\left({\bf d}^3\right)\right]}+
z^{1/2\left[\langle{\bf a},{\bf d}\rangle+J\left({\bf 
d}^3\right)\right]}\;,\label{ari1}\\
J^2\left({\bf d}^3\right)&=&\langle{\bf a},{\bf d}\rangle^2-
4\sum_{i>j}^3a_{ii}a_{jj}d_id_j+4d_1d_2d_3\;,\;\;\;\langle{\bf a},{\bf d}
\rangle=\sum_{i=1}^3a_{ii}d_i\;.\label{ari1a}
\end{eqnarray}
\end{theorem}
\noindent
{\sf Proof} $\;\;\;$Making use of (\ref{herznon3aaa}) and (\ref{herznon3aa}) for
the matrix $\widehat{\cal A}_3^{(n)}$ of the 1st minimal relation for 
non--symmetric semigroup ${\sf S}\left({\bf d}^3\right)$ observe that the 
integers $L_1$ and $L_2$, defined in (\ref{horiz19}), satisfy
\begin{eqnarray}
L_1L_2=\sum_{i>j}^3a_{ii}a_{jj}d_id_j-d_1d_2d_3\;,\;\;\;L_1+L_2=\langle{\bf a},
{\bf d}\rangle\;,\;\;\;\langle{\bf a},{\bf d}\rangle=\sum_{i=1}^3a_{ii}d_i\;.
\label{ari1e}
\end{eqnarray}
In other words, $L_1$ and $L_2$ are the solutions of quadratic equation
\begin{eqnarray}
L_{1,2}^2-\langle{\bf a},{\bf d}\rangle L_{1,2}+\sum_{i>j}^3a_{ii}a_{jj}d_id_j-
d_1d_2d_3=0\;.
\label{ari1d} 
\end{eqnarray}
Therefore we have
\begin{eqnarray}
L_{1,2}=\frac{1}{2}\left[\langle{\bf a},{\bf d}\rangle\pm J\left({\bf d}^3
\right)\right]\;,\;\;J\left({\bf d}^3\right)=\sqrt{\langle{\bf a},{\bf d}
\rangle^2-4\sum_{i>j}^3a_{ii}a_{jj}d_id_j+4d_1d_2d_3}\;,\label{ari1ee}
\end{eqnarray}
that implies
\begin{eqnarray}
\langle{\bf a},{\bf d}\rangle-J\left({\bf d}^3\right)=\min\{2L_1,2L_2\}\;,\;\;\;
\langle{\bf a},{\bf d}\rangle+J\left({\bf d}^3\right)=\max\{2L_1,2L_2\}\;.
\label{ark1} 
\end{eqnarray}
Recalling the expression (\ref{integr1c}) for the numerator $Q({\bf d}^m;z)$ of 
Hilbert series and inserting (\ref{ari1e}) into (\ref{horiz18}) we come to 
(\ref{ari1}) that proves the Theorem.$\;\;\;\;\;\;\Box$
\begin{theorem}
\label{theo8}
The Frobenius number $F\left({\bf d}^3\right)$ and genus $G\left({\bf d}^3
\right)$ of non--symmetric semigroup ${\sf S}\left({\bf d}^3\right)$ read, 
respectively
\begin{eqnarray}
F\left({\bf d}^3\right)=\frac{1}{2}\left[\langle{\bf a},{\bf d}\rangle+
J\left({\bf d}^3\right)\right]-\sum_{i=1}^3d_i\;,\;\;\;\;
G\left({\bf d}^3\right)=\frac{1}{2}\left(1+\langle{\bf a},{\bf d}\rangle-
\sum_{i=1}^3d_i-\prod_{i=1}^3a_{ii}\right)\;.
\label{ari2}
\end{eqnarray}
\end{theorem}
\noindent
{\sf Proof} $\;\;\;$Implementation of formulas (\ref{integr1a})
and (\ref{integr1b}) for $Q({\bf d}^3;z)$ given by (\ref{ari1}) leads to
(\ref{ari2}). $\;\;\;\;\;\;\Box$

In Theorem \ref{theo7} the number $J\left({\bf d}^3\right)$ was considered as 
a positive integer such that the degrees $L_1$ and $L_2$ of two 
last terms in (\ref{ari1}) are positive integers. In other words, it was 
also presumed that the numbers $\langle{\bf a},{\bf d}\rangle\pm J\left({\bf d}
^3\right)$ are even positive integers. Here we are going to prove these 
statements in the case of non--symmetric semigroup ${\sf S}\left({\bf d}^3
\right)$.
\begin{lemma}
\label{lem11}
Let ${\bf d}^3$ be given, ${\bf d}^3=(d_1,d_2,d_3)$, and the 1st minimal 
relation ${\cal R}_1\left({\bf d}^3\right)$ for non--symmetric semigroup 
${\sf S}\left({\bf d}^3\right)$ be defined by (\ref{herznon1a}). Then the 
numbers $J\left({\bf d}^3\right)$ and $1/2\left[\langle{\bf a},{\bf d}\rangle
\pm J\left({\bf d}^3\right)\right]$ are non--negative and positive 
integers, respectively.
\begin{eqnarray}
J\left({\bf d}^3\right)&=&\arrowvert a_{12}a_{23}a_{31}-a_{13}a_{32}a_{21}
\arrowvert\;,\label{ari12a}\\
\frac{1}{2}\left[\langle{\bf a},{\bf d}\rangle\pm J\left({\bf d}^3\right)
\right]&=&
a_{11}a_{22}a_{33}+\frac{1}{2}\left(a_{12}a_{23}a_{31}+a_{13}a_{32}a_{21}\pm
\arrowvert a_{12}a_{23}a_{31}-a_{13}a_{32}a_{21}\arrowvert\right).\label{ari12c}
\end{eqnarray}
\end{lemma}
{\sf Proof} $\;\;\;$Inserting relations (\ref{herznon3aaa}) and 
(\ref{herznon3aa}) into the expressions $J^2\left({\bf d}^3\right)$ and 
$\langle{\bf a},{\bf d}\rangle$ given in (\ref{ari1a}) and making use of 
equality $\det \widehat {\cal A}_{3}^{(n)}=0$ which leads to
$$
a_{11}a_{22}a_{33}-(a_{11}a_{23}a_{32}+a_{22}a_{13}a_{31}+a_{33}a_{12}a_{21})
=a_{12}a_{23}a_{31}+a_{13}a_{32}a_{21}\;,
$$
we obtain
\begin{eqnarray}
J^2\left({\bf d}^3\right)=\left(a_{12}a_{23}a_{31}-a_{13}a_{32}a_{21}
\right)^2\;,\;\;\;\langle{\bf a},{\bf d}\rangle=2a_{11}a_{22}a_{33}+
a_{13}a_{32}a_{21}+a_{12}a_{23}a_{31}\;.\label{ari13}
\end{eqnarray}
The last relations lead to (\ref{ari12a}) and (\ref{ari12c}).$\;\;\;\;\;\;\Box$
\begin{corollary}
\label{corol4}
Let $\{d_1,d_2,d_3\}$ be a minimal generating set of a non--symmetric semigroup.
Then
\begin{eqnarray}
J\left({\bf d}^3\right)\geq 1\;.
\label{jneqj}
\end{eqnarray}
\end{corollary}
{\sf Proof} $\;\;\;$According to Lemma \ref{lem11} the number $J\left({\bf d}^3
\right)$ is a non--negative integer. On the other hand, due to (\ref{lneql}) and
(\ref{ari1d}) $J\left({\bf d}^3\right)$ does not vanish, that leads to 
(\ref{jneqj}).$\;\;\;\;\;\;\Box$

The unity in (\ref{jneqj}) is best possible, as the following Example shows.
\begin{example}
\label{ex1}
The triple $\{3,4,5\}$ generates a non--symmetric semigroup (with minimal 
possible $d_i$).
{\footnotesize 
\begin{eqnarray}
&&\widehat {\cal A}_3=\left(\begin{array}{rrr}3&-1&-1\\-1&2&-1\\-2&-1&2
\end{array}\right),\;\;\Delta(3,4,5)=\{1,2\}\;,\;\;\langle{\bf a},
{\bf d}\rangle=27\;,\;\;J\left({\bf d}^3\right)=1\;,\nonumber\\
&&Q\left({\bf d}^3;z\right)=1-z^8-z^9-z^{10}+z^{13}+z^{14}\;,\;\;
G({\bf d}^3)=2\;,\;\;F({\bf d}^3)=2\;.\nonumber
\end{eqnarray}}
\end{example}
Note that $d_1\geq 3$ by (\ref{assump3}) and since $\widehat{\cal A}_3$ is 
containing no zeroes, hence the semigroup is non--symmetric (see Section 
\ref{nonsymsub}).

Theorem \ref{theo8} and Lemma \ref{lem11} make it possible to express $F\left(
{\bf d}^3\right)+\sum_{i=1}^3d_i$ and $2G\left({\bf d}^3\right)+\sum_{i=1}^3
d_i$ through the elements of the matrix $\widehat{\cal A}_3^{(n)}$ of the 
1st minimal relation for a non--symmetric semigroup ${\sf S}\left({\bf d}^3
\right)$ only.
\begin{eqnarray}
F\left({\bf d}^3\right)+\sum_{i=1}^3d_i&=&a_{11}a_{22}a_{33}+\max\left\{
a_{12}a_{23}a_{31},a_{13}a_{32}a_{21}\right\}\;,\label{ari13a}\\
2G\left({\bf d}^3\right)+\sum_{i=1}^3d_i&=&1+a_{11}a_{22}a_{33}+
a_{13}a_{32}a_{21}+a_{12}a_{23}a_{31}\;.\label{ari13b}
\end{eqnarray}
Formula (\ref{ari13a}) is in full agreement with formula (\ref{herznon1}) for 
the Frobenius number obtained in \cite{herz70}, \cite{frob94}. This can be seen 
if one substitutes the relations (\ref{herznon3aaa}) and (\ref{herznon3aa}) 
into (\ref{herznon1}). Let us point out the following inequality for 
non--symmetric semigroups.
\begin{lemma} 
\label{lem12}
Let $\{d_1,d_2,d_3\}$ be the minimal generating set for non--symmetric
semigroup. Then
\begin{eqnarray}
G\left({\bf d}^3\right)\geq 1+\frac{1}{2}F\left({\bf d}^3\right)\;.
\label{winonsym}
\end{eqnarray}
\end{lemma}
{\sf Proof} $\;\;\;$Consider (\ref{ari13b}) and (\ref{ari13a}) and take 
their difference 
\begin{eqnarray}
2G\left({\bf d}^3\right)-F\left({\bf d}^3\right)&=&1+a_{13}a_{32}a_{21}+
a_{12}a_{23}a_{31}-\max\left\{a_{12}a_{23}a_{31},a_{13}a_{32}a_{21}\right\}
\nonumber\\
&=&1+\min\left\{a_{12}a_{23}a_{31},a_{13}a_{32}a_{21}\right\}\geq 2\;.
\nonumber
\end{eqnarray}
The last inequality proves the Lemma.$\;\;\;\;\;\;\Box$

Note that (\ref{winonsym}) is slightly stronger than a similar inequality 
obtained by Nijenius and Wilf \cite{wilf72} for the mD Frobenius problem.

Below we illustrate formulas (\ref{ari1}) and (\ref{ari2}) obtained for the 3D 
Frobenius problem in example for three triples, (23,29,44), (137,251,256) and 
(1563,2275,2503), which were considered numerically in \cite{barv02}, 
\cite{john60} and \cite{davis94}, respectively.
\begin{example} 
\label{ex2}
{\footnotesize
\begin{eqnarray}
&&\left(\begin{array}{c}d_1\\d_2\\d_3\end{array}\right)=\left(\begin{array}{c}
23\\29\\44\end{array}\right),\;\;\widehat {\cal A}_3=\left(\begin{array}{rrrr}
7&-1&-3\\-5&7&-2\\-2&-6&5\end{array}\right),\;\;\langle{\bf a},{\bf d}
\rangle=584\;,\;\;\;J\left({\bf d}^3\right)=86\;,\nonumber\\
&&Q({\bf d}^3;z)=1-z^{161}-z^{203}-z^{220}+z^{249}+z^{335}\;,\;\;
F\left({\bf d}^3\right)=239\;,\;\;G\left({\bf d}^3\right)=122\;.\nonumber\\
&&\left(\begin{array}{c}d_1\\d_2\\d_3\end{array}\right)=\left(\begin{array}{c} 
137\\251\\256\end{array}\right),\;\;\widehat {\cal A}_3=\left(\begin{array}
{rrrr}24&-8&-5\\-7&13&-9\\-17&-5&14\end{array}\right),\;\;\langle{\bf 
a},{\bf d}\rangle=10135\;,\;\;\;J\left({\bf d}^3\right)=1049\;,\nonumber\\
&&Q({\bf d}^3;z)=1-z^{3263}-z^{3288}-z^{3584}+z^{4543}+z^{5592}\;,\;\;
F\left({\bf d}^3\right)=4948\;,\;\;G\left({\bf d}^3\right)=2562\;.\nonumber\\
&&\left(\begin{array}{c}d_1\\d_2\\d_3\end{array}\right)=\left(\begin{array}{c}
1563\\2275\\2503\end{array}\right),\;\;\widehat {\cal A}_3=\left(\begin{array}
{rrrr}23&-7&-8\\-17&114&-93\\-6&-107&101\end{array}\right),\;\;\langle{\bf a},
{\bf d}\rangle=548102\;,\;\;\;J\left({\bf d}^3\right)=10646\;,\nonumber\\
&&Q({\bf d}^3;z)=1-z^{35949}-z^{252803}-z^{259350}+z^{268728}+z^{279374}\;,
\;\;F\left({\bf d}^3\right)=273033\;,\;\;G\left({\bf d}^3\right)=
138470\;.\nonumber
\end{eqnarray}} 
\end{example}
In the next Example we present a special kind of non-symmetric semigroups, 
{\em the Pythagorean semigroups} \cite{kraf85}. Their generators 
$d_1,d_2,d_3$ satisfy $d_1^2+d_2^2=d_3^2$.
\begin{example}
\label{ex2a}
Suppose $1\leq k_2<k_1$ such that $\gcd(k_1,k_2)=1$. Then
{\footnotesize
\begin{eqnarray}
&&\left(\begin{array}{c}d_1\\d_2\\d_3\end{array}\right)=\left(\begin{array}{c}
k_1^2-k_2^2\\2k_1k_2\\k_1^2+k_2^2\end{array}\right),\;\;\widehat {\cal 
A}_3^{Pth}=\left(\begin{array}{rrrr}k_1+k_2&-k_1+k_2&-k_1+k_2\\-k_2&k_1&-k_2
\\-k_1&-k_2&k_1\end{array}\right)\;,\;\left.\begin{array}{l}
\langle{\bf a},{\bf d}\rangle=(2k_1-k_2)(k_1+k_2)^2\\
J^{Pth}\left({\bf d}^3\right)=k_2(k_1-k_2)^2\end{array}\right.,\nonumber\\
&&Q^{Pth}\left({\bf d}^3;z\right)=1-z^{(k_1+k_2)(k_1^2-k_2^2)}-z^{2k_1^2k_2}-
z^{k_1(k_1^2+k_2^2)}+z^{k_1(k_1+k_2)^2-k_2(k_1^2+k_2^2)}+z^{k_1(k_1^2+
2k_1k_2-k_2^2)}\;,\nonumber\\
&&F^{Pth}\left({\bf d}^3\right)=k_1[k_1^2-k_2^2+2(k_1k_2-k_1-k_2)]\;,\;\;\;
G^{Pth}\left({\bf d}^3\right)=\frac{1+k_1^3-k_2^3}{2}+k_1(k_1k_2-k_1-k_2)\;.
\nonumber
\end{eqnarray}}
\end{example}
Note that the triple $\{3,4,5\}$ from Example \ref{ex1} generates the 
Pythagorean semigroup with minimal generators.
\subsection{Frobenius problem for symmetric semigroup ${\sf S}\left({\bf d}^3
\right)$}
\label{case1}
Being a special type of non--symmetric semigroup ${\sf S}\left({\bf d}^3\right)$
the case of symmetric semigroup essentially simplifies formulas (\ref{ari1}) 
and (\ref{ari2}) for $Q\left({\bf d}^3;z\right)$, $F\left({\bf d}^3\right)$ 
and $G\left({\bf d}^3\right)$.

\begin{figure}[h]
\centerline{\psfig{figure=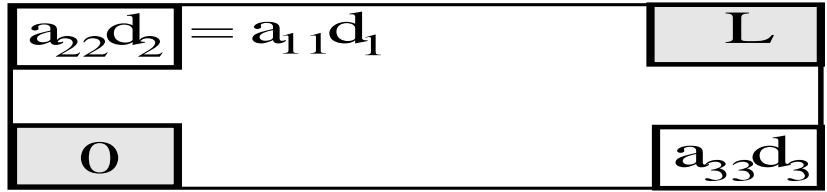,height=2.6cm,width=11cm}}
\caption{Matrix representation of the set $\tau\left[\left(1-z^{d_3}\right)
\left(1-z^{d_2}\right)\left(\sum_{k=0}^{d_1-1}z^k-\left(1-z^{d_1}\right)\Phi
\left({\bf d}^3;z\right)\right)\right]$ for symmetric semigroup ${\sf S}\left(
{\bf d}^3\right)$. Positive and negative contributions of four terms $z^{
\lambda},\;\lambda=0,a_{ii}d_i,L=a_{11}d_1+a_{33}d_3$ to numerator $Q\left(
{\bf d}^3;z\right)$ are marked by {\em gray} and {\em white} boxes,
respectively.}\label{repr9}
\end{figure}

First, a matrix representation of the set $\tau\left[\left(1-z^{d_3}\right)
\left(1-z^{d_2}\right)\left(\sum_{k=0}^{d_1-1}z^k-\left(1-z^{d_1}\right)\Phi
\left({\bf d}^3;z\right)\right)\right]$ looks much simpler (see Figure 
\ref{repr9}) and leads to the known Hilbert series (\ref{sylves10}) with four 
non--zero terms in the numerator $Q\left({\bf d}^3;z\right)$. Denote the 
Frobenius number and the genus for symmetric semigroup ${\sf S}\left({\bf d}^3
\right)$ by $F_s\left({\bf d}^3\right)$ and $G_s\left({\bf d}^3\right)$, 
respectively, and derive their expressions. The 1st minimal relation 
(\ref{herznon3d}) for the given symmetric semigroup together with 
(\ref{herznon1a}) and (\ref{herznon3aa}) yield 
\begin{eqnarray}
a_{11}d_1=a_{22}d_2\;,\;\;\;a_{22}a_{33}=d_1\;,\;\;\;a_{11}a_{33}=d_2\;,
\;\;\;\;a_{ii}\geq 2\;,\;\;i=1,2,3\;.
\label{ari2c}
\end{eqnarray}

\noindent
Substituting the relations (\ref{ari2c}) into (\ref{ari1ee}) and (\ref{ari2}) 
we obtain
\begin{eqnarray}
J_s\left({\bf d}^3\right)=a_{33}d_3\;,\;\;\;\;L_{1,2}=\frac{1}{2}\left[
2a_{11}d_1+a_{33}d_3\pm a_{33}d_3\right]\;,\nonumber\\
F_s\left({\bf d}^3\right)=a_{11}d_1+a_{33}d_3-\sum_{i=1}^3d_i\;,\;\;\;\;
G_s\left({\bf d}^3\right)=\frac{1}{2}\left[1+F_s\left({\bf d}^3\right)\right]\;.
\label{ari2a}
\end{eqnarray}
The latter formula in (\ref{ari2a}) has the following Corollary.
\begin{corollary}
\label{corol5}
Let $\{d_1,d_2,d_3\}$ be the minimal generating set for symmetric semigroup
${\sf S}\left({\bf d}^3\right)$. Then the Frobenius number $F_s\left({\bf d}^3
\right)$ is always an odd integer.
\end{corollary}
We finish this Section with an interesting observation. Recall that due to  
(\ref{assump3}) all elements $d_i$ of the minimal generating set $\{d_1,d_2,
d_3\}$ for a non--symmetric semigroup ${\sf S}\left({\bf d}^3\right)$ exceed 2. 
It appears that this restriction becomes even stronger for symmetric semigroup.
\begin{lemma}
Let $\{d_1,d_2,d_3\}$ be the minimal generating set for a symmetric semigroup 
${\sf S}\left({\bf d}^3\right)$ and the 1st minimal relation be defined by
(\ref{herznon3d}). Then all elements $d_i$ of the minimal set exceed 3.
\label{lem13}
\end{lemma}
{\sf Proof} $\;\;\;$Let ${\sf S}\left({\bf d}^3\right)$ be a symmetric semigroup
and the 1st minimal relation be defined by (\ref{herznon3d}). Then
according to (\ref{ari2c}) we have $d_1,d_2\geq 4$, otherwise the generating 
set $d_1,d_2,d_3$ would be not minimal. Inserting the  expressions (\ref{ari2c})
for $d_1,d_2$ into one of the 1st minimal relation $a_{33}d_3=a_{31}d_1+a_{32}d_2$ 
and keeping in mind $a_{31},a_{32}\geq 1$ we get
$$
a_{33}d_3=a_{31}a_{22}a_{33}+a_{32}a_{11}a_{33}\;\;\;\rightarrow\;\;\;
d_3=a_{31}a_{22}+a_{32}a_{11}\geq 4\;.
$$
Combining all restrictions $d_i\geq 4$ we come to the proof of the Lemma.
$\;\;\;\;\;\;\Box$
\begin{example}
\label{ex3}
The triple $\{4,5,6\}$ generates a symmetric semigroup (with minimal possible 
elements $d_i$).
{\footnotesize
\begin{eqnarray}
\widehat {\cal A}_3^{(s)}=\left(\begin{array}{rrr}
3 & 0 & -2 \\
-1& 2 & -1 \\
-3& 0 & 2\end{array}\right),\;\;\Delta(4,5,6)=
\{1,2,3,7\}\;,\;\;H({\bf d}^3;z)=\frac{(1-z^{10})(1-z^{12})}{(1-z^4)(1-z^5)
(1-z^6)}\;,\;\;G({\bf d}^3)=4\;,\;\;F({\bf d}^3)=7\;.\nonumber
\end{eqnarray}}
\end{example}
\subsection{Lower bounds of the Frobenius number $F\left({\bf d}^3\right)$ 
and genus $G\left({\bf d}^3\right)$}
\label{boun1}
The history of bounds for the Frobenius number $F\left({\bf d}^3\right)$ dates 
back to Schur (see Theorem A in \cite{schok62}) and has been the subject of 
intensive study for the last 30 years (see \cite{selm77}, \cite{erdo72}, 
\cite{hami98} and references therein). The subject is a very active research 
area till now. In particular, the main interest was devoted to the upper bound 
$F^{+}({\bf d}^3)$ of the Frobenius number
\footnote{Two conjectures on the upper bound $F^{+}({\bf d}^3)$ were put 
forward recently \cite{beck05}. Detailed description of the conjectures and 
their disproof will be given in Appendix \ref{appendix2}.}. Concerning the 
lower bound $F^{-}({\bf d}^3)$, in 1994, Davison \cite{davis94} obtained
\begin{eqnarray}
F\left({\bf d}^3\right)\geq F^{-}_{Dav}({\bf d}^3)\;,\;\;\;
F^{-}_{Dav}({\bf d}^3)=\sqrt{3}\sqrt{d_1d_2d_3}-\sum_{i=1}^3d_i\;,
\label{litl0}
\end{eqnarray}
where 
`{\em the constant $\sqrt{3}$ cannot be replaced by a larger value with the 
inequality remaining true for all $d_1,d_2,d_3$}' (\cite{davis94}, Theorem 2.3).
Being obtained by combinatorial means it does not distinguish between the 
triples generating the non--symmetric and {\em symmetric} semigroups. In fact, 
the lower bound of $F\left({\bf d}^3\right)$ for the set $\{d_1,d_2,d_3\}$ 
generating symmetric semigroups is stronger than (\ref{litl0}). Moreover, it 
appears that the case of non--symmetric semigroups permits also to enhance 
slightly the Davison's bound (\ref{litl0}). In order to show this we apply here 
the results of Sections \ref{case2} and \ref{case1}, and start with the lower 
bound for non--symmetric semigroups ${\sf S}\left({\bf d}^3\right)$.
\begin{lemma}
\label{lem14}  
Let $\{d_1,d_2,d_3\}$ be the minimal generating set for non--symmetric
semigroup. Then
\begin{eqnarray}
F\left({\bf d}^3\right)\geq \sqrt{3}\sqrt{d_1d_2d_3+1}-\sum_{i=1}^3d_i\;.
\label{dlow}
\end{eqnarray}
\end{lemma}
{\sf Proof} $\;\;\;$First, we find the lower bound for $\langle{\bf a},{\bf d}
\rangle$. We start with inequalities which follow from (\ref{herznon3aa})
\begin{eqnarray}
a_{11}a_{22}>d_3\;,\;\;a_{22}a_{33}>d_1\;,\;\;a_{33}a_{11}>d_2\;\;\;  
\rightarrow\;\;\;a_{11}^2a_{22}^2a_{33}^2>d_1d_2d_3\;.\label{hard0}
\end{eqnarray}
According to (\ref{hard0}) and inequality for symmetric polynomials  
\cite{hardy59} we obtain
\begin{eqnarray}
\langle{\bf a},{\bf d}\rangle\geq 3\prod_{i=1}(a_{ii}d_i)^{1/3}>
3\sqrt{d_1d_2d_3}\;.\label{hard01}
\end{eqnarray}  
Making use of (\ref{horiz7i}) we can write
$$
(a_{11}d_1-a_{22}d_2)^2+(a_{11}d_1-a_{33}d_3)^2+(a_{22}d_2-a_{33}d_3)^2\geq 
1^2+1^2+2^2=6\;,
$$
or, in other words, 
\begin{eqnarray}
\sum_{i>j}^3a_{ii}a_{jj}d_id_j\leq \frac{1}{3}\langle{\bf a},{\bf d}\rangle^2-
1\;.\label{hard1}
\end{eqnarray}
Consider the lower bound of $F\left({\bf d}^3\right)$ in 2 regions for 
$\langle{\bf a},{\bf d}\rangle$ : 
$$
1)\;\;\langle{\bf a},{\bf d}\rangle>2\sqrt{3}\sqrt{d_1d_2d_3+1}\;\;\;\;\;
\mbox{and}\;\;\;\;\;2)\;\;
\langle{\bf a},{\bf d}\rangle\leq 2\sqrt{3}\sqrt{d_1d_2d_3+1}\;.
$$
In the 1st region we immediately arrive at (\ref{dlow}) according to the 
expression (\ref{ari2}) for $F\left({\bf d}^3\right)$. Consider the 2nd region 
and observe that due to (\ref{hard1}), 
\begin{eqnarray}
J^2\left({\bf d}^3\right)=4d_1d_2d_3+\langle{\bf a},{\bf d}\rangle^2-
4\sum_{i>j}^3a_{ii}a_{jj}d_id_j\geq 4d_1d_2d_3+4-\frac{1}{3}\langle{\bf a}
,{\bf d}\rangle^2\geq 0\;.\nonumber
\end{eqnarray}
Thus, we arrive at
\begin{eqnarray}
F\left({\bf d}^3\right)+\sum_{i}^3d_i\geq \frac{1}{2}\left(\langle{\bf  a},{\bf
d}\rangle+\sqrt{4d_1d_2d_3+4-\frac{1}{3}\langle{\bf a},{\bf d}\rangle^2}
\right)\;.\label{hard2}
\end{eqnarray}
Denote $x=\langle{\bf a},{\bf d}\rangle, c=\sqrt{d_1d_2d_3+1}$ and consider a
function $f(x)=1/2(x+\sqrt{4c^2-x^2/3})$ in the interval $3\sqrt{c^2-1}<x
\leq 2\sqrt{3}c$. It is easy to find its minimum: $\min f(x)=\sqrt{3}c$ when 
$x=2\sqrt{3}c$.
Comparing this with (\ref{hard2}) we come to (\ref{dlow}) in the 2nd region.
Combining the bounds in both regions finishes the proof of the Lemma.
$\;\;\;\;\;\;\Box$

In the next Lemma we find the lower bound of $F_s\left({\bf d}^3\right)$ 
for a symmetric semigroup ${\sf S}\left({\bf d}^3\right)$.
\begin{lemma}
\label{lem15}
Let $\{d_1,d_2,d_3\}$ be the minimal generating set for symmetric semigroup
${\sf S}\left({\bf d}^3\right)$. Then
\begin{eqnarray}
F_s\left({\bf d}^3\right)\geq 2\sqrt{d_1d_2d_3}-\sum_{i=1}^3d_i\;.\label{davisym}
\end{eqnarray}
\end{lemma}
{\sf Proof} $\;\;\;$Substituting $a_{11}a_{33}=d_2$ from (\ref{ari2a}) into the
expression (\ref{ari2c}) for $F_s\left({\bf d}^3\right)$ obtain
$$
F_s\left({\bf d}^3\right)+\sum_{i=1}^3d_i=a_{11}d_1+a_{33}d_3=
\frac{d_1d_2}{a_{33}}+a_{33}d_3\geq 2\sqrt{d_1d_2d_3}\;,
$$
that proves the Lemma.$\;\;\;\;\;\;\Box$

The lower bound (\ref{davisym}) is stronger than the Davison's lower bound 
(\ref{dlow}) for $F\left({\bf d}^3\right)$ in the generic case of a 
non--symmetric semigroup ${\sf S}\left({\bf d}^3\right)$.

As for the lower bound $G^-\left({\bf d}^3\right)$ of the genus, to our 
knowledge, this question was not discussed earlier (see \cite{wilf72} though). 
Combining (\ref{winonsym}) with Lemma \ref{lem14} for non--symmetric 
semigroups and (\ref{ari2a}) with Lemma \ref{lem15} for symmetric semigroups 
one can obtain the lower bounds $G^-\left({\bf d}^3\right)$ and 
$G^-_s\left({\bf d}^3\right)$, respectively
\begin{corollary}
\label{corol6}
\begin{eqnarray}
G^-\left({\bf d}^3\right)=1+\frac{\sqrt{3}}{2}\sqrt{d_1d_2d_3+1}
-\frac{1}{2}\sum_{i=1}^3d_i\;,\;\;\;
G^-_s\left({\bf d}^3\right)=\frac{1}{2}+\sqrt{d_1d_2d_3}-
\frac{1}{2}\sum_{i=1}^3d_i\;.\label{dagennsym}
\end{eqnarray}
\end{corollary}
\section{On semigroups ${\sf S}\left({\bf d}^m\right)$ of higher dimensions, 
$m\geq 4$.}
\label{appl0}
In 1975, Bresinsky \cite{ber752} has shown that the complexity of the Frobenius 
problem changes qualitatively once $m$ exceeds 3: there exists monomial curve
in mD space, $m\geq 4$, requiring arbitrary large number of generators for 
its defining ideal ${\cal I}_m$ (see Introduction). This led Sz\'ekely 
and Wormald \cite{worm86} to the following statement,
\begin{theorem}
\label{theo9} {\rm (\cite{worm86})}
The number of non--zero coefficients in the polynomials $Q({\bf d}^m;z)$ is not 
bounded by any function of $m$ for $m\geq 4$, although it is finite for every 
choice of the generators $d_i$.
\end{theorem}
Later this Theorem was interpreted in \cite{denh03}: `{\em for any $m\geq 4$, 
there is no way to write $H({\bf d}^m;z)$ so that the polynomial $Q({\bf d}^m;z)$ 
has a bounded number of non--zero terms for all choices of $d_1,\ldots,d_m$}'. 
Making use of diagrammatic calculation developed for ${\sf S}\left({\bf d}^3
\right)$ in Section \ref{qqq1} we are going to refine the above statements
here.

Denote the number of non--zero coefficients in the polynomial $P({\bf d}^m;z)$ 
by $\#\left\{P({\bf d}^m;z)\right\}$. Thus, following (\ref{sylves1}), 
(\ref{sylves10}) and (\ref{ari1}) respectively
\begin{eqnarray}
\#\left\{Q({\bf d}^2;z)\right\}=2\;,\;\;\;
\#\left\{Q({\bf d}^3;z)\right\}=\left\{\begin{array}{l}
4\;,\;\mbox{if}\;\;{\sf S}\left({\bf d}^3\right)\;\mbox{is symmetric}\;,\\
6\;,\;\mbox{if}\;\;{\sf S}\left({\bf d}^3\right)\;\mbox{is non--symmetric}\;.
\end{array}\right.\nonumber
\end{eqnarray}
Estimate $\#\left\{Q({\bf d}^m;z)\right\}$ for non--symmetric semigroup 
${\sf S}\left({\bf d}^m\right),\;m\geq 4$. Before going to determination of an 
upper bound of $\#\left\{Q({\bf d}^m;z)\right\}$ (see Section \ref{numenu1}) we
give a brief description of basic properties of the set $\Delta({\bf d}^m)$ and
its matrix representation. 
\subsection{Basic properties of the matrix representation of the set 
$\Delta({\bf d}^m)$}
\label{basi1}
Let ${\bf d}^m$ be given, ${\bf d}^m=(d_1,\ldots,d_m)$, and let $\Delta({\bf 
d}^m)$ be the set of integers which are unrepresentable by $d_1,\ldots,d_m$, 
and let $\Phi\left({\bf d}^m;z\right)$ be a generating function for this set, 
$\tau\left[\Phi\left({\bf d}^m;z\right)\right]=\Delta\left({\bf d}^m\right)$. 
In order to construct its matrix representation $M\left\{\Delta({\bf d}^{m})
\right\}$ we have {\em to delete} from $\Delta({\bf d}^2)$ a set $\Xi\left(
{\bf d}^m\right)$ of all integers $s$ representable by $d_1,\ldots,d_m$
\begin{eqnarray}
\Delta\left({\bf d}^m\right)=\Delta({\bf d}^2)\setminus \Xi\left({\bf d}^m
\right)\;,\;\;\;\Xi\left({\bf d}^m\right)=\left\{s\;\bracevert\;s=\sum_{j=1}^m
r_jd_j,\;r_j\in{\mathbb N}\cup \{0\}\right\}\;.\label{dele1}
\end{eqnarray}
In 3D case (see Section \ref{repr3}, Theorem \ref{theo2}) this procedure was 
reduced to the construction of the complement of the union of $kd_3$--associated 
sets $\bigcup_{k=1}^{a_{33}-1}\Omega_{d_3}^k({\bf d}^2)$ in $\Delta({\bf d}^2)$.
However, in higher dimensions, $m\geq 4$, a construction of $\Delta({\bf d}^m)$ 
is not exhausted by the complement of the {\em unions} of all $kd_j$--associated
sets $\bigcup_{j=3}^m\left\{\bigcup_{k=1}^{a_{jj}-1}\Omega_{d_j}^k({\bf d}^2)
\right\}$ in $\Delta({\bf d}^2)$, where $a_{jj}$ is an uniquely defined diagonal
element of the matrix $\widehat {\cal A}_m$ of the 1st minimal relation 
${\cal R}_1\left({\bf d}^m\right)$ for given ${\bf d}^m$
\begin{eqnarray} 
\widehat {\cal A}_m\left(\begin{array}{r}
d_1\\ \ldots\\ d_j\\ \ldots \\d_m\end{array}\right)=\left(\begin{array}{c}
0\\ \ldots\\ 0\\ \ldots\\0\end{array}\right)\;,\;\;\;
\widehat {\cal A}_m=\left(\begin{array}{rrrrr}
a_{11} & \ldots & -a_{1j} & \ldots & -a_{1m} \\
\ldots & \ldots & \ldots & \ldots & \ldots \\
-a_{j1} & \ldots &a_{jj} & \ldots & -a_{jm} \\
\ldots & \ldots & \ldots & \ldots & \ldots \\ 
-a_{m1} & \ldots & -a_{mj} & \ldots & a_{mm}\end{array}\right),\label{dele2}
\end{eqnarray}
and
\begin{eqnarray}
a_{jj}=\min\left\{v_{jj}\;\bracevert\;v_{jj}\geq 2,\;v_{jj}d_j=
\sum_{i=1}^{j-1}v_{ji}d_i+\sum_{i=j+1}^mv_{ji}d_i,\;v_{ji}\in {\mathbb 
N}\cup\{0\}\right\}\;,\;\;j=1,\ldots,m\;.\label{dele3}
\end{eqnarray}
The defined values of $v_{ij},i\neq j$ which give $a_{ii}$ will be denoted by 
$a_{ij},i\neq j$. Due to minimality of the set $(d_1,\ldots,d_m)$ the elements 
$a_{ij}$ satisfy $\gcd(a_{j1},\ldots,a_{jm})=1,1\leq j\leq m$.

The reason of the complexity of the mD Frobenius problem in higher dimensions, 
$m\geq 4$, relies on the fact that there can appear \cite{kraf85} {\em 
additional minimal relations} ${\cal R}_n\left({\bf d}^m\right),n\geq 2$, 
which are linearly independent. The problem is also complicated due to the 
reason that the off--diagonal matrix elements $a_{ij},i\neq j$ are not 
necessarily unique (see Example \ref{ex5} in Section \ref{numenu1}). We omit 
here the discussion of these properties which are unimportant for further 
consideration.

All this makes the construction of the matrix representation $M\left\{\Delta(
{\bf d}^{m})\right\},m\geq 4$, extremely difficult and therefore such 
construction will not be a subject of the present paper. Nevertheless we are 
in a position to get some positive answer to the question about the number of 
non--zero coefficients in the polynomials $Q({\bf d}^m;z)$.

\begin{figure}[h]
\centerline{\psfig{figure=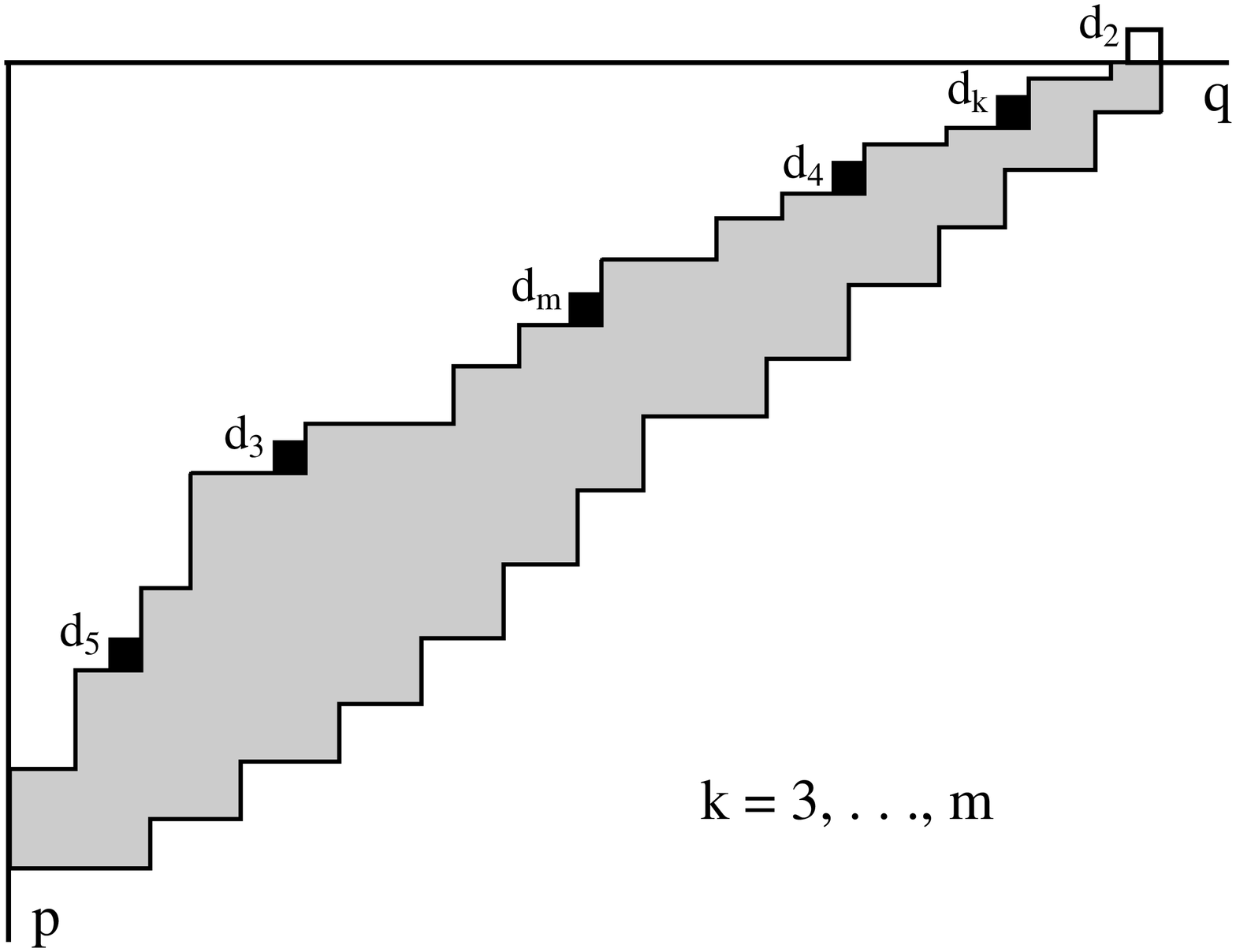,height=5.5cm,width=14cm}}
\caption{Typical matrix representation $M\left\{\Delta({\bf d}^m)\right\}$ of
the set $\Delta({\bf d}^m)$ ({\em gray color}) inside $M\left\{\Delta(
{\bf d}^2)\right\}$. The integers $d_k\not\in\Delta({\bf d}^m),\;3\leq k\leq m$
({\em black boxes}) and the integer $d_2\not\in\Delta({\bf d}^2)$ ({\em white 
box}) are adjacent to the top layer ${\sf TL}_M\left\{\Delta({\bf d}^m)
\right\}$.}
\label{repr10}
\end{figure}

Following (\ref{def1b}), recall the first containment of the sets $\Delta(
{\bf d}^{m})\subset\Delta({\bf d}^2)$ and construct the diagram of matrix 
representation $M\left\{\Delta({\bf d}^{m})\right\}$ on the basis of 
$M\left\{\Delta(d_1,d_2)\right\}$ by {\em deletion procedure} described in 
(\ref{dele1}). The obtained representation $M\left\{\Delta({\bf d}^m)\right\}$ 
(see Figure \ref{repr10}) is similar to $M\left\{\Delta({\bf d}^3)\right\}$ in 
Figure \ref{repr3} and has two common features which are important to us. 

Before discussing this let us generalize ${\sf BL}_M\left\{\Delta({\bf d}^2)
\right\}$ and ${\sf TL}_M\left\{\Delta({\bf d}^2)\right\}$ for $m\geq 4$. 
Call the totalities of the lowest and top cells in every column of $M\left\{
\Delta({\bf d}^m)\right\}$ {\em the bottom} and {\em top layers} of $M\left\{
\Delta({\bf d}^m)\right\}$, respectively, with the corresponding notations, 
${\sf BL}_M\left\{\Delta({\bf d}^m)\right\}$ and ${\sf TL}_M\left\{\Delta(
{\bf d}^m)\right\}$. 

\noindent
We also preserve the definition (\ref{map1b}) of an upward shift operator 
$\widehat {\sf U}_1$ by its action on the matrix representation of the set 
$\Delta({\bf d}^m)$: $\widehat {\sf U}_1\;\Delta({\bf d}^m)=\bigcup_{(p,q)\in 
M\left\{\Delta({\bf d}^m)\right\}}\widehat{\sf U}_1\;\sigma(p,q)$.

First, it is clear that the $d_1-1$ integers $1,\ldots,d_1-1$ are 
unrepresentable by $d_1,\ldots,d_m$ and therefore, in accordance with 
(\ref{map2b}), we have 
\begin{eqnarray}
\sigma\left\{{\sf BL}_M\left\{\tau\left[\Phi\left({\bf d}^m\right)\right]
\right\}\right\}=\left\{1,\ldots,d_1-1\right\}=\sigma\left\{{\sf BL}_M\left\{
\Delta({\bf d}^2)\right\}\right\}\;.\label{nomi17a}
\end{eqnarray}
Second, consider the top layer ${\sf TL}_M\left\{\tau\left[\Phi\left({\bf 
d}^m\right)\right]\right\}$.  It is given by
\begin{eqnarray}
\sigma\left\{{\sf TL}_M\left\{\tau\left[\Phi\left({\bf d}^m\right)\right]
\right\}\right\}=\left\{\sigma\left(p_{t{\sf m}}(q),q\right)\right\}\;,\;\;\;
q=1,\ldots,d_1-1\;,\label{nomi17e}
\end{eqnarray}
where the subscript $"t{\sf m}"$ stands for the top of $M\left\{\Delta({\bf 
d}^m)\right\}$ and $p_{t{\sf m}}(q)$ is defined as
$$
p_{t{\sf m}}(q)=\min\left\{1\leq p\;\bracevert\;\sigma(p,q)\in\Delta({\bf d}^m)
\right\}\;.
$$
Making use of an upward shift operator $\widehat {\sf U}_1$ which shifts the 
diagram of the matrix representation $M\left\{\Delta({\bf d}^m)\right\}$ one 
step upwards (\ref{map1b}) it is easy to generalize Lemma \ref{lem7} for 
$m\geq 4$
\begin{eqnarray}
\sigma\left\{{\sf TL}_M\left\{\tau\left[z^{d_1}\Phi\left({\bf d}^m\right)\right]
\right\}\right\}=\left\{\sigma\left(p_{t{\sf m}}(q)-1,q\right)\right\}\;,\;\;\;
q=1,\ldots,d_1-1\;.\label{nomi17c}
\end{eqnarray}
Apply diagrammatic calculation described in Section \ref{nomina0} in order 
to obtain the matrix representation of the set $\tau\left[\left(1-z^{d_1}
\right)\Phi\left({\bf d}^m;z\right)\right]$. In full analogue with Theorem 
\ref{theo3} its corresponding generating function looks like
\begin{eqnarray}
\left(1-z^{d_1}\right)\Phi\left({\bf d}^m;z\right)&=&\tau^{-1}\left[\sigma
\left\{{\sf BL}_M\left\{\tau\left[\Phi\left({\bf d}^m\right)\right]\right\}
\right\}\right]-\tau^{-1}\left[\sigma\left\{{\sf TL}_M\left\{\tau\left[z^{d_1}
\Phi\left({\bf d}^m\right)\right]\right\}\right\}\right]\nonumber\\
&=&\sum_{q=1}^{d_1-1}z^q-\sum_{q=1}^{d_1-1}z^{\lambda_q}\;,\;\;\;\;\;
\lambda_q\in \sigma\left\{{\sf TL}_M\left\{\tau\left[z^{d_1}\Phi\left(
{\bf d}^m\right)\right]\right\}\right\}\;,\label{nomi18}
\end{eqnarray}
where we have used (\ref{nomi17a}) and (\ref{nomi17c}). Denoting by $\Lambda
\left({\bf d}^m\right)$ a set
$$
\Lambda\left({\bf d}^m\right)=\sigma\left\{{\sf TL}_M\left\{\tau\left[z^{d_1}
\Phi\left({\bf d}^m\right)\right]\right\}\right\}\;\cup\;\left\{0\right\}\;,
$$
we finally arrive at
\begin{eqnarray}
\sum_{k=0}^{d_1-1}z^k-\left(1-z^{d_1}\right)\Phi\left({\bf d}^m;z\right)=
\sum_{q=0}^{d_1-1}z^{\lambda_q}\;,\;\;\;\lambda_0=1\;,\;\;\;\lambda_q\in
\Lambda\left({\bf d}^m\right)\;.\label{nomi17}
\end{eqnarray}
which follows when one substitutes (\ref{nomi18}) into the left hand side of 
(\ref{nomi17}). The structure of the set $\Lambda\left({\bf d}^m\right),m\geq 
4$, is built in a much more sophisticated manner than $\Lambda\left({\bf d}^3
\right)$, e.g. the matrix representation of $\Lambda\left({\bf d}^m\right)$ is 
not reduced to the simple form (\ref{juh1}). Nevertheless, one can prove the 
following general statement which will be of high importance for the estimation 
of the upper bound of $\#\left\{Q({\bf d}^m;z)\right\}$.
\begin{lemma}
\label{lem16}
Let ${\bf d}^m$ be given, ${\bf d}^m=(d_1,\ldots,d_m)$, and the 1st minimal 
relation ${\cal R}_1\left({\bf d}^m\right)$ be defined by (\ref{dele2}) 
and (\ref{dele3}). Then 
\begin{eqnarray}
kd_j\in \sigma\left\{{\sf TL}_M\left\{\tau\left[z^{d_1}\Phi\left({\bf d}^m
\right)\right]\right\}\right\}\;,\;\;\;k=1,\ldots,a_{jj}-1\;,\;\;\;
j=2,\ldots,m\;.\label{nomi18a}
\end{eqnarray}
\end{lemma}
{\sf Proof} $\;\;\;$Let ${\cal R}_1\left({\bf d}^m\right)$ be the 1st minimal 
relation defined by (\ref{dele2}) and (\ref{dele3}). Then $kd_j\not\in\Delta
\left({\bf d}^m\right)$ where $j=2,\ldots,m$ and $1\leq k<a_{jj}$. Consider one 
of such integers $kd_j$. Let, by way of contradiction, $kd_j-d_1\not\in\Delta
\left({\bf d}^m\right)$, then there exist $\rho_1,\ldots,\rho_m\in{\mathbb N}
\cup\{0\}$ such that
$$
kd_j-d_1=\sum_{i=1}^m\rho_id_i\;\;\rightarrow\;\;(k-\rho_j)d_j=(\rho_1+1)d_1+
\sum_{i=2}^{j-1}\rho_id_i+\sum_{i=j+1}^{m}\rho_id_i\;,\;\;1\leq k\leq a_{jj}-1\;,
$$
that violates the minimality of the 1st relation ${\cal R}_1\left({\bf d}^m
\right)$ given by (\ref{dele2}) and (\ref{dele3}). Hence, for every $k=1,
\ldots,a_{jj}-1$ and $j=2,\ldots,m$ we have the following pair of relations
\begin{eqnarray}
kd_j\not\in\Delta\left({\bf d}^m\right)\;\;\;\;\mbox{and}\;\;\;\;
kd_j-d_1\in \Delta\left({\bf d}^m\right)\;.
\label{nomi18b}
\end{eqnarray}
Comparing (\ref{nomi18b}) with (\ref{nomi17e}) and (\ref{nomi17c}) we conclude 
that the integers $kd_j-d_1,1\leq k<a_{jj},j=2,\ldots,m$ occupy ${\sf TL}_M
\left\{\tau\left[\Phi\left({\bf d}^m\right)\right]\right\}$ while the integers 
$kd_j,1\leq k<a_{jj},j=2,\ldots,m$ occupy ${\sf TL}_M\left\{\tau\left[z^{d_1}
\Phi\left({\bf d}^m\right)\right]\right\}$. This proves the Lemma.
$\;\;\;\;\;\;\Box$


\subsection{Upper bound for the number of non--zero coefficients in $Q({\bf d}
^m;z)$}
\label{numenu1}
Now we are ready to prove the main Theorem of this Section.
\begin{theorem}
\label{theo10} 
The number of non--zero coefficients in the polynomial $Q({\bf d}^m;z),m\geq 4$
is bounded 
\begin{eqnarray}
\#\left\{Q({\bf d}^m;z)\right\}\leq 2^{m-1}\left(d_1-\sum_{j=2}^m(a_{jj}-2)
\right)-2(m-1)\;.\label{nomi19}
\end{eqnarray}
\end{theorem}
\noindent
{\sf Proof} $\;\;\;$We will prove the Theorem in several steps. 

First, consider the expression (\ref{integr1c}) for $Q({\bf d}^m;z)$ and take 
into account (\ref{nomi17}) which implies $\#\Lambda\left({\bf d}^m\right)=
d_1$. By assumption, that a successive multiplication in (\ref{integr1c}) does 
not lead to the partial cancellation of the terms, we can get $2^{m-1}d_1$ 
non--zero terms contributing to $Q({\bf d}^m)$ that gives the first 
preliminary bound
\begin{eqnarray}
\#\left\{Q({\bf d}^m;z)\right\}\leq 2^{m-1}d_1\;.\label{nomi19a}
\end{eqnarray}

Next, if some of the diagonal elements $a_{ii}$ of the matrix $\widehat{\cal A}
_m$ of the 1st minimal relation ${\cal R}_1\left({\bf d}^m\right)$ exceeds 2, 
this bound (\ref{nomi19a}) can be actually enhanced. According to Lemma
\ref{lem16}, the polynomial (\ref{nomi17}) can be presented as follows 
\begin{eqnarray}
\sum_{k=0}^{d_1-1}z^k-\left(1-z^{d_1}\right)\Phi\left({\bf d}^m;z\right)=1+
R_1({\bf d}^m;z)+R_2({\bf d}^m;z)\;,\label{nomi19x}
\end{eqnarray}
where 
\begin{eqnarray}
R_1({\bf d}^m;z)=\sum_{k=1}^{a_{22}-1}z^{kd_2}+\ldots+\sum_{k=1}^{a_{mm}-1}
z^{kd_m}\;,\;\;\;R_2({\bf d}^m;z)=\sum_{q=1}^Nz^{\lambda_q}\;,\;\;\;
N=d_1-1-\sum_{j=2}^m(a_{jj}-1).\nonumber
\end{eqnarray}
$N$ exponents $\lambda_q$, contributing to the term $R_2({\bf d}^m;z)$, 
do not have a simple representation $\lambda_q=kd_j,1\leq k<a_{jj},2\leq j
\leq m$. Denote by $Q_1({\bf d}^m;z)$ the following part of the numerator 
$Q({\bf d}^m;z)$
\begin{eqnarray}
Q_1({\bf d}^m;z)=R_1({\bf d}^m;z)\prod_{j=2}^{m}\left(1-z^{d_j}\right)\;.
\label{nomi19b}
\end{eqnarray}
A straightforward calculation in (\ref{nomi19b}) gives
\begin{eqnarray}
Q_1({\bf d}^m;z)=\left(z^{d_2}-z^{a_{22}d_2}\right)\prod_{j=3}^{m}
\left(1-z^{d_j}\right)+\ldots+\left(z^{d_m}-z^{a_{mm}d_m}\right)
\prod_{j=2}^{m-1}\left(1-z^{d_j}\right)\;.\label{nomi19c}
\end{eqnarray}
Comparing the number of non--zero terms on the right hand sides of 
(\ref{nomi19b}) and (\ref{nomi19c}) we come to the conclusion that the entire 
number (\ref{nomi19a}) of non--zero terms of the numerator $Q({\bf d}^m;z)$ 
can be diminished by $2^{m-1}\left(\sum_{j=2}^ma_{jj}-(m-1)\right)-2^{m-1} 
(m-1)=2^{m-1}\left(\sum_{j=2}^m(a_{jj}-2)\right)$ that gives the second 
preliminary bound
\begin{eqnarray}
\#\left\{Q({\bf d}^m;z)\right\}\leq 2^{m-1}\left(d_1-\sum_{j=2}^m(a_{jj}-2)
\right)\;.\label{nomi19e}
\end{eqnarray}

Finally, return to (\ref{nomi19x}) and consider the partial cancellation of the 
terms in the polynomial $\left(1+R_1({\bf d}^m;z)\right)\prod_{j=2}^{m}\left(1-
z^{d_j}\right)$. Due to (\ref{nomi19c}) one can establish at least $m-1$ such 
terms which appear twice with different signs. Namely, these are $\sum_{j=2}^m
z^{d_j}$. Thus, one can diminish the second preliminary bound (\ref{nomi19e})
by $2(m-1)$. This proves the Theorem.$\;\;\;\;\;\;\Box$

Note that, independently of the structure of the matrix $\widehat {\cal A}_m$ 
of the 1st minimal relation ${\cal R}_1\left({\bf d}^m\right)$, the following 
always holds
\begin{eqnarray}
\#\left\{Q({\bf d}^m;z)\right\}\leq 2^{m-1}d_1-2(m-1)\;.
\label{nomi19k}
\end{eqnarray}
Theorem \ref{theo10} leads to the following restriction on the diagonal 
elements $a_{ii}$ of a matrix $\widehat{\cal A}_m$ of the 1st minimal relation 
${\cal R}_1\left({\bf d}^m\right)$.
\begin{corollary}
\label{corol7}
\begin{eqnarray}
\sum_{j=2}^ma_{jj}\leq 
d_1+2(m-1)\left(1-\frac{1}{2^{m-1}}\right)\;.\label{nomi19d}
\end{eqnarray}
\end{corollary}
{\sf Proof} $\;\;\;$The proof follows immediately from the fact that the right 
hand side in the formula (\ref{nomi19}) is positive.$\;\;\;\;\;\;\Box$

It is interesting to compare (\ref{nomi19d}) for $m=3$ with the left hand side 
of the inequality (\ref{ari14}) which is definitely stronger.

We finish this Section with one more observation about the restrictions 
imposed by the dimension of a non--symmetric semigroup ${\sf S}\left({\bf 
d}^m\right)$ on the minimal generating set $\{d_1,\ldots,d_m\}$. 
\begin{theorem} 
\label{theo11}
Let $\{d_1,\ldots,d_m\}$ be a minimal set generating semigroup ${\sf S}\left(
{\bf d}^m\right)$ such that $d_1<\ldots<d_m$, and let the 1st minimal relation 
${\cal R}_1\left({\bf d}^m\right)$ be defined by (\ref{dele2}) and 
(\ref{dele3}). Then the minimal element $d_1$ of ${\bf d}^m$ exceeds $m-1$,
\begin{eqnarray}
d_1\geq m\;.\label{nomi20}
\end{eqnarray} 
\end{theorem}
{\sf Proof} $\;\;\;$We will prove the Theorem in several steps. First, 
consider the distribution of the generators $d_j,j=3,\ldots,m$ inside the 
matrix representation $M\left\{\Delta(d_1,d_2)\right\}$ in Figure \ref{repr1}. 
According to the definition (\ref{defin2}) of the minimal generating set 
$\{d_1,\ldots,d_m\}$ we have $d_j\in \Delta({\bf d}^{2}),j=3,\ldots,m$. Let $d_u,d_w,3\leq u<w\leq m$ be two of 
such generators with the corresponding representations (\ref{sylves2})
\begin{eqnarray}
\left\{\begin{array}{l}
d_u=d_1 d_2-p_ud_1-q_ud_2\\
d_w=d_1 d_2-p_wd_1-q_wd_2\end{array}\right.,\;\;\;\;\;
1\leq p_u,p_w\leq \left\lfloor d_2-\frac{d_2}{d_1}\right\rfloor\;,\;\;\;
1\leq q_u,q_w\leq d_1-1\;.\label{nomi20a}
\end{eqnarray}
One can show that the minimality of the set $\{d_1,\ldots,d_m\}$ does not 
allow to have at least one of the equalities, $p_u=p_w$ or $q_u=q_w$. 
Indeed, assume, by way of 
contradiction, that the first equality holds, $p_u=p_w$. Then due to 
(\ref{nomi20a}) we have 
$$
d_u=d_w+(q_w-q_u)d_2\;,
$$
which leads to the linear dependence of the three elements $d_2,d_u,d_w$ 
that contradicts (\ref{defin2}). The other equality, $q_u=q_w$, is also 
forbidden for the same reason. Thus, we come to the conclusion that the 
number 
$d_1-1$ of columns in the diagram in Figure \ref{repr1} is at least not less 
than the number $m-2$ of such elements, i.e. 
$$
d_1-1\geq m-2\;,\;\;\;\mbox{or}\;\;\;d_1\geq m-1\;.
$$
Observe that this non--strict inequality was obtained by the assumption 
that every element $d_j,j=3,\ldots,m$ gives rise solely to one associated set 
$\Omega_{d_j}^1({\bf d}^2)$, i.e. $a_{jj}=2,j=3,\ldots,m$. 

Next, in order to prove (\ref{nomi20}) we have to show the existence of at least 
one generator $d_j,3\leq j\leq m$ such that $a_{jj}\geq 3$. Indeed, if $d_h$ 
is such a generator, then $d_h$ gives rise to at least 2 associated sets, 
$\Omega_{d_h}^1({\bf d}^2)$ and $\Omega_{d_h}^2({\bf d}^2)$. Distributing all 
generators $d_j$ inside $M\left\{\Delta(d_1,d_2)\right\}$ we must account for 
$d_h$ twice ($d_h$ and $2d_h$, respectively). The final comparison between the
number of columns in the diagram $M\left\{\Delta(d_1,d_2)\right\}$ and the 
number of distributing generators gives $d_1-1\geq m-1$, or $d_1\geq m$.

Finally, it remains to prove the existence of an integer $d_j,\;3\leq j\leq 
m$, such that $a_{jj}\geq 3$. Consider two last columns, $q=d_1-1$ and 
$q=d_1-2$, of the diagram $M\left\{\Delta(d_1,d_2)\right\}$ and determine the 
numbers $H_{d_1-1}$ and $H_{d_1-2}$ of integers within, respectively. According 
to (\ref{sylves2}) the integers occupying these columns are of the form 
\begin{eqnarray}
\sigma(p,d_1-1)=d_2-pd_1\;\;\;\;\mbox{and}\;\;\;\;
\sigma(p,d_1-2)=2d_2-pd_1\;.\label{nomi20b}
\end{eqnarray}
Therefore, the restriction $\sigma(p,q)>0$ gives
\begin{eqnarray}
\sigma(H_{d_1-1},d_1-1)> 0\;\rightarrow\;H_{d_1-1}=\left\lfloor
\frac{d_2}{d_1}-1\right\rfloor\;\;\mbox{and}\;\;\;
\sigma(H_{d_1-2},d_1-2)> 0\;\rightarrow\;H_{d_1-2}=
\left \lfloor 2\frac{d_2}{d_1}-1\right\rfloor\;.\nonumber
\end{eqnarray}
Denote an integer from the last column $\sigma(p,d_1-1)=d_h$ and prove 
that $2d_h\in\Delta(d_1,d_2)$. In accordance with (\ref{sylves2}) we have
$$
2d_h=2(d_1d_2-pd_1-(d_1-1)d_2)=d_1d_2-pd_1-(d_1-2)d_2=\sigma(2p,d_1-2)\;.
$$
Show that the integer $2d_h$ is contained in the last but one column of the 
diagram $M\left\{\Delta(d_1,d_2)\right\}$. In order to verify this we must 
prove that $H_{d_1-2}\geq 2H_{d_1-1}$.
$$
H_{d_1-2}-2H_{d_1-1}=\left \lfloor 2\frac{d_2}{d_1}-1\right\rfloor-
2\left\lfloor\frac{d_2}{d_1}-1\right\rfloor\geq 2\frac{d_2}{d_1}-2-
2\left\lfloor\frac{d_2}{d_1}-1\right\rfloor=2\left\{\frac{d_2}{d_1}-1
\right\}\geq  0\;.
$$
Thus, $2d_h\in\Delta(d_1,d_2)$. Applying Lemma \ref{lem5} we obtain $2d_3\not
\in\Omega_{d_3}^1({\bf d}^2)$, and therefore $a_{jj}\geq 3$. In fact, we have 
proved a stronger statement, namely, that all integers $d_h$ which belong to 
the last column of the diagram $M\left\{\Delta(d_1,d_2)\right\}$ give rise to 
at least 2 associated sets, and therefore they all have $a_{jj}\geq 3$. 
This completes the proof of the Theorem. $\;\;\;\;\;\;\Box$

It is easy to see that Theorem \ref{theo11} generalizes the restriction 
(\ref{assump3}) obtained for 3D  non--symmetric semigroup. Combining now 
Theorem \ref{theo10} and Theorem \ref{theo11} we can find the minimal bound 
for non--symmetric  semigroup ${\sf S}\left(4,d_2,d_3,d_4\right)$, namely, 
$\#\left\{Q(4,d_2,d_3,d_4;z)\right\}\leq 26$. 

Below we present the results of numerical calculations for two tetrads, 
(4,21,26,43) and (4,31,37,50), which give rise to the corresponding 
$\Delta({\bf d}^4)$--sets and Hilbert series $H({\bf d}^4;z)$. 
\begin{example} $\{d_1,d_2,d_3,d_4\}=\{4,21,26,43\}$
\label{ex4}
{\footnotesize
\begin{eqnarray}
\left\{\begin{array}{rcr}
a_{11}d_1&=&52\\
a_{22}d_2&=&42\\
a_{33}d_3&=&52\\
a_{44}d_4&=&86\end{array}\right.,\;\;
\widehat {\cal A}_4&=&\left(\begin{array}{rrrrr}
13 &0&-2&0\\-4&2&-1&0\\-13&0&2&0\\-11&-2&0&2\end{array}\right),\nonumber\\
\Delta(4,21,26,43)&=&\{1,2,3,5,6,7,9,10,11,13,14,15,17,18,19,22,23,
27,31,35,39\}\;,\nonumber\\
H(4,21,26,43;z)&=&\frac{Q(4,21,26,43;z)}{(1-z^{4})(1-z^{21})(1-z^{26})
(1-z^{43})}\;,\;\;\;\;G(4,21,26,43)=21\;,\nonumber\\
Q(4,21,26,43;z)&=&1-z^{42}-z^{47}-z^{52}-z^{64}+z^{68}-z^{69}+z^{73}+z^{85}-
z^{86}+2z^{90}+z^{95}+\nonumber\\
&&z^{107}-z^{111}+z^{112}-z^{116}-z^{133}\;,\;\;\;\;\;\;\;
\#\left\{Q(4,21,26,43;z)\right\}=18\;.\nonumber
\end{eqnarray}}
\end{example}
\begin{example}$\{d_1,d_2,d_3,d_4\}=\{4,31,37,50\}$
\label{ex5}
{\footnotesize
\begin{eqnarray}
\left\{\begin{array}{rcr}
a_{11}d_1&=&68\\
a_{22}d_2&=&62\\
a_{33}d_3&=&74\\
a_{44}d_4&=&100\end{array}\right.,\;\;\widehat{\cal A}_4&=&\left(\begin{array}
{rrrrr}17&-1&-1&0\\-3&2&0&-1\\-6&0&2&-1\\-8&-1&-1&2\end{array}\right),
\nonumber\\
\Delta(4,31,37,50)&=&\{1,2,3,5,6,7,9,10,11,13,14,15,17,18,19,21,22,
23,25,26,27,29,30,33,34,38,42,46\},\nonumber\\
H(4,31,37,50;z)&=&\frac{Q(4,31,37,50;z)}{(1-z^{4})(1-z^{31})(1-z^{37})
(1-z^{50})}\;,\;\;\;\;G(4,31,37,50)=28\;,\nonumber\\
Q(4,31,37,50;z)&=&1 -z^{62}-z^{68}-z^{74}-z^{81}-z^{87}+z^{99}-z^{100}+
z^{105}+z^{112}+2z^{118}+z^{124}+z^{131}+\nonumber\\
&&z^{137}-z^{149}-z^{155}-z^{168}\;,\;\;\;\;\;\;\;\;\;\;\;
\#\left\{Q(4,31,37,50;z)\right\}=18\;.\nonumber
\end{eqnarray}}
\end{example}
Note that the matrix $\widehat {\cal A}_4$ of the 1st minimal relation ${\cal R}
_1\left({\bf d}^4\right)$ in Example \ref{ex5} is not unique. Indeed, $2d_4=8d_1
+d_2+d_3$ and $2d_4=25d_1$. Nevertheless, this does not affect the final result 
for the Hilbert series $H(4,31,37,50;z)$. 

Both numerators, $Q(4,21,26,35;z)$ and $Q(4,31,37,50;z)$, have exactly 18 
terms satisfying the above restriction (\ref{nomi19e}), $\#\left\{Q(4,d_2,d_3,
d_4;z)\right\}\leq 26$. On the other hand, this may indicate
\footnote{A number 18 appears for $\#\left\{Q(4,d_2,d_3,d_4;z)\right\}$ in  
numerical calculations for a dozen of tetrads $4,d_2,d_3,d_4$ such that a tuple 
$(d_2,d_3,d_4)$ is built out of three pairwise relatively prime elements and 
the only one of them is an even integer not divisible by 4, e.g., $(4,13,15,18)$, 
$(4,17,23,26)$, $(4,29,31,34)$, $(4,41,42,51)$ {\em etc}. The author thanks 
G. Tchernikov for help with numerical calculations.}
that in the 4D Frobenius problem there exist more strong universal properties 
than the upper bound (\ref{nomi19}) for the number of non--zero coefficients 
in the polynomial $Q({\bf d}^4;z)$. 
\section{Genera of higher orders}
\label{gen0}
The generating function of unrepresentable integers $\Phi\left({\bf d}^m;z
\right)$ is a source of another information about the set $\Delta\left({\bf d}^
m\right)$. We show how $\Phi\left({\bf d}^m;z\right)$ can be used 
by computing the power series
\begin{equation}
g_n\left({\bf d}^m\right)=\sum_{s\;\in\;\Delta\left({\bf d}^m\right)}s^n\;,\;
\;\;g_0\left({\bf d}^m\right)=G\left({\bf d}^m\right)\;.\label{gen1}
\end{equation}
For the first time, the simplest series $g_1\left({\bf d}^2\right)$ was 
calculated
in \cite{jau93}. In this Section we give a regular approach to that problem and 
compute some of $g_n$ for 2D and 3D semigroups based on the results obtained 
in Section \ref{qqq1}. Denoting the derivative $d^n/dz^n=\partial_z^n$ we find
\begin{eqnarray}
\partial^n_z\Phi\left({\bf d}^m;1\right)=g_n\left({\bf d}^m\right)-I_1
g_{n-1}\left({\bf d}^m\right)+I_2g_{n-2}\left({\bf d}^m\right)-\ldots\pm
I_{n-1}g_1\left({\bf d}^m\right)\;,
\end{eqnarray}
where the coefficients $I_k$ appear as symmetric invariants of the set 
of the integers $\{1,2,\ldots,n-1\}$
$$
I_1=\sum_{j=1}^{n-1}j=\frac{n(n-1)}{2}\;,\;\;I_2=\sum_{j_1>j_2=1}^{n-1}j_1j_2\;,
\;\;\ldots\;\;I_{n-2}=I_{n-1}\sum_{j=1}^{n-1}\frac{1}{j}\;,\;\;
I_{n-1}=\prod_{j=1}^{n-1}j=(n-1)!\;
$$
and $\partial_z^n\Phi\left({\bf d}^m;1\right)=\partial_z^n\Phi\left({\bf 
d}^m;z\right)_{|\;z=1}$. Successive calculation of the first three terms gives
\begin{eqnarray}
g_1\left({\bf d}^m\right)=\partial_z\Phi\left({\bf d}^m;1\right),\;\;\;
g_2\left({\bf d}^m\right)=\left(\partial_z^2+\partial_z\right)\Phi\left({\bf d}
^m;1\right),\;\;\;g_3\left({\bf d}^m\right)=\left(\partial_z^3+
3\partial_z^2+\partial_z\right)\Phi\left({\bf d}^m;1\right)\;.
\nonumber
\end{eqnarray}
Below we present the first three genera of higher orders for the 
semigroup ${\sf S}\left({\bf d}^2\right)$
\begin{eqnarray}
g_1\left({\bf d}^2\right)&=&\frac{G\left({\bf d}^2\right)}{6}\left(2d_1 d_2-
d_1-d_2-1\right)\;,\;\;\;\;g_2\left({\bf d}^2\right)=\frac{d_1 d_2}{6}\;
G\left({\bf d}^2\right)F\left({\bf d}^2\right)\;,\label{gen2}\\
g_3\left({\bf d}^2\right)&=&\frac{G\left({\bf d}^2\right)}{60}
\left[(1+d_1d_2)\left(1+d_1^2+d_2^2+6d_1^2d_2^2\right)+(d_1+d_2)   
\left(1+d_1^2+d_2^2-9d_1^2d_2^2\right)\right]\;,\nonumber
\end{eqnarray}
and the first genus $g_1\left({\bf d}^3\right)$ for the non--symmetric 
semigroup ${\sf S}\left({\bf d}^3\right)$
\begin{eqnarray}
g_1\left({\bf d}^3\right)&=&\frac{1}{12}
\left(-1+\prod_{i=1}^3d_i+\sum_{i=1}^3A_id_i^2+\sum_{i>j}^3B_{ij}d_id_j-
\prod_{i=1}^3a_{ii}\sum_{i=1}^3C_jd_j\right)\;,\label{gen3}\\
&&A_i=(a_{ii}-1)(2a_{ii}-1)\;,\;\;
B_{ij}=3(a_{ii}-1)(a_{jj}-1)-a_{ii}a_{jj}\;,\;\;C_j=2a_{jj}-3\;.
\nonumber
\end{eqnarray}  
In Example \ref{ex6} we calculate $g_1\left({\bf d}^3\right)$ for the 
triples presented in Example \ref{ex2}.
\begin{example}
\label{ex6}
{\footnotesize
\begin{eqnarray}
g_1(23,29,44)=9526\;,\;\;\;g_1(137,251,256)=2380976\;,\;\;\;
g_1(1563,2275,2503)=12178811815\;.
\nonumber
\end{eqnarray}}
\end{example}
\section*{Acknowledgement}
The usefull discussions with Arye Juhasz and his help in preparing the paper 
are highly appreciated. The author thanks M. Morales, N. Wormald and J. Ramirez 
Alfonsin for sending their papers and M. Beck for his communication.

\noindent
The research was supported by the Kamea Fellowship.

The part of the paper has been written during my stay at the Isaac Newton 
Institute for Mathematical Sciences and its hospitality is highly appreciated.

\newpage
\appendix
\renewcommand{\theequation}{\thesection\arabic{equation}}
\section{Matrix $\widehat {\cal A}_3^{(n)}$ of the 1st minimal relation with 
high degeneration}
\label{appendix1}
\setcounter{equation}{0}
Consider a non--symmetric semigroup ${\sf S}\left({\bf d}^3\right)$ which 
is minimally generated by a triple $d_1,d_2,d_3$ with the matrix $\widehat 
{\cal A}_3^{(n)}$ of the 1st minimal relation where its diagonal elements 
$a_{ii}$ 
completely coincide, $a_{ii}=a$. For such kind of non--symmetric semigroups 
the expressions for $F\left({\bf d}^3\right)$ and $G\left({\bf d}^3\right)$ 
can be represented in a simple form
\begin{eqnarray}
F\left({\bf 
d}^3\right)=\frac{a-2}{2}D_1+\frac{1}{2}\sqrt{(D_1^2-4D_2)a^2+4D_3}\;,\;\;\;
G\left({\bf d}^3\right)=\frac{1}{2}\left[1+(a-1)D_1-a^3\right]\;,
\label{equalB}
\end{eqnarray}
where $D_i$ denote the basic invariants of symmetric group $S_3$ acting on 
$d_1,d_2,d_3$
$$
D_1=d_1+d_2+d_3\;,\;\;D_2=d_1d_2+d_2d_3+d_3d_1\;,\;\;D_3=d_1d_2d_3\;.
$$
In this Appendix we present all possible different admissible triples $d_1,d_2,
d_3$ generating the non--symmetric semigroups and corresponding to the matrix 
$\widehat {\cal A}_3^{(n)}$ with a complete coincidence of their diagonal 
elements $a_{ii}=a$ for the first three values $a=3,4,5$.
\begin{itemize}
\item a=3
\begin{eqnarray}
{\footnotesize
\left.\begin{array}{r}
5\\7\\8\end{array}\right.\rightarrow\left(\begin{array}{rrr}
 3 &-1 &-1\\
-1 & 3 &-2 \\
-2 &-2 & 3 \end{array}\right).}
\label{equalA3}
\end{eqnarray}
{\footnotesize
$$
F(5,7,8)=11\;,\;\;G(5,7,8)=7\;.\nonumber
$$}
\item a=4
\begin{eqnarray}
{\footnotesize
\left.\begin{array}{r}
7\\13\\15\end{array}\right.\rightarrow\left(\begin{array}{rrr}
 4 &-1 &-1\\
-1 & 4 &-3 \\
-3 &-3 & 4 \end{array}\right),\;\left.\begin{array}{r}
10\\13\\14\end{array}\right.\rightarrow\left(\begin{array}{rrr}
 4 &-2 &-1\\ 
-1 & 4 &-3\\
-3 &-2 & 4 \end{array}\right)\;.}
\label{equalA4} 
\end{eqnarray}
{\footnotesize
$$
F(7,13,15)=38\;,\;\;G(7,13,15)=21\;,\;\;\;\;\;F(10,13,14)=45\;,\;\;
G(10,13,14)=24\;.
$$}
\item a=5
\begin{eqnarray}
{\footnotesize
\left.\begin{array}{r}
9\\22\\23\end{array}\right.\rightarrow\left(\begin{array}{rrr}
 5 &-1 &-1\\
-2 & 5 &-4\\
-3&-4 & 5 \end{array}\right),\;\left.\begin{array}{r}
16\\17\\23 \end{array}\right.\rightarrow\left(\begin{array}{rrr}
 5 &-2 &-2\\
-1 & 5 &-3\\
-4 &-3 & 5 \end{array}\right),\;\left.\begin{array}{r}
17\\19\\22\end{array}\right.\rightarrow\left(\begin{array}{rrr}
 5 &-1 &-3\\
-3 & 5 &-2\\
-2 &-4 & 5 \end{array}\right),}\nonumber\\
{\footnotesize
\left.\begin{array}{r}
13\\19\\23 \end{array}\right.\rightarrow\left(\begin{array}{rrr}
 5 &-1 &-2\\
-2 & 5 &-3 \\
-3 &-4 & 5 \end{array}\right),\;\left.\begin{array}{r}
13\\21\\23\end{array}\right.\rightarrow\left(\begin{array}{rrr} 
 5 &-2 &-1\\
-1 & 5 &-4\\
-4 &-3 & 5 \end{array}\right),\;\left.\begin{array}{r}
13\\21\\22\end{array}\right.\rightarrow\left(\begin{array}{rrr}
 5 &-1 &-2\\
-3 & 5 &-3\\ 
-2 &-4 & 5 \end{array}\right),}\label{equalA5}\\
{\footnotesize
\left.\begin{array}{r}
13\\17\\24\end{array}\right.\rightarrow\left(\begin{array}{rrr}
 5 &-1 &-2\\
-1 & 5 &-3\\
-4 &-4 & 5 \end{array}\right),\;\left.\begin{array}{r}
16\\19\\21\end{array}\right.\rightarrow\left(\begin{array}{rrr}
 5 &-2 &-2\\
-2 & 5 &-3\\
-3 &-3 & 5\end{array}\right),\;\left.\begin{array}{r}
17\\21\\22\end{array}\right.\rightarrow\left(\begin{array}{rrr}
 5 &-3 &-1\\
-1 & 5 &-4\\
-4 &-2 & 5 \end{array}\right).}\nonumber
\end{eqnarray}
{\footnotesize
\begin{eqnarray}
\left.\begin{array}{l}
F(9,22,23)=83\\G(9,22,23)=46\end{array}\right.,
\left.\begin{array}{l}
F(16,17,23)=93\\G(16,17,23)=50\end{array}\right.,
\left.\begin{array}{l}
F(17,19,22)=103\\G(17,19,22)=54\end{array}\right.,
\left.\begin{array}{l}
F(13,19,23)=86\\G(13,19,23)=48\end{array}\right.,
\left.\begin{array}{l}
F(13,21,23)=100\\G(13,21,23)=52\end{array}\right.\nonumber\\
\left.\begin{array}{l}
F(13,21,22)=93\\G(13,21,22)=50\end{array}\right.,
\left.\begin{array}{l}
F(13,17,24)=83\\G(13,17,24)=46\end{array}\right.,
\left.\begin{array}{l}
F(16,19,21)=87\\G(16,19,21)=50\end{array}\right.,
\left.\begin{array}{l}
F(17,21,22)=113\\G(17,21,22)=58\end{array}\right.\nonumber
\end{eqnarray}}
\end{itemize}
In fact, there is one more, the 10th tuple, $(9,21,24)$, generating a 
semigroup ${\sf S}\left({\bf d}^3\right)$ with the matrix $\widehat{\cal 
A}_3^{(n)}$ of the 1st minimal relation which is distinguished from those 
presented in (\ref{equalA5}). However, the set $\{9,21,24\}$ is not minimal 
since $\gcd(9,21,24)=3$ and therefore is not included into (\ref{equalA5}).

Note that there exist only two pairs of triples -- (9,22,23), (13,17,24) and 
(16,17,23), (13,21,22) -- which have at the same time the equal Frobenius 
numbers and genera in every pair. This fact may be interesting in the sense of 
the question posed in \cite{frob87} on the number of semigroups ${\sf S}\left(
{\bf d}^3\right)$ with the prescribed Frobenius number $F\left({\bf d}^3\right)=
const$. Here we have two constraints, $F\left({\bf d}^3\right)=const_1$ and 
$Q\left({\bf d}^3\right)=const_2$, that must essentially diminish the number of 
addimisible semigroups ${\sf S}\left({\bf d}^3\right)$.
\renewcommand{\theequation}{\thesection\arabic{equation}}
\section{On two conjectures about the upper bound for $F({\bf d}^3)$}
\label{appendix2}
\setcounter{equation}{0}
A recent paper \cite{beck05} asserts two conjectures based on numerical 
calculations for more than ten thousands randomly chosen admissible triples 
$(d_1,d_2,d_3)$ such
\footnote{In fact, the typical values of $d_i$ were not exceeding 750 
\cite{beck04}.}
that $\sqrt{d_1d_2d_3}<2\cdot 10^4$. We quote from \cite{beck05}:

`{\em For all admissible triples $(d_1,d_2,d_3)$ the Frobenius number 
$F({\bf d}^3)$ can be bounded from above},
\begin{equation}
F({\bf d}^3)\leq F^+_{C;\nu}({\bf d}^3)\;,\;\;\; 
F^+_{C;\nu}({\bf d}^3)=C(d_1d_2d_3)^{\nu}-(d_1+d_2+d_3)\;.
\label{zak1}
\end{equation}

{\em where $C=const$ and $\nu<2/3$}`,

and further,

`{\em In fact, our data suggests, more precisely, that for all admissible 
triples $(d_1,d_2,d_3)$,}
\begin{equation}
F({\bf d}^3)\leq F^+_{1;5/8}({\bf d}^3)\;,\;\;\mbox{i.e.}\;\;
F^+_{1;5/8}({\bf d}^3)=(d_1d_2d_3)^{5/8}-(d_1+d_2+d_3)\;.
\label{zak2}
\end{equation}

\noindent
In this Appendix we are going to falsify these both conjectures. 

We start with (\ref{zak2}) by showing two counterexamples. Following 
\cite{beck05} recall the terms which are necessary to discuss this conjecture. 
First, call the triple $(d_1,d_2,d_3)$ {\it constituting an almost arithmetic 
sequence} if there exist the integers $a,b$ such that
\begin{eqnarray}
d_2=a d_1+b,\;\;d_3=a d_1+2b,\;\;\;a\geq 1\;,\;\;b\geq 1,\;\;\gcd(d_1,b)=1\;.
\label{adddm0}
\end{eqnarray}
Next, call the triple $(d_1,d_2,d_3)$ {\it excluded} if at least one   
of the following holds:
\begin{eqnarray}
&&\mbox{1) one of the elements}\; d_i\; \mbox{being representable by the
other two;}\label{adddm1}\\
&&\mbox{2) one element}\;d_i\;\mbox{dividing the sum of the other two;}
\label{adddm2}\\
&&\mbox{3) the elements}\; d_i\; \mbox{represent an almost arithmetic 
sequence}. \label{adddm3}
\end{eqnarray}
Following \cite{beck05} define an admissible triple $(d_1,d_2,d_3)$ as {\em 
a triple of pairwise coprime integers that is not excluded}. 

Now consider the triple of pairwise coprime integers $d_1,d_2,d_3$ generating 
a non--symmetric semigroup  ${\sf S}\left({\bf d}^3\right)$ with the matrix 
$\widehat {\cal A}_{3}$ (see (\ref{herznon2}))
\begin{eqnarray}
\left\{\begin{array}{rcr}
d_1&=&10001=73\cdot 137\\
d_2&=&10003=7\cdot 1429\\
d_3&=&20003=83\cdot 241\end{array}\right.,\;\;\;\;
\widehat {\cal A}_3=\left(\begin{array}{rrrr}
5003&-5000&-1\\-5000&5001&-1\\-3&-1&2\end{array}\right)\;,\label{new1}
\end{eqnarray}
where $d_i$ are uniquely factorized into a product of primes. Notice that 
\begin{eqnarray}
2d_3=3d_1+d_2\;,\;\;\;\;d_2-d_1\ll d_1\;,\label{new1o}
\end{eqnarray}
and 
\begin{eqnarray}
d_1>2^{13}\;.\label{new1ko}
\end{eqnarray}
Show that the triple (\ref{new1}) is admissible. 
First, (\ref{adddm1}) is not satisfied due to minimality of the set $\{10001,
10003,20003\}$ according to the matrix $\widehat {\cal A}_3$ of minimal 
relation in (\ref{new1}) (see (\ref{defin2}) and (\ref{herznon3})). Next, 
(\ref{adddm2}) is not satisfied, since 
\begin{eqnarray}
(10001+10003)/20003\not\in {\mathbb N}\;,\;
(10003+20003)/10001\not\in {\mathbb N}\;,\;
(20003+10001)/10003\not\in {\mathbb N}\;.
\label{new1b}
\end{eqnarray}
In order to prove that (\ref{adddm3}) is also not satisfied observe that it is 
sufficient to show, in accordance with (\ref{adddm0}), that $a=(2d_2-d_3)/d_1$ 
is not an integer. Indeed, a straightforward calculation gives $a=
\frac{3}{10001}$. Thus, the triple (\ref{new1}) is not excluded and according 
to definition \cite{beck05} is admissible.

Calculate the Frobenius number for the triple (\ref{new1}). By (\ref{ari13a}) 
we obtain 
$$
F(10001,10003,20003)=50014999\;,
$$
while the conjectured bound (\ref{zak2}) reads
$$
F^+_{1;5/8}({\bf d}^3)=(10001\cdot 10003\cdot 20003)^{5/8}-
(10001+10003+20003)=48745742.422\;.
$$
Thus, the conjecture (\ref{zak2}) is disproved. The contradiction becomes 
even stronger if we increase the values of $d_1,d_2,d_3$ preserving 
(\ref{new1o}) and (\ref{new1ko}), e.g.
\begin{eqnarray}
&&F(100001,100003,200003)=5000149999\;,\;\;\;
F^+_{1;5/8}(100001,100003,200003)=3656883908.3\;,\nonumber\\
&&\;\;\;\mbox{where}\;\;\;\;\;100001=11\cdot 9091\;,\;\;100003=1\cdot 
100003\;,\;\;200003=1\cdot 200003\;.\nonumber
\end{eqnarray}
The relation between the degrees "5/8" in (\ref{zak2}) and "13" in 
(\ref{new1ko}) is not accidental and will be clarified below. 

Move on to the main conjecture (\ref{zak1}) and consider the triple 
$(d_1,d_2,d_3)$ such that $d_1$ is a prime number and 
\begin{eqnarray}
d_1=2l+1\;,\;\;\;d_2=d_1+2=2l+3\;,\;\;\;d_3=2d_1+1=4l+3\;,\;\;\;
l\gg 1\;.\label{new2}
\end{eqnarray}
Note that (\ref{new1o}) is satisfied and all elements $d_1,d_1+2,2d_1+1$ of 
the triple are pairwise coprime integers. The minimal set $\{2l+1,2l+3,4l+3\}$ 
generates a non--symmetric semigroup ${\sf S}\left({\bf d}^3\right)$ with the 
matrix $\widehat {\cal A}_{3}$ of the 1st minimal relation
\begin{eqnarray}
\widehat {\cal A}_3=\left(\begin{array}{ccc}
l+3&-l&-1\\-l&l+1&-1\\-3&-1&2\end{array}\right)\;,\label{new2a}
\end{eqnarray}
and the Frobenius number (see (\ref{ari13a}))
\begin{eqnarray}
F(2l+1,2l+3,4l+3)=2l^2+3l-1\;.\label{new2b}
\end{eqnarray}
We prove an auxiliary Lemma.
\begin{lemma}
\label{lem17}
The triple (\ref{new2}) is admissible.
\end{lemma}
{\sf Proof} $\;\;\;$
First, (\ref{adddm1}) is not satisfied due to minimality of the set $\{2l+1,2l+3,4l+3\}$ according 
to (\ref{new2a}). Next, (\ref{adddm2}) is not satisfied, since
\begin{eqnarray}
\frac{4l+4}{4l+3}\not\in {\mathbb N}\;,\;\;\;
\frac{6l+4}{2l+3}\not\in {\mathbb N}\;,\;\;\;
\frac{6l+6}{2l+1}\not\in {\mathbb N}\;.\label{new2c}
\end{eqnarray}
Finally, calculating $a=(2d_2-d_3)/d_1$ we get $a=3/(2l+1)$ and conclude that 
(\ref{adddm3}) is also not satisfied. Thus, the triple (\ref{new2}) is 
admissible and the Lemma is proved.$\;\;\;\;\;\;\Box$

Finally we are ready to prove the main Lemma of this Appendix.
\begin{lemma}
\label{lem18}
Let ${\bf d}^3$ be given admissible triple, ${\bf d}^3=(d_1,d_2,d_3)$. The 
Frobenius number $F({\bf d}^3)$ can not be bounded from above by $F^+_{C;\nu}
({\bf d}^3)$ given by 
\begin{equation}
F^+_{C;\nu}({\bf d}^3)=C(d_1d_2d_3)^{\nu}-(d_1+d_2+d_3)\;,\;\;\;\nu<2/3\;,
\label{zak3} 
\end{equation}
for any $C=const$.
\end{lemma}
{\sf Proof} $\;\;\;$Consider the triple (\ref{new2}) which is admissible 
according to Lemma \ref{lem17} and denote by $\delta_{C;\nu}(l)$ the ratio
\begin{eqnarray}
\delta_{C;\nu}(l)=\frac{F^+_{C;\nu}(2l+1,2l+3,4l+3)}{F(2l+1,2l+3,4l+3)}\;.
\label{new2d}
\end{eqnarray}
In order to verify the conjecture \ref{zak1} we have to find $C=const$ and
$\nu<2/3$ such that
\begin{eqnarray}
\delta_{C;\nu}(l)\geq 1\label{new2f}
\end{eqnarray}
holding for all $l>1$. However, this is not true. Indeed, find a leading 
term of the asymptotics of $\delta_{C;\nu}(l)$ when $l\to \infty$  
\begin{eqnarray}
\delta_{C;\nu}(l)\simeq C2^{4\nu-1}l^{3\nu-2}
\;.\label{new2e}
\end{eqnarray}
Observe that its growth with $l\to \infty$ is enough to break (\ref{new2f}) 
when $l$ exceeds a critical value $l_{cr}$
\begin{eqnarray}
l>l_{cr}\;,\;\;\;
\lg_2 l_{cr}=\frac{4\nu-1}{2-3\nu}+\frac{\lg_2 C}{2-3\nu}
\label{new2k}
\end{eqnarray}
for all $\nu <2/3$. This is true for arbitrary large finite $C$.
$\;\;\;\;\;\;\Box$

Hence there follows the critical value $l_{cr}=2^{12}$ for $\nu=5/8$ and 
$C=1$ that leads to (\ref{new1ko}). 

\end{document}